\documentclass{article}
\usepackage{amsmath, amssymb, amsthm, graphicx, cite, comment}
\usepackage[arrow,curve,matrix,arc,2cell]{xy}
\UseAllTwocells
\DeclareFontFamily{U}{rsfs}{} \DeclareFontShape{U}{rsfs}{n}{it}{<->
rsfs10}{} \DeclareSymbolFont{mscr}{U}{rsfs}{n}{it}
\DeclareSymbolFontAlphabet{\scr}{mscr}
\def\mathscr{\scr}
\begin{document}
\def\e#1\e{\begin{equation}#1\end{equation}}
\def\ea#1\ea{\begin{align}#1\end{align}}
\def\eq#1{{\rm(\ref{#1})}}
\theoremstyle{plain}
\newtheorem{thm}{Theorem}[section]
\newtheorem{prop}[thm]{Proposition}
\newtheorem{lem}[thm]{Lemma}
\newtheorem{cor}[thm]{Corollary}
\newtheorem{quest}[thm]{Question}
\newtheorem{conj}[thm]{Conjecture}
\newtheorem{princ}[thm]{Principle}
\newtheorem{prob}[thm]{Problem}
\theoremstyle{definition}
\newtheorem{dfn}[thm]{Definition}
\newtheorem{ex}[thm]{Example}
\newtheorem{rem}[thm]{Remark}
\numberwithin{equation}{section}
\numberwithin{figure}{section}
\def\dim{\mathop{\rm dim}\nolimits}
\def\id{\mathop{\rm id}\nolimits}
\def\SO{\mathop{\rm SO}\nolimits}
\def\Hom{\mathop{\rm Hom}\nolimits}
\def\Hol{\mathop{\rm Hol}\nolimits}
\def\Re{\mathop{\rm Re}}
\def\Im{\mathop{\rm Im}}
\def\SU{\mathop{\rm SU}}
\def\ind{\mathop{\rm ind}}
\def\area{\mathop{\rm area}}
\def\U{{\rm U}}
\def\vol{\mathop{\rm vol}\nolimits}
\def\coh{\mathop{\rm coh}}
\def\nov{{\rm nov}}
\def\bs{\boldsymbol}
\def\ge{\geqslant}
\def\le{\leqslant\nobreak}
\def\O{{\mathbin{\cal O}}}
\def\cA{{\mathbin{\cal A}}}
\def\oA{{\mathbin{\smash{\,\overline{\!\mathcal A}}}}}
\def\cB{{\mathbin{\cal B}}}
\def\cC{{\mathbin{\cal C}}}
\def\cD{{\mathbin{\scr D}}}
\def\cE{{\mathbin{\cal E}}}
\def\cF{{\mathbin{\cal F}}}
\def\cG{{\mathbin{\cal G}}}
\def\cH{{\mathbin{\cal H}}}
\def\cL{{\mathbin{\cal L}}}
\def\cM{{\mathbin{\cal M}}}
\def\oM{{\mathbin{\smash{\,\,\overline{\!\!\mathcal M\!}\,}}}}
\def\cO{{\mathbin{\cal O}}}
\def\cP{{\mathbin{\cal P}}}
\def\cS{{\mathbin{\cal S}}}
\def\cT{{\mathbin{\cal T}}}
\def\PV{{\mathbin{\cal{PV}}}}
\def\TS{{\mathbin{\cal{TS}}}}
\def\cQ{{\mathbin{\cal Q}}}
\def\cW{{\mathbin{\cal W}}}
\def\C{{\mathbin{\mathbb C}}}
\def\CP{{\mathbin{\mathbb{CP}}}}
\def\F{{\mathbin{\mathbb F}}}
\def\Q{{\mathbin{\mathbb Q}}}
\def\R{{\mathbin{\mathbb R}}}
\def\Z{{\mathbin{\mathbb Z}}}
\def\sF{{\mathbin{\mathscr F}}}
\def\al{\alpha}
\def\be{\beta}
\def\ga{\gamma}
\def\de{\delta}
\def\io{\iota}
\def\ep{\epsilon}
\def\la{\lambda}
\def\ka{\kappa}
\def\th{\theta}
\def\ze{\zeta}
\def\up{\upsilon}
\def\vp{\varphi}
\def\si{\sigma}
\def\om{\omega}
\def\De{\Delta}
\def\La{\Lambda}
\def\Si{\Sigma}
\def\Tau{{\rm T}}
\def\Th{\Theta}
\def\Om{\Omega}
\def\Ga{\Gamma}
\def\Up{\Upsilon}
\def\sSi{{\smash{\sst\Si}}}
\def\pd{\partial}
\def\ts{\textstyle}
\def\st{\scriptstyle}
\def\sst{\scriptscriptstyle}
\def\w{\wedge}
\def\sm{\setminus}
\def\bu{\bullet}
\def\op{\oplus}
\def\ot{\otimes}
\def\ov{\overline}
\def\ul{\underline}
\def\bigop{\bigoplus}
\def\bigot{\bigotimes}
\def\iy{\infty}
\def\es{\emptyset}
\def\ra{\rightarrow}
\def\Ra{\Rightarrow}
\def\Longra{\Longrightarrow}
\def\ab{\allowbreak}
\def\longra{\longrightarrow}
\def\hookra{\hookrightarrow}
\def\dashra{\dashrightarrow}
\def\t{\times}
\def\ci{\circ}
\def\ti{\tilde}
\def\d{{\rm d}}
\def\ha{{\ts\frac{1}{2}}}
\def\md#1{\vert #1 \vert}
\def\bmd#1{\big\vert #1 \big\vert}
\def\ms#1{\vert #1 \vert^2}
\def\nm#1{\Vert #1 \Vert}
\title{Conjectures on Bridgeland stability for Fukaya categories of Calabi--Yau manifolds, special Lagrangians, and Lagrangian mean curvature flow}
\author{Dominic Joyce}
\date{}
\maketitle

\begin{abstract} Let $(M,J,g,\Om)$ be a Calabi--Yau $m$-fold, and consider compact, graded Lagrangians $L$ in $M$. Thomas and Yau \cite{Thom,ThYa} conjectured that there should be a notion of `stability' for such $L$, and that if $L$ is stable then Lagrangian mean curvature flow $\{L^t:t\in[0,\iy)\}$ with $L^0=L$ should exist for all time, and $L^\iy=\lim_{t\ra\iy}L^t$ should be the unique special Lagrangian in the Hamiltonian isotopy class of $L$. This paper is an attempt to update the Thomas--Yau conjectures, and discuss related issues. 

It is a folklore conjecture, extending \cite{Thom}, that there exists a Bridgeland stability condition $(Z,\cP)$ on the derived Fukaya category $D^b\sF(M)$, such that an isomorphism class in $D^b\sF(M)$ is $(Z,\cP)$-semistable if (and possibly only if) it contains a special Lagrangian, which must then be unique.

In brief, we conjecture that if $(L,E,b)$ is an object in an enlarged version of $D^b\sF(M)$, then there is a unique family $\{(L^t,E^t,b^t):t\in[0,\iy)\}$ such that $(L^0,E^0,b^0)=(L,E,b)$, and $(L^t,E^t,b^t)\cong(L,E,b)$ in $D^b\sF(M)$ for all $t$, and $\{L^t:t\in[0,\iy)\}$ satisfies Lagrangian MCF {\it with surgeries\/} at singular times $T_1,T_2,\ldots,$ and in graded Lagrangian integral currents we have $\lim_{t\ra\iy}L^t=L_1+\cdots+L_n$, where $L_j$ is a special Lagrangian integral current of phase $e^{i\pi\phi_j}$ for $\phi_1>\cdots>\phi_n$, and $(L_1,\phi_1),\ldots,(L_n,\phi_n)$ correspond to the decomposition of $(L,E,b)$ into $(Z,\cP)$-semistable objects. 

We also give detailed conjectures on the nature of the singularities of Lagrangian MCF that occur at the finite singular times $T_1,T_2,\ldots.$
\end{abstract}

\setcounter{tocdepth}{2}
\tableofcontents

\section{Introduction}
\label{bs1}

Thomas \cite{Thom} and Thomas and Yau \cite{ThYa} proposed some interesting conjectures on graded Lagrangians $L$ in Calabi--Yau manifolds $(M,J,g,\Om)$: they defined a notion of `stability' \cite[Def.~5.1]{Thom} for Hamiltonian isotopy classes $[L]$ of (almost calibrated) graded Lagrangians $L$, and conjectured \cite[Conj.~5.2]{Thom} that $[L]$ contains a (unique) special Lagrangian $L'$ if and only if $[L]$ is stable. Furthermore, they conjectured \cite[Conj.~7.3]{ThYa} that if $[L]$ is stable and $L$ satisfies an extra condition \cite[(7.1) or (7.2)]{ThYa} then Lagrangian mean curvature flow $\{L^t:t\in[0,\iy)\}$ with $L^0=L$ exists for all time, and~$\lim_{t\ra\iy}L^t=L'$.

Thomas and Yau's papers \cite{Thom,ThYa} are remarkably prescient, as they predate (and motivated) much important mathematics relevant to their picture, including the invention of Bridgeland stability on triangulated categories \cite{Brid1}, the publication of Fukaya, Oh, Ohta and Ono's \cite{Fuka1,Fuka2,FOOO} and Seidel's work \cite{Seid3} on Lagrangian Floer cohomology and Fukaya categories, and progress on singularities of Lagrangian MCF such as Neves \cite{Neve1,Neve2,Neve3}. I believe their big picture is correct, although I think they are too optimistic in expecting Lagrangian MCF to exist for all time without singularities even in the stable case, and want to substitute Lagrangian MCF {\it with surgeries\/} instead (see (iv) below).

The aim of this paper is to update the Thomas--Yau conjectures in the light of subsequent discoveries, to add more detail to the picture, and to extend their scope. Thomas and Yau's papers are clearly intended as a programme for future research rather than as precise conjectures; they have many caveats on points they are uncertain about, and the conjectures they actually state are fairly cautious (for instance, the inclusion of the strong condition \cite[(7.1) or (7.2)]{ThYa} on stable Lagrangians for Lagrangian MCF to converge to a special Lagrangian).

I am going to be a lot less cautious, and will make conjectures on unique long-time existence of Lagrangian MCF with surgeries starting from {\it any\/} compact graded Lagrangian with unobstructed Lagrangian Floer cohomology in the sense of \cite{FOOO} (`$HF^*$ unobstructed' for short). Nonetheless, I ask readers to take the conjectures in the spirit they are intended: as provisional, quite probably false in their current form, to be refined (or discarded) as our understanding improves, but in the mean time, as (hopefully) a useful guide and motivation for research in the area. I will say more on this in the introduction to~\S\ref{bs3}.

In reading Thomas and Yau \cite{Thom,ThYa}, I think it is helpful to impose the standing assumption that all graded Lagrangians $L$ considered are {\it almost calibrated}, that is, have phase variation less than $\pi$. This is not clearly articulated in \cite{Thom,ThYa}, although bounds on the phase variation are assumed in several places, with the almost calibrated condition used in \cite[\S 5.3]{ThYa}. We need $L$ to be almost calibrated since otherwise the `global phase' $\phi(L)$ in \cite[\S 3]{Thom} is not well-defined, and so `stability' in \cite[Def.~5.1]{Thom} does not make sense. 

Including the almost calibrated assumption, I am not aware of any counterexamples to the precise conjectures stated in \cite{Thom,ThYa} (although I do expect such counterexamples to exist, see (iv) below). In particular, Neves' examples \cite{Neve3} of finite time singularities to Lagrangian MCF discussed in Example \ref{bs3ex5} below are not almost calibrated, and so not counterexamples to~\cite[Conj.~7.3]{ThYa}.

Here are the main differences between our programme and that of~\cite{Thom,ThYa}:
\begin{itemize}
\setlength{\parsep}{0pt}
\setlength{\itemsep}{0pt}
\item[(i)] We work in the `derived Fukaya category' $D^b\sF(M)$ of $M$, as in Fukaya, Oh, Ohta and Ono \cite{Fuka1,Fuka2,FOOO} (see also Seidel \cite{Seid3}). Objects of $D^b\sF(M)$ include triples $(L,E,b)$, where $L$ is a compact, graded Lagrangian in $M$ and $E\ra L$ a rank one local system such that $(L,E)$ has `$HF^*$ unobstructed', and $b$ is a `bounding cochain' for $(L,E)$, as in~\cite{FOOO}. 

Rather than working in a Hamiltonian isotopy class $[L]$ as in \cite{Thom,ThYa}, we work in an isomorphism class $[(L,E,b)]$ in~$D^b\sF(M)$.
\item[(ii)] The derived Fukaya category $D^b\sF(M)$ must be enlarged to include immersed Lagrangians as in \cite{AkJo} in dimension $m\ge 2$, and certain classes of singular Lagrangians in dimension $m\ge 3$, for the programme to work.
\item[(iii)] Our notion of `stability' of Lagrangians is a `Bridgeland stability condition' $(Z,\cP)$ on the triangulated category $D^b\sF(M)$, as in Bridgeland~\cite{Brid1}.
\item[(iv)] Even for a `stable' object $(L,E,b)$ with small phase variation, I do not expect Lagrangian MCF $\{L^t:t\in[0,\iy)\}$ with $L^0=L$ to exist without singularities, as hoped in \cite{ThYa}. Instead, in a similar way to the proof of the Poincar\'e Conjecture by Perelman and others using Ricci flow (see \cite{MoTi,Pere1,Pere2,Pere3}), I expect there to exist a unique family $\{(L^t,E^t,b^t):t\in[0,\iy)\}$ of objects in the isomorphism class of $(L,E,b)$ in $D^b\sF(M)$ with $(L^0,E^0,b^0)=(L,E,b)$, where $\{L^t:t\in[0,\iy)\}$ satisfies Lagrangian MCF {\it with surgeries}. 

That is, at a discrete series of `singular times' $t=T_1,T_2,\ldots$ the flow develops a singularity, but one can continue the flow uniquely for $t>T_i$ in a way which is continuous at $t=T_i$ in a weak sense. The $L^t$ for $T_i-\ep<t<T_i$ and for $T_i<t<T_i+\ep$ may have different topologies.
\item[(v)] Lagrangians $L$ or pairs $(L,E)$ in $M$ `with $HF^*$ obstructed' do not give objects of $D^b\sF(M)$, and `stability' does not make sense for them.

For $L$ with $HF^*$ obstructed, the author expects that Lagrangian MCF $\{L^t:t\in[0,T)\}$ with $L^0=L$ may develop finite time singularities at $t=T$ after which it is not possible to continue the flow, even with a surgery. So, the long time existence of Lagrangian MCF with surgeries in (iv) should apply only for Lagrangians with $HF^*$ unobstructed.
\end{itemize}

Part (iv), our insistence on including finite time singularities of Lagrangian MCF and surgeries, is the greatest divergence between our picture and that of \cite{Thom,ThYa}. As some justification, note that Neves \cite{Neve3} proves that every Hamiltonian isotopy class $[L]$ of compact Lagrangians $L$ in a Calabi--Yau $m$-fold $(M,J,g,\Om)$ for $m\ge 2$ contains (not almost calibrated) representatives $\ti L$ such that Lagrangian MCF $\{L^t:t\in[0,T)\}$ with $L^0=\ti L$ develops a finite time singularity at $t=T$, so without (strong) extra assumptions, finite time singularities of Lagrangian MCF are unavoidable.

One of the goals of this paper is to persuade mathematicians working on Lagrangian MCF that obstructions to $HF^*$ are important in understanding finite time singularities of Lagrangian MCF, that the flow should be better behaved if $HF^*$ is unobstructed, and that tools from symplectic topology such as $J$-holomorphic curves, Lagrangian Floer cohomology, and Fukaya categories, should be used to make the next generation of advances in the field.

Some evidence for this is provided by Imagi, Oliveira dos Santos and the author \cite{IJO}, in which, motivated by this paper, we use Lagrangian Floer cohomology and Fukaya categories to prove that the unique special Lagrangians in $\C^m$ asymptotic at infinity to the union $\Pi_0\cup\Pi_{\bs\phi}$ of two transverse Lagrangian planes $\Pi_0,\Pi_{\bs\phi}$ are the `Lawlor necks' of \cite{Lawl}, and the unique Lagrangian MCF expanders in $\C^m$ asymptotic at infinity to $\Pi_0\cup\Pi_{\bs\phi}$ are the examples in Lee, Tsui and the author \cite[Th.s C \& D]{JLT}, as in Theorems \ref{bs2thm3} and \ref{bs2thm6} below.

Section \ref{bs2} explains some background material, \S\ref{bs3} states the conjectures, and \S\ref{bs4} discusses some generalizations.
\smallskip

\noindent{\it Acknowledgements.} I would like to thank
Mohammed Abouzaid, Joana Amorim, Lino Amorim, Denis Auroux, Mark Haskins, Yohsuke Imagi, Yng-Ing Lee, Andr\'e Neves, Paul Seidel, Richard Thomas, and Ivan Smith for useful conversations. This research was supported by EPSRC grant EP/H035303/1.

\section{Background material}
\label{bs2}

We now summarize the background material we will need to state our conjectures in \S\ref{bs3}. We discuss Calabi--Yau $m$-folds, graded Lagrangians and special Lagrangians in \S\ref{bs21} and Lagrangian mean curvature flow in \S\ref{bs23}, giving examples of SL $m$-folds in $\C^m$ in \S\ref{bs22} and solitons for Lagrangian MCF in \S\ref{bs24}. Section \ref{bs25} explains Lagrangian Floer cohomology, obstructions to $HF^*$, and derived Fukaya categories $D^b\sF(M)$ for embedded Lagrangians in Calabi--Yau $m$-folds, and \S\ref{bs26} considers the extension to immersed Lagrangians.

Some references are McDuff and Salamon \cite{McSa} for symplectic geometry, the author \cite{Joyc11} and Harvey and Lawson \cite{HaLa} for Calabi--Yau $m$-folds and special Lagrangians, Mantegazza \cite{Mant}, Smoczyk \cite{Smoc2} and Neves \cite{Neve2} for (Lagrangian) MCF, Fukaya \cite{Fuka1,Fuka2}, Fukaya, Oh, Ohta and Ono \cite{FOOO} and Seidel \cite{Seid3} for Lagrangian Floer cohomology and Fukaya categories for embedded Lagrangians, and Akaho and the author \cite{AkJo} for the extension to immersed Lagrangians.  

\subsection{Calabi--Yau $m$-folds and special Lagrangians}
\label{bs21}

We define Calabi--Yau $m$-folds, graded Lagrangians, and special Lagrangians.

\begin{dfn} A {\it Calabi--Yau $m$-fold\/} is a quadruple
$(M,J,g,\Om)$ such that $(M,J)$ is an $m$-dimensional complex
manifold, $g$ is a K\"ahler metric on $(M,J)$ with K\"ahler form
$\om$, and $\Om$ is a holomorphic $(m,0)$-form on $(M,J)$ satisfying
\e
\om^m/m!=(-1)^{m(m-1)/2}(i/2)^m\Om\w\bar\Om.
\label{bs2eq1}
\e
Then $g$ is Ricci-flat and its holonomy group is a subgroup of
$\mathop{\rm SU}(m)$. We do not require $M$ to be compact, or $g$ to
have holonomy $\mathop{\rm SU}(m)$, although many authors make these
restrictions.

If $(M,J,g,\Om)$ is a Calabi--Yau $m$-fold with K\"ahler form $\om$,
then $(M,\om)$ is a symplectic manifold. A {\it Lagrangian\/} $L$ in
$M$ is a real $m$-dimensional submanifold (embedded or immersed)
with $\om\vert_L=0$.

Let $L$ be a Lagrangian in $M$. Then $\Om\vert_L$ is a complex
$m$-form on $L$. Equation \eq{bs2eq1} implies that
$\bmd{\Om\vert_L}=1$, where $\md{\,.\,}$ is computed using the
Riemannian metric $g\vert_L$. Suppose $L$ is oriented. Then we have
a volume form $\d V_L$ on $L$ defined using the metric $g\vert_L$
and orientation with $\md{\d V_L}=1$, so $\Om\vert_L=\Th_L\cdot \d
V_L$, where $\Th_L:L\ra\U(1)$ is a unique smooth function, and
$\U(1)=\{z\in\C:\md{z}=1\}$.

There is an induced morphism of cohomology groups
$\Th_L^*:H^1(\U(1),\Z)\ra H^1(L,\Z)$. The {\it Maslov class\/}
$\mu_L\in H^1(L;\Z)$ of $L$ is the image under $\Th_L^*$ of the
generator of $H^1(\U(1),\Z)\cong\Z$. If $H^1(M,\R)=0$ then $\mu_L$
depends only on $(M,\om),L$ and not on $g,J,\Om$. We call $L$ {\it
Maslov zero} if $\mu_L=0$.

A {\it grading\/} or {\it phase function\/} of an oriented
Lagrangian $L$ is a smooth function $\th_L:L\ra\R$ with
$\Th_L=\exp(i\th_L)$, so that $\Om\vert_L=e^{i\th_L}\d V_L$. That
is, $i\th_L$ is a continuous choice of logarithm for $\Th_L$.
Gradings exist if and only if $L$ is Maslov zero. If $L$ is
connected then gradings are unique up to addition of $2\pi n$ for
$n\in\Z$. A {\it graded Lagrangian\/} $(L,\th_L)$ in $M$ is an
oriented Lagrangian $L$ with a grading $\th_L$. Usually we refer to
$L$ as the graded Lagrangian, leaving $\th_L$ implicit.

An oriented Lagrangian $L$ in $M$ is called {\it almost
calibrated\/} if $(\cos\phi\,\Re\Om-\sin\phi\,\Im\Om)\vert_L$ is a
positive $m$-form on $L$ for some $\phi\in\R$. Then $L$ admits a
unique grading $\th_L$ taking values in $(\phi-\frac{\pi}{2},
\phi+\frac{\pi}{2})$. If a graded Lagrangian $L$ has phase variation
less than $\pi$, then it is almost calibrated.

An oriented Lagrangian $L$ in $M$ is called {\it special Lagrangian
with phase\/} $e^{i\phi}$ if $\Th_L$ is constant with value
$e^{i\phi}\in\U(1)$. If we do not specify a phase, we usually mean
phase 1. We will write {\it SL\/} for special Lagrangian, and {\it
SL\/ $m$-fold\/} for special Lagrangian submanifold. SL $m$-folds
with phase $e^{i\phi}$ are Maslov zero, and graded with phase
function $\th_L=\phi$. They are minimal submanifolds in $(M,g)$.
Compact SL $m$-folds are volume-minimizing in their homology class.
\label{bs2def1}
\end{dfn}

Special Lagrangians were introduced by Harvey and Lawson \cite[\S
III]{HaLa}. The deformation theory of SL $m$-folds was studied by
McLean \cite[\S 3]{McLe}:

\begin{thm} Let\/ $(M,J,g,\Om)$ be a Calabi--Yau $m$-fold, and\/ $L$
a compact SL\/ $m$-fold in $M$. Then the moduli space $\cM_{\sst L}$
of special Lagrangian deformations of\/ $L$ is a smooth manifold of
dimension $b^1(L),$ the first Betti number of\/~$L$.
\label{bs2thm1}
\end{thm}

\subsection{Special Lagrangian $m$-folds in $\C^m$}
\label{bs22}

\begin{dfn} Let $\C^m$ have coordinates $(z_1,\dots,z_m)$ and
complex structure $J$, and define a K\"ahler metric $g$, K\"ahler
form $\om$ and $(m,0)$-form $\Om$ on $\C^m$ by
\e
\begin{split}
g=\ms{\d z_1}+\cdots+\ms{\d z_m},\quad
\om&=\ts\frac{i}{2}(\d z_1\w\d\bar z_1+\cdots+\d z_m\w\d\bar z_m),\\
\text{and}\quad\Om&=\d z_1\w\cdots\w\d z_m.
\end{split}
\label{bs2eq2}
\e
Then $(\C^m,J,g,\Om)$ is the simplest example of a Calabi--Yau
$m$-fold.

Define a real 1-form $\la$ on $\C^m$ called the {\it Liouville
form\/} by
\begin{equation*}
\la=-\ha\Im(z_1\d\bar z_1+\cdots+z_m\d\bar z_m).
\end{equation*}
Then $\d\la=\om$. Thus, if $L$ is a Lagrangian in $\C^m$ then
$\d(\la\vert_L)=0$. We call $L$ an {\it exact\/} Lagrangian if
$\la\vert_L=\d f$ for some smooth $f:L\ra\R$.

A (singular) Lagrangian $C$ in $\C^m$ is called a {\it cone} if
$C=tC$ for all $t>0$, where $tC=\{t\,{\bf z}:{\bf z} \in C\}$. Let
$C$ be a closed Lagrangian cone in $\C^m$ with an isolated
singularity at 0. Then $\Si=C\cap{\cal S}^{2m-1}$ is a compact,
nonsingular Legendrian $(m\!-\!1)$-submanifold of ${\cal S}^{2m-1}$,
not necessarily connected. Let $g_\sSi$ be the metric on $\Si$
induced by the metric $g$ on $\C^m$ in \eq{bs2eq2}, and $r$ the
radius function on $\C^m$. Define $\iota:\Si\t(0,\iy)\ra\C^m$ by
$\iota(\si,r)=r\si$. Then the image of $\iota$ is $C\sm\{0\}$, and
$\iota^*(g)=r^2g_\sSi+\d r^2$ is the cone metric on~$C\sm\{0\}$.

Let $L$ be a closed, nonsingular Lagrangian $m$-fold in $\C^m$, e.g.
$L$ could be special Lagrangian, or a Lagrangian LMCF expander, as in Definition \ref{bs2def3} below. We call $L$
{\it asymptotically conical (AC)} with {\it rate} $\rho<2$ and {\it
cone} $C$ if there exists a compact subset $K\subset L$ and a
diffeomorphism $\vp:\Si\t(T,\iy)\ra L\sm K$ for some $T>0$, such
that
\begin{equation*}
\bmd{\nabla^k(\vp-\iota)}=O(r^{\rho-1-k}) \quad\text{as $r\ra\iy$,
for all $k=0,1,2,\ldots.$}
\end{equation*}
Here $\nabla,\md{\,.\,}$ are computed using the cone metric $\iota^*(g)$. Note that if $\rho<\si<2$ and $L$ is AC with rate $\rho$, then $L$ is also AC with rate~$\si$. 
\label{bs2def2}
\end{dfn}

Asymptotically conical special Lagrangians are an important class of
SL $m$-folds in $\C^m$. McLean's Theorem, Theorem \ref{bs2thm1},
was generalized to AC SL $m$-folds by Marshall \cite{Mars} and
Pacini \cite{Paci}. Here is a special case of their results:

\begin{thm} Let\/ $L$ be an asymptotically conical SL\/ $m$-fold in
$\C^m$ for $m\ge 3$ with cone $C$ and rate $\rho\in (2-m,0),$ and
write $\cM_{\sst L}^\rho$ for the moduli space of deformations of
$L$ as an AC SL\/ $m$-fold in $\C^m$ with cone $C$ and rate $\rho$.
Then $\cM_{\sst L}^\rho$ is a smooth manifold of dimension~$b^1_{\rm
cs}(L)=b^{m-1}(L)$.
\label{bs2thm2}
\end{thm}

Here $b^k(L)=\dim H^k(L,\R)$, $b^k_{\rm cs}(L)=\dim H^k_{\rm cs}(L,\R)$  for $H^*(L,\R),\ab H^*_{\rm cs}(L,\R)$ the (compactly-supported) cohomology of $L$, and $b^k_{\rm cs}(L)=b^{m-k}(L)$ by Poincar\'e duality. The next family of AC SL $m$-folds in $\C^m$ was first found by
Lawlor \cite{Lawl}, and rewritten by Harvey
\cite[p.~139--140]{Harv}. They are often called {\it Lawlor necks}.

\begin{ex} Let $m>2$ and $a_1,\ldots,a_m>0$, and define
polynomials $p,P$ by
\e
p(x)=(1+a_1x^2)\cdots(1+a_mx^2)-1 \quad\text{and}\quad
P(x)=\frac{p(x)}{x^2}.
\label{bs2eq3}
\e
Define real numbers $\phi_1,\ldots,\phi_m$ and $A$ by
\begin{equation*}
\phi_k=a_k\int_{-\iy}^\iy\frac{\d x}{(1+a_kx^2)\sqrt{P(x)}}
\quad\text{and}\quad A=\int_{-\iy}^\iy\frac{\d x}{2\sqrt{P(x)}}\,.
\end{equation*}
Clearly $\phi_k,A>0$. But writing $\phi_1+\cdots+\phi_m$ as one
integral gives
\begin{equation*}
\phi_1+\cdots+\phi_m=\int_0^\iy\frac{p'(x)\d x}{(p(x)+1)\sqrt{p(x)}}
=2\int_0^\iy\frac{\d w}{w^2+1}=\pi,
\end{equation*}
making the substitution $w=\sqrt{p(x)}$. So $\phi_k\in(0,\pi)$ and
$\phi_1+\cdots+\phi_m=\pi$. This yields a 1-1 correspondence between
$m$-tuples $(a_1,\ldots,a_m)$ with $a_k>0$, and $(m\!+\!1)$-tuples
$(\phi_1,\ldots,\phi_m,A)$ with $\phi_k\in (0,\pi)$,
$\phi_1+\cdots+\phi_m=\pi$ and~$A>0$.

For $k=1,\ldots,m$, define a function $z_k:\R\ra\C$ by
\begin{equation*}
z_k(y)={\rm e}^{i\psi_k(y)}\sqrt{a_k^{-1}+y^2}, \quad\text{where}\quad
\psi_k(y)=a_k\int_{-\iy}^y\frac{\d x}{(1+a_kx^2)\sqrt{P(x)}}\,.
\end{equation*}
Now write ${\bs\phi}=(\phi_1,\ldots,\phi_m)$, and define a
submanifold $L_{{\bs\phi},A}$ in $\C^m$ by
\begin{equation*}
L_{{\bs\phi},A}=\bigl\{(z_1(y)x_1,\ldots,z_m(y)x_m): y\in\R,\;
x_k\in\R,\; x_1^2+\cdots+x_m^2=1\bigr\}.
\end{equation*}

Then $L_{{\bs\phi},A}$ is closed, embedded, and diffeomorphic to
${\cal S}^{m-1}\t\R$, and Harvey \cite[Th.~7.78]{Harv} shows that
$L_{{\bs\phi},A}$ is special Lagrangian. Also $L_{{\bs\phi},A}$ is
asymptotically conical, with rate $\rho=2-m$ and cone $C$ the union
$\Pi_0\cup\Pi_{\bs\phi}$ of two special Lagrangian $m$-planes
$\Pi_0,\Pi_{\bs\phi}$ in $\C^m$ given by
\begin{equation*}
\Pi_0=\bigl\{(x_1,\ldots,x_m):x_j\in\R\bigr\},\;\>
\Pi_{\bs\phi}=\bigl\{({\rm e}^{i\phi_1}x_1,\ldots, {\rm
e}^{i\phi_m}x_m):x_j\in\R\bigr\}.
\end{equation*}

Apply Theorem \ref{bs2thm2} with $L=L_{{\bs\phi},A}$ and
$\rho\in(2-m,0)$. As $L\cong{\cal S}^{m-1}\t\R$ we have $b^1_{\rm
cs}(L)=1$, so Theorem \ref{bs2thm2} shows that $\dim\cM_{\sst
L}^\rho=1$. This is consistent with the fact that when $\bs\phi$ is
fixed, $L_{{\bs\phi},A}$ depends on one real parameter $A>0$. Here
$\bs\phi$ is fixed in $\cM_{\sst L}^\rho$ as the cone
$C=\Pi_0\cup\Pi_{\bs\phi}$ of $L$ depends on $\bs\phi$, and all
$\hat L\in\cM_{\sst L}^\rho$ have the same cone $C$, by definition.
\label{bs2ex1}
\end{ex}

Imagi, Oliveira dos Santos and the author \cite[Th.~1.1]{IJO} prove a
uniqueness theorem for Lawlor necks. The proof involves deep results on Lagrangian Floer cohomology and Fukaya categories, and was motivated by the
ideas of this paper. The exactness assumption is needed to apply the Fukaya category results. Note that special Lagrangians $L$ in $\C^m$ occurring as blow-ups of graded Lagrangian MCF in a Calabi--Yau $m$-fold (as in \S\ref{bs33} below) are automatically exact, so the exactness requirement does not matter for our programme. 

\begin{thm} Suppose $L$ is a closed, embedded, exact,
asymptotically conical special Lagrangian in $\C^m$ for $m\ge 3,$
asymptotic at rate $\rho<0$ to a union $\Pi_1\cup\Pi_2$ of two
transversely intersecting special Lagrangian planes $\Pi_1,\Pi_2$ in
$\C^m$. Then $L$ is equivalent under an $\SU(m)$ rotation to one of
the `Lawlor necks' $L_{\bs\phi,A}$ found by Lawlor {\rm\cite{Lawl},}
and described in Example\/ {\rm\ref{bs2ex1}}.
\label{bs2thm3}
\end{thm}

Here is an example based on Harvey and Lawson~\cite[\S III.3.A]{HaLa}:

\begin{ex} Define a special Lagrangian $T^2$-cone $C$ in $\C^3$ by
\e
C=\bigl\{(z_1,z_2,z_3)\in\C^3:\md{z_1}=\md{z_2}=\md{z_3},\;\>
z_1z_2z_3\in[0,\iy)\bigr\}.
\label{bs2eq4}
\e
This will be important in \S\ref{bs36} as it is a `stable' special
Lagrangian singularity in the sense of  \cite[Def.~3.6]{Joyc2}.
There are three families of explicit asymptotically conical SL
3-folds $L^A_1,L^A_2,L^A_3$ for $A>0$ in $\C^3,$ each diffeomorphic
to $\cS^1\t\R^2$ and asymptotic at rate $\rho=0$ to the cone $C$,
where
\e
L^A_1=\bigl\{(z_1,z_2,z_3)\in\C^3:\ms{z_1}-A=\ms{z_2}=\ms{z_3},\;\>
z_1z_2z_3\in[0,\iy)\bigr\},
\label{bs2eq5}
\e
and $L^A_2,L^A_3$ are obtained from $L^A_1$ by cyclic permutation
of~$z_1,z_2,z_3$.
\label{bs2ex2}
\end{ex}

\begin{ex} In \cite{Joyc6,Joyc7,Joyc8} we study SL 3-folds in
$\C^3$ invariant under the $\U(1)$-action
\begin{equation*}
{\rm e}^{i\th}:(z_1,z_2,z_3)\longmapsto
({\rm e}^{i\th}z_1,{\rm e}^{-i\th}z_2,z_3)
\quad\text{for ${\rm e}^{i\th}\in\U(1)$.}
\end{equation*}
The three papers are surveyed in \cite{Joyc9}. A $\U(1)$-invariant
SL 3-fold $N$ may locally be written in the form
\e
\begin{split}
N=\bigl\{(z_1,z_2,z_3)\in\C^3:\,& z_1z_2=v(x,y)+iy,\quad z_3=x+iu(x,y),\\
&\ms{z_1}-\ms{z_2}=2a,\quad (x,y)\in S\bigr\},
\end{split}
\label{bs2eq6}
\e
where $S$ is a domain in $\R^2$, $a\in\R$ and $u,v:S\ra\R$ satisfy
(in a weak sense if $a=0$) the {\it nonlinear Cauchy--Riemann
equations}
\e
\frac{\pd u}{\pd x}=\frac{\pd v}{\pd y}\quad\text{and}\quad
\frac{\pd v}{\pd x}=-2\bigl(v^2+y^2+a^2\bigr)^{1/2}\frac{\pd u}{\pd
y}.
\label{bs2eq7}
\e

If $S$ is simply-connected, as $\frac{\pd u}{\pd x}=\frac{\pd v}{\pd
y}$ there exists a potential $f$ for $u,v$ with $\frac{\pd f}{\pd
y}=u$, $\frac{\pd f}{\pd x}=v$, satisfying
\e
\Bigl(\Bigl(\frac{\pd f}{\pd x}\Bigr)^2+y^2+a^2\Bigr)^{-1/2}
\frac{\pd^2f}{\pd x^2}+2\,\frac{\pd^2f}{\pd y^2}=0.
\label{bs2eq8}
\e
In \cite{Joyc6,Joyc7}, for suitable strictly convex domains
$S\subset\R^2$ and boundary data $\phi:\pd S\ra\R$, we prove the
existence of a unique $f:S\ra\R$ satisfying \eq{bs2eq8} and
$f\vert_{\pd S}=\phi$, and then $u=\frac{\pd f}{\pd y}$,
$v=\frac{\pd f}{\pd x}$ satisfy \eq{bs2eq7} (possibly in a weak
sense if $a=0$), and $N$ in \eq{bs2eq6} is special Lagrangian.

When $v=y=a=0$, equations \eq{bs2eq7}--\eq{bs2eq8} become singular,
and the SL 3-fold $N$ in \eq{bs2eq6} has a singularity at
$(0,0,z_3)=\bigl(0,0,x+iu(x,0)\bigr)$ in $\C^3$. In the simplest
cases $N$ is locally modelled on the cone $C$ in \eq{bs2eq4} near
$(0,0,z_3)$, but there are also infinitely many other topological
types of singularities not locally modelled on cones. Note that the
existence and uniqueness results for $N$ are {\it entirely
independent\/} of the singularities appearing in the interior
of~$N$.

The following will be important in \S\ref{bs36}. Using the results
of \cite{Joyc6,Joyc7,Joyc8,Joyc9}, by choosing a suitable family
$\phi^t:t\in(-\ep,\ep)$ of boundary conditions for the potential
$f^t$, we can construct a family $N^t:t\in(-\ep,\ep)$ of exact
$\U(1)$-invariant SL 3-folds in $\C^3$ of the form \eq{bs2eq6} with
$a=0$, with the following properties:
\begin{itemize}
\setlength{\parsep}{0pt}
\setlength{\itemsep}{0pt}
\item[(i)] $N^t$ depends continuously on $t\in(-\ep,\ep)$ in a suitable sense, for instance as special Lagrangian integral currents in Geometric Measure Theory.
\item[(ii)] $N^t$ is nonsingular for $t<0$.
\item[(iii)] $N^0$ has one singular point at $(0,0,0)\in\C^3$, which
has tangent cone $\Pi_1\cup\Pi_2$, where $\Pi_1,\Pi_2$ are
$\U(1)$-invariant special Lagrangian planes in $\C^3$
intersecting non-transversely with $\Pi_1\cap\Pi_2=\R$.
\item[(iv)] $N^t$ for $t>0$ has two singular points at $(0,0,\pm
z (t))$, where $z(t)$ depends smoothly on $t$ and $z(t)\ra 0$ as
$t\ra 0$. Each singular point is locally modelled on the special
Lagrangian $T^2$-cone $C$ in \eq{bs2eq4}.
\end{itemize}
Thus, {\it isolated singular points of SL\/ $3$-folds modelled on
the $T^2$-cone $C$ in \eq{bs2eq4} can appear or disappear in pairs
under continuous deformation}.
\label{bs2ex3}
\end{ex}

\subsection{Lagrangian mean curvature flow}
\label{bs23}

Next we discuss (Lagrangian) mean curvature flow. A book on mean
curvature flow (MCF) for hypersurfaces in $\R^n$ is Mantegazza
\cite{Mant}. Two useful surveys on Lagrangian MCF are Smoczyk
\cite{Smoc2} and Neves~\cite{Neve2}.

Let $(M,g)$ be a Riemannian manifold, and $N$ a compact manifold
with $\dim N<\dim M$, and consider embeddings or immersions
$\io:N\hookra M$, so that $\io(N)$ is a submanifold of $M$. {\it
Mean curvature flow\/} ({\it MCF\/}) is the study of smooth
1-parameter families $\io_t$, $t\in[0,T)$ of such $\io_t:N\hookra M$
satisfying
\begin{equation*}
\frac{\d\io_t}{\d t}=H_{\io_t},
\end{equation*}
where $H_{\io_t}\in C^\iy(\io_t^*(TM))$ is the mean curvature of the
submanifold $\io_t:N\hookra M$. We usually write $N^t$ rather than $\io_t:N\hookra M$, suppressing the immersion, so that $\{N^t:t\in[0,T)\}$ is a family of submanifolds satisfying MCF. 

Mean curvature flow is the gradient flow of the volume functional for compact submanifolds $N$ in $M$. It has a unique short-time solution starting from any compact submanifold $N$.

Now let $(M,J,g,\Om)$ be a Calabi--Yau $m$-fold, and $L$ a compact Lagrangian submanifold in $M$. Then the mean curvature of $L$
is $H=J\nabla\Th_L$, where $\Th_L:L\ra\U(1)$ is the phase function
from Definition \ref{bs2def1}. Thus $H$ is an infinitesimal
deformation of $L$ as a Lagrangian. Smoczyk \cite{Smoc1} shows that
MCF starting from $L$ preserves the Lagrangian condition, yielding a
1-parameter family of Lagrangians $\{L^t:t\in[0,\ep)\}$ with $L^0=L$, which are all in the same Hamiltonian isotopy class if $L$ is Maslov zero. This is {\it Lagrangian mean curvature flow\/} ({\it LMCF\/}). Special Lagrangians are stationary points of Lagrangian MCF.

We will be especially interested in Lagrangian MCF for {\it
graded\/} Lagrangians. Suppose $\{L^t:t\in[0,T)\}$ is a family of
compact, graded Lagrangians satisfying Lagrangian MCF. Then $L^t$
are all Hamiltonian isotopic, that is, {\it graded Lagrangian MCF
stays within a fixed Hamiltonian isotopy class}. Also, if the phase
function $\th_{L^0}$ takes values in an interval $[a,b]$ or $(a,b)$,
then so does $\th_{L^t}$ for $t\in[0,T)$. Thus, {\it Lagrangian MCF
preserves the almost calibrated condition}.

It is an important problem to understand the singularities which
arise in Lagrangian mean curvature flow. Singularities in Lagrangian
MCF are often locally modelled on {\it soliton solutions},
Lagrangians in $\C^m$ which move by rescaling or translation under
Lagrangian MCF.

\begin{dfn} A closed Lagrangian $L$ in $\C^m$ is called an {\it LMCF
expander\/} if $H=\al F^\perp$ in $C^\iy(T\C^m\vert_L)$, where $H$
is the mean curvature of $L$ and $F^\perp$ is the orthogonal
projection of the position vector $F$ (that is, the inclusion
$F:L\hookra\C^m$) to the normal bundle $TL^\perp\subset
T\C^m\vert_L$, and $\al>0$ is constant.

This implies that (after reparametrizing by diffeomorphisms of $L$)
the family of Lagrangians $L^t:=\sqrt{2\al t}\,L$ for $t\in(0,\iy)$
satisfy Lagrangian mean curvature flow. That is, Lagrangian MCF
expands $L$ by dilations.

Similarly, we call $L$ an {\it LMCF shrinker\/} if $H=\al F^\perp$
for $\al<0$, and then $L^t:=\sqrt{2\al t}\,L$ for $t\in(-\iy,0)$
satisfy LMCF, so LMCF shrinks $L$ by dilations.

We call $L$ an {\it LMCF translator\/} if $H=v^\perp$, where
$v\in\C^m$ is the {\it translating vector\/} of $L$, and $v^\perp$
the orthogonal projection of $v$ to $TL^\perp$. Then $L^t:=L+tv$ for
$t\in\R$ satisfy LMCF, so Lagrangian MCF translates $L$ in $\C^m$.
\label{bs2def3}
\end{dfn}

Finite time singularities of MCF  have a fundamental division into
`type I' and `type II' singularities:

\begin{dfn} Let $(M,g)$ be a compact Riemannian manifold (e.g. a
Calabi--Yau $m$-fold) and $\{L^t:t\in[0,T)\}$ a family of compact
immersed submanifolds in $M$ (e.g. Lagrangians) satisfying mean
curvature flow. We say that the family {\it has a finite time
singularity at\/} $t=T$ if the flow cannot be smoothly continued to
$[0,T+\ep)$ for any $\ep>0$. As in Wang \cite[Lem.~5.1]{Wang} this
implies that $\mathop{\rm lim\,sup}_{t\ra T}\nm{A^t}_{C^0}\ra\iy$,
where $A^t$ is the second fundamental form of~$L^t$.

We call such a finite time singularity {\it of type I\/} if
$\nm{A^t}_{C^0}^2\le C/(T-t)$ for some $C>0$ and all $t\in [0,T)$.
Otherwise we call the singularity {\it of type II}.

We call $x\in M$ a {\it singular point\/} of the flow if
$\mathop{\rm lim\,sup}_{t\ra T}\nm{A^t\vert_{U\cap L^t}}_{C^0}=\iy$
for all open neighbourhoods $U$ of $x$ in $M$.
\label{bs2def4}
\end{dfn}

Huisken \cite{Huis} showed that type I singularities developing a
singularity at $x\in M$ are locally modelled in a strong sense on
MCF shrinkers in $\R^n=T_xM$, through a process known as `type I
blow up', as in Smoczyk \cite[Prop.~3.17]{Smoc2} or
Mantegazza~\cite[\S 3]{Mant}.

However, we are interested in MCF of {\it graded\/} Lagrangians in
Calabi--Yau $m$-folds, and it turns out that type I singularities do
not occur in graded Lagrangian MCF, as was proved by Wang
\cite[Rem.~5.1]{Wang} and Chen and Li \cite[Cor.~6.7]{ChLi} in the
almost calibrated case (i.e.\ Lagrangians $L^t$ with phase
variation less than $\pi$) and by Neves \cite[Th.~A]{Neve1} in the
graded (or Maslov zero) case.

\begin{thm} Let\/ $(M,J,g,\Om)$ be a compact Calabi--Yau $m$-fold
and\/ $\{L^t:t\in[0,T)\}$ a family of compact, immersed, graded
Lagrangians in $M$ satisfying Lagrangian MCF. Then the flow cannot
develop a type I singularity.
\label{bs2thm4}
\end{thm}

A parallel result of Neves \cite[Cor.~3.5]{Neve2} says that there
exist no nontrivial, immersed, graded Lagrangian MCF shrinkers in
$\C^m$ (satisfying a few extra conditions such as closed in $\C^m$
and of bounded Lagrangian angle), so there are no possible local
models for type I blow ups of graded Lagrangian MCF. Examples of
Lagrangian MCF shrinkers in $\C^m$ can be found in Abresch and
Langer \cite{AbLa} for $m=1$ and in Anciaux \cite{Anci} and Joyce,
Lee and Tsui \cite[Th.~F]{JLT} in higher dimensions, but none of
them are graded.

So, for graded Lagrangian MCF, all finite time singularities are of
type II. It is a well known `folklore' theorem that type II
singularities of MCF admit `type II blow ups', eternal smooth
solutions of MCF in $\R^n$ modelling the formation of the
singularity in the small region where the second fundamental form
$A^t$ is largest as $t\ra T$. The idea of type II blow ups is due to
Hamilton, and explanations can be found in Smoczyk \cite[\S 3.4]{Smoc2} and Mantegazza \cite[\S 4.1]{Mant}, and for Lagrangian
MCF in Han and Li \cite[\S 2]{HaLi}. We state it for graded LMCF:

\begin{thm} Let\/ $(M,J,g,\Om)$ be a compact Calabi--Yau $m$-fold
and\/ $\{L^t:t\in[0,T)\}$ a family of compact, immersed, graded
Lagrangians in $M$ satisfying Lagrangian MCF, with a finite time
singularity at\/ $t=T$. Then at some singular point\/ $x\in M$ of
the flow there exists a \begin{bfseries}type II blow
up\end{bfseries}.

That is, identifying $M$ near $x$ with\/ $T_xM\cong\C^m$ near $0,$
there exist sequences $(t_i)_{i=1}^\iy$ in $[0,T),$
$(x_i)_{i=1}^\iy$ in $M$ and\/ $(\la_i)_{i=1}^\iy$ in $(0,\iy),$
such that\/ $t_i\ra T,$ $x_i\ra x,$ $\la_i\ra\iy$ and\/
$\la_i^2(T-t_i)\ra 0$ as $i\ra\iy,$ and for each\/ $s\in\R$ the
limit
\begin{equation*}
\ti L^s=\lim_{i\ra\iy} \la_i\cdot (L^{t_i+\la_i^{-2}s}-x_i)
\end{equation*}
exists as a nonempty, noncompact, smooth, closed, immersed, exact,
graded Lagrangian in $\C^m$ whose mean curvature $\ti A^s$ is
nonzero (so that\/ $\ti L^s$ is not a union of Lagrangian planes
in\/ $\C^m$). All derivatives of\/ $\ti A^s,$ and the phase function
$\th_{\smash{\ti L^s}},$ are uniformly bounded independently of\/
$s\in\R$. Also $\ti L^s$ depends smoothly on $s\in\R,$ and\/ $\{\ti
L^s:s\in\R\}$ satisfies Lagrangian MCF in~$\C^m$.
\label{bs2thm5}
\end{thm}

A solution $\{\ti L^s:s\in\R\}$ of MCF for all $s\in\R$ is called an
{\it eternal solution}. Two obvious classes of eternal solutions of
Lagrangian MCF in $\C^m$ are
\begin{itemize}
\setlength{\parsep}{0pt}
\setlength{\itemsep}{0pt}
\item[(a)] $\ti L^s=L$ is independent of $s\in\R$, and is an SL
$m$-fold in $\C^m$.
\item[(b)] $\ti L^s=L+sv$ for $s\in\R$, where $L$ is a
Lagrangian MCF translator in $\C^m$ with translating vector
$v\in\C^m$.
\end{itemize}
Many examples of special Lagrangian $m$-folds in $\C^m$ are known
suitable for use in (a), but for (b) there are few, as we explain
in~\S\ref{bs24}.

\subsection{Examples of solitons for Lagrangian MCF}
\label{bs24}

We now give examples of solitons for Lagrangian MCF. We are
interested in graded Lagrangians, and as in \S\ref{bs23} there are
no graded Lagrangian MCF shrinkers. The next example describes a
family of LMCF expanders from Joyce, Lee and Tsui \cite[Th.s C \&
D]{JLT}, generalizing the `Lawlor necks' of Example~\ref{bs2ex1}.

\begin{ex} Let $m>2$, $\al\ge 0$ and $a_1,\ldots,a_m>0$, and define
a smooth function $P:\R\ra\R$ by $P(0)=\al+a_1+\cdots+a_m$ and
\e
P(x)=\ts\frac{1}{x^2}\bigl(e^{\al
x^2}\prod_{k=1}^m(1+a_kx^2)-1\bigl),\quad x\ne 0.
\label{bs2eq9}
\e
Define real numbers $\phi_1,\ldots,\phi_m$ by
\begin{equation*}
\phi_k=a_k\int_{-\iy}^\iy\frac{\d x}{(1+a_kx^2)\sqrt{P(x)}}\,,
\end{equation*}
For $k=1,\ldots,m$ define a function $z_k:\R\ra\C$ by
\begin{equation*}
z_k(y)={\rm e}^{i\psi_k(y)}\sqrt{a_k^{-1}+y^2}, \;\>\text{where}\;\>
\psi_k(y)=a_k\int_{-\iy}^y\frac{\d x}{(1+a_kx^2)\sqrt{P(x)}}\,.
\end{equation*}
Now write ${\bs\phi}=(\phi_1,\ldots,\phi_m)$, and define a
submanifold $L_{\bs\phi}^\al$ in $\C^m$ by
\begin{equation*}
L_{\bs\phi}^\al=\bigl\{(z_1(y)x_1,\ldots,z_m(y)x_m): y\in\R,\;
x_k\in\R,\; x_1^2+\cdots+x_m^2=1\bigr\}.
\end{equation*}
Then $L_{\bs\phi}^\al$ is a closed, embedded Lagrangian
diffeomorphic to $\cS^{m-1}\t\R$ and satisfying $H=\al F^\perp$. If
$\al>0$ it is an LMCF expander, and if $\al=0$ it is one of the
Lawlor necks $L_{{\bs\phi},A}$ from Example \ref{bs2ex1}. It is
graded, with Lagrangian angle
\begin{equation*}
\th_{L_{\bs\phi}^\al}\bigl((z_1(y)x_1,\ldots,z_m(y)x_m)\bigr)
=\ts\sum_{k=1}^m\psi_k(y)+\arg\bigl(-y-iP(y)^{-1/2}\bigr).
\end{equation*}
Note that the only difference between the constructions of
$L_{{\bs\phi},A}$ in Example \ref{bs2ex1} and $L_{\bs\phi}^\al$
above is the term $e^{\al x^2}$ in \eq{bs2eq9}, which does not
appear in \eq{bs2eq3}. If $\al=0$ then $e^{\al x^2}=1$, and the two
constructions agree.

As in \cite[Th.~D]{JLT}, $L_{\bs\phi}^\al$ is asymptotically conical, with cone $C$ the union $\Pi_0\cup\Pi_{\bs\phi}$ of two
Lagrangian $m$-planes $\Pi_0,\Pi_{\bs\phi}$ in $\C^m$ given by
\begin{equation*}
\Pi_0=\bigl\{(x_1,\ldots,x_m):x_j\in\R\bigr\},\;\>
\Pi_{\bs\phi}=\bigl\{({\rm e}^{i\phi_1}x_1,\ldots, {\rm
e}^{i\phi_m}x_m):x_j\in\R\bigr\}.
\end{equation*}
But in contrast to Example \ref{bs2ex1}, for $\al>0$ we do not have
$\phi_1+\cdots+\phi_m=\pi$, so $\Pi_{\bs\phi}$ and $C$ are not
special Lagrangian.

In \cite[Th.~D]{JLT} we prove that for fixed $\al>0$, the map
$\Phi^\al:(a_1,\ldots,a_m) \mapsto(\phi_1,\ldots,\phi_m)$ gives a
diffeomorphism
\begin{equation*}
\Phi^\al:(0,\iy)^m\longra\bigl\{(\phi_1,\ldots,\phi_m)\in(0,\pi)^m:
0<\phi_1+\cdots+\phi_m<\pi\bigr\}.
\end{equation*}
That is, for all $\al>0$ and ${\bs\phi}=(\phi_1,\ldots,\phi_m)$ with
$0<\phi_1,\ldots,\phi_m<\pi$ and $0<\phi_1+\cdots+\phi_m<\pi$, the
above construction gives a unique LMCF expander $L_{\bs\phi}^\al$ asymptotic to~$\Pi_0\cup\Pi_{\bs\phi}$.
\label{bs2ex4}
\end{ex}

Motivated by the ideas of this paper, Imagi, Oliveira dos Santos and the author \cite[Th.~1.1]{IJO} prove a uniqueness theorem for these LMCF expanders when $m\ge 3$. The case $m=2$ was already proved by Lotay and Neves \cite{LoNe}. For LMCF expanders, being exact is equivalent to being graded, or Maslov zero.

\begin{thm} Suppose $L$ is a closed, embedded, exact,
asymptotically conical Lagrangian MCF expander in $\C^m$ for $m\ge
2,$ satisfying the expander equation $H=\al F^\perp$ for $\al>0,$
and asymptotic at rate $\rho<2$ to a union $\Pi_1\cup\Pi_2$ of two
transversely intersecting Lagrangian planes $\Pi_1,\Pi_2$ in $\C^m$.
Then $L$ is equivalent under a $\U(m)$ rotation to one of the LMCF
expanders $L_{\bs\phi}^\al$ found by Joyce, Lee and Tsui\/
{\rm\cite[Th.s C \& D]{JLT},} and described in Example\/
{\rm\ref{bs2ex4}}.
\label{bs2thm6}
\end{thm}

\begin{ex} In dimension $m=1$, the unique connected Lagrangian MCF
translator in $\C$, up to rigid motions and rescalings, is the `grim
reaper'
\begin{equation*}
\bigl\{x+iy\in\C: y\in(-\pi/2,\pi/2),\quad x=-\log\cos y\bigr\},
\end{equation*}
with translating vector $v=1\in\C$, which is sketched in Figure
\ref{bs2fig1}.
\begin{figure}[htb]
\centerline{$\splinetolerance{.8pt}
\begin{xy}
0;<1mm,0mm>:
,(60,14);(60,-14)**\crv{(15,13.5)&(-5,10)&(-13,5)&(-15,0)(-13,-5)
&(-5,-10)&(15,-13.5)}
,(25,0)*{\text{MCF translates in this direction $\longra$}}
\end{xy}$}
\caption{`Grim reaper' Lagrangian MCF translating soliton in $\C$}
\label{bs2fig1}
\end{figure}
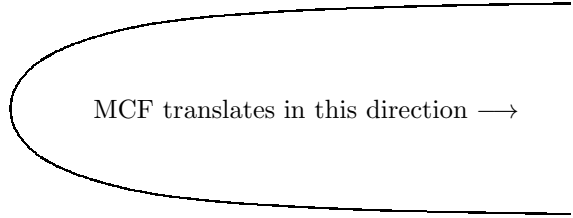
\label{bs2ex5}
\end{ex}

Here is a family of LMCF translators from Joyce, Lee and
Tsui~\cite[Cor.~I]{JLT}:

\begin{ex} For given constants $\al>0$ and
$a_1,\ldots,a_{m-1}>0,$ define
\begin{equation*}
\psi_j(y)=\int_{-\iy}^y\frac{\d t}{(\frac{1}{a_j}+t^2)\sqrt{P(t)}}\,,
\;\>\text{where}\;\>
P(t)=\frac{1}{t^2}\bigg(\prod_{k=1}^{m-1}(1+a_kt^2)e^{\al
t^2}-1\bigg),
\end{equation*}
for $j=1,\ldots,m-1$ and $y\in\R$. Then
\begin{align}
L=&\bigl\{\bigl(x_1\ts\sqrt{\frac{1}{a_1}\!+\!y^2}\,e^{i\psi_1(y)},
\ldots,x_{m-1}\sqrt{\frac{1}{a_{m-1}}\!+\!y^2}\,e^{i\psi_{m-1}(y)},
\ts\ha y^2\!-\!\ha\sum_{j=1}^{m-1}x_j^2
\nonumber\\
&-\ts\frac{i}{\al}\sum_{j=1}^{m-1}\psi_j(y)-\ts\frac{i}{\al}
\arg(y+iP(y)^{-1/2})\bigr):x_1,\ldots,x_{m-1},y\in\R\bigr\}
\label{bs2eq10}
\end{align}
is a closed, embedded Lagrangian in $\C^m$ diffeomorphic to $\R^m,$
which is a Lagrangian MCF translator with translating vector
$(0,\ldots,0,\al)\in\C^m$.

Define $\phi_1,\ldots,\phi_{m-1}\in\R$ by 
\begin{equation*}
\phi_j=\int_{-\iy}^\iy\frac{\d t}{(\frac{1}{a_j}+t^2)\sqrt{P(t)}}\,.
\end{equation*}
Then $\phi_1,\ldots,\phi_{m-1}\in(0,\pi)$ with $\phi_1+\cdots+\phi_{m-1}<\pi$, and $\psi_j(y)\ra\phi_j$ as $y\ra\iy$, and $\psi_j(y)\ra 0$ as $y\ra -\iy$. For fixed $\al>0,$ the map
$(a_1,\ldots,a_{m-1})\mapsto (\phi_1,\ldots,\phi_{m-1})$ is
a {\rm 1-1} correspondence from $(0,\iy)^{m-1}$ to
$\bigl\{(\phi_1,\ldots,\phi_{m-1})\in
(0,\pi)^{m-1}:\phi_1+\cdots+\phi_{m-1}<\pi\bigr\}$.

The phase function $\th_L$ of $L$ in \eq{bs2eq10} is a monotone decreasing function of $y$ only, with limits $\pi$ as $y\ra -\iy$ and $\sum_{j=1}^{m-1}\phi_j$ as $y\ra +\iy$. Thus, by choosing $\sum_{j=1}^{m-1}\phi_j$ close to $\pi,$ the phase variation of $L$ can be made arbitrarily small.

We can give the following heuristic description of $L$ in
\eq{bs2eq10}. If $y\gg 0$ then $\psi_j(y)\approx\phi_j$ and
$\sqrt{\frac{1}{a_j}+y^2}\approx y$, and the terms $-\frac{i}{\al}
\sum_{j=1}^n\psi_j(y)-\frac{i}{\al}\arg(y+iP(y)^{-1/2})$ are
negligible compared to $\ha y^2$ in the last coordinate. Thus, the
region of $L$ with $y\!\gg\! 0$ is in a weak sense approximate to
\begin{equation*}
\bigl\{\bigl(x_1ye^{i\phi_1},\ldots,
x_{m-1}ye^{i\phi_{m-1}},\ts\ha y^2-\ha\sum_{j=1}^{m-1}x_j^2
\bigr):x_1,\ldots,x_{m-1}\in\R,\; y>0\bigr\}.
\end{equation*}
But this is just an unusual way of parametrizing
\begin{equation*}
\Pi_{\bs\phi}=\bigl\{\bigl(y_1e^{i\phi_1},\ldots,
y_{m-1}e^{i\phi_{m-1}},y_m\bigr):y_j\in\R\bigr\}\sm
\bigl\{(0,\ldots,0,y_m):y_m\le 0\bigr\},
\end{equation*}
the complement of a ray in a Lagrangian plane. Similarly, the region
of $L$ with $y\ll 0$ is in a weak sense approximate to
\begin{equation*}
\Pi_0=\bigl\{(y_1,\ldots,
y_{m-1},y_m):y_j\in\R\bigr\}\sm
\bigl\{(0,\ldots,0,y_m):y_m\le 0\bigr\}.
\end{equation*}

So, $L$ can be roughly described as asymptotic to the union of two
Lagrangian planes $\Pi_0,\Pi_{\bs\phi}\cong\R^m$ which intersect in an $\R$ in $\C^m$, the $y_m$-axis $\bigl\{(0,\ldots,0,y_m):y_m\in\R\bigr\}$. To
make $L$, we glue these Lagrangian planes by a kind of `connect sum'
along the negative $y_m$-axis $\bigl\{(0,\ldots,0,y_m):y_m\le
0\bigr\}$. Under Lagrangian mean curvature flow, $\Pi_0,\Pi_{\bs\phi}$ remain fixed, but the gluing region translates in the positive $y_m$
direction, as though $\Pi_0,\Pi_{\bs\phi}$ are being `zipped together'.

A slightly more accurate description of the ends of $L$ for large $y$ is that $L$ approximates $\ti\Pi_{\bs\phi}$ when $y\gg 0$ and $\ti\Pi_0$ when $y\ll 0$, where $\ti\Pi_{\bs\phi}$ and $\ti\Pi_0$ are the non-intersecting affine Lagrangian planes in $\C^m$ 
\e
\begin{split}
\ti\Pi_{\bs\phi}&=\bigl\{\bigl(y_1e^{i\phi_1},\ldots,
y_{m-1}e^{i\phi_{m-1}},y_m\!-\!\ts\frac{i}{\al}(\phi_1\!+\!\cdots\!+\!\phi_{m-1})\bigr):y_j\!\in\!\R\bigr\},\\
\ti\Pi_0&=\bigl\{\bigl(y_1,\ldots,
y_{m-1},y_m-\ts\frac{i\pi}{\al}\bigr):y_j\in\R\bigr\}.
\end{split}
\label{bs2eq11}
\e
We will discuss these Lagrangian MCF translators further in
Example~\ref{bs3ex7}.
\label{bs2ex6}
\end{ex}

Castro and Lerma \cite{CaLe} give more examples of Lagrangian MCF
translators in $\C^2$. Neves and Tian \cite{NeTi} prove some
nonexistence results.

\subsection{Lagrangian Floer cohomology and Fukaya categories}
\label{bs25}

Let $(M,J,g,\Om)$ be a Calabi--Yau $m$-fold, which may be compact or
noncompact, with K\"ahler form $\om$. We now explain a little about
(embedded) {\it Lagrangian branes\/} $(L,E)$ in $(M,\om)$, {\it bounding cochains\/} $b$ for $(L,E)$ and {\it obstructions to\/} $HF^*$, {\it Lagrangian Floer cohomology\/} $HF^*\bigl((L,E,b),(L',E',b')\bigr)$, the {\it Fukaya category\/} $\sF(M)$, and the {\it derived Fukaya category\/} $D^b\sF(M)$. Section \ref{bs26} discusses the extension of all this to {\it immersed\/} Lagrangians.

The construction of $D^b\sF(M)$ in the generality we need may not
yet be available in the literature. As this paper is wholly
conjecture anyway, and clearly the theory will eventually work, this
does not matter very much.

The version of bounding cochains, obstructions to $HF^*$, and
Lagrangian Floer cohomology we need is in Fukaya, Oh, Ohta and Ono
\cite{FOOO}. An early explanation of how to define the (derived)
Fukaya category $\sF(M),D^b\sF(M)$ is Fukaya \cite{Fuka1}, and a
more recent survey is Fukaya \cite{Fuka2}. Floer \cite{Floe}
originally introduced Lagrangian Floer cohomology.

For {\it exact\/} Lagrangians in Liouville manifolds (a class of
noncompact, exact symplectic manifolds), a simpler, more complete,
and more satisfactory theory of Lagrangian Floer cohomology and
Fukaya categories is given in Seidel \cite{Seid3}, which we used in
\cite{IJO} to prove Theorems \ref{bs2thm3} and \ref{bs2thm6}. In
Seidel's theory there are no bounding cochains or obstructions
to~$HF^*$.

However, for our purposes Seidel's theory (in its current form) will not do: we need to extend the theory to immersed Lagrangians, and even for exact Lagrangians, bounding cochains and obstructions to $HF^*$ will then
appear. This extension will be discussed briefly in \S\ref{bs41}. Also, we wish to stress the idea that Lagrangian MCF is better behaved for Lagrangians with $HF^*$ unobstructed, and in Seidel's framework this issue is hidden by restricting to exact, embedded Lagrangians, for which $HF^*$ is automatically unobstructed.

\begin{dfn} Fix a field $\F$, in which we will do `counting' of $J$-holomorphic curves. If nontrivial $J$-holomorphic $\CP^1$'s can exist in the symplectic manifold $(M,\om)$ we are interested in, the virtual counts can be rational, so $\F$ must have characteristic zero, and $\F=\Q,\R$ or $\C$ are the obvious possibilities. If $M$ has no $J$-holomorphic $\CP^1$'s (for example, if $\om$ is exact, or if $\pi_2(M)=0$) then $\F$ can be arbitrary, so we can take $\F=\Z_2$, for instance, which means we do not have to worry about orientations on moduli spaces of $J$-holomorphic curves.

The {\it Novikov ring\/} $\La_\nov$ is the field of formal
power series $\sum_{i=0}^\iy a_iP^{\la_i}$ for $a_i\in\F$ and
$\la_i\in\R$ with $\la_i\ra+\iy$ as $i\ra\iy$, for $P$ a formal
variable. Write $\La_\nov^{\ge 0}$ for the subring of $\sum_{i=0}^\iy
a_iP^{\la_i}$ in $\La_\nov$ with all $\la_i\ge 0$, and $\La_\nov^+$ for the ideal of $\sum_{i=0}^\iy a_iP^{\la_i}$ in $\La_\nov$ with all~$\la_i>0$.
\label{bs2def5}
\end{dfn}

\begin{dfn} Let $(M,J,g,\Om)$ be a Calabi--Yau $m$-fold. A {\it Lagrangian brane\/} in $M$ is a pair $(L,E)$, where $L$ is a compact, spin, graded Lagrangian in $M$, and $E\ra L$ is a {\it rank one\/ $\F$-local system\/} on $L$, for $\F$ as in Definition \ref{bs2def5}. That is, $E$ is a locally constant rank one $\F$-vector bundle over $L$, so that if $p\in L$ then $E\vert_p$ is a dimension one $\F$-vector space, which is locally independent of~$p$.

In this section we take $L$ to be embedded, but in \S\ref{bs26} $L$ can be immersed, and in \S\ref{bs3} we will (conjecturally) allow $L$ to have certain kinds of singularities. 
\label{bs2def6}
\end{dfn}

\begin{rem} `Lagrangian branes' are the objects for which we will define Lagrangian Floer cohomology and Fukaya categories; the term is used in the same way by Seidel \cite[\S 12a]{Seid3} and Haug \cite[\S 3.1]{Haug}, for instance, although with different definitions. Our definition is designed to try to make the programme of \S\ref{bs3} work. The precise details of Definition \ref{bs2def6} will be important in Remark \ref{bs3rem2} and \S\ref{bs34}, and are discussed in Remark~\ref{bs3rem5}.

If we take $\F=\C$ then $E\ra L$ is a complex line bundle on $L$ with a flat connection $\nabla_E$, which is determined up to isomorphism by its holonomy $\Hol(\nabla_E):\pi_1(L)\ra\C^*$. In Mirror Symmetry it is natural to suppose that $\nabla_E$ preserves a unitary metric on $E$, so that $\Hol(\nabla_E)$ takes values in $\U(1)\subset\C^*$. One can also allow $E$ to be an $\F$-local system of higher rank. Kontsevich \cite{Kont} and Fukaya \cite[\S 2.1]{Fuka1} include a unitary local system $E\ra L$ of arbitrary rank in objects of their Fukaya categories.

We need to restrict to $E$ of rank one, and not to impose the unitary condition.

Much of the literature on Lagrangian Floer cohomology and Fukaya categories including \cite{FOOO,Seid3} omits the local system $E\ra L$, which is equivalent to taking $E$ to be trivial, $E=\F\t L\ra L$. As in \S\ref{bs34}, we cannot do this, since in the programme of \S\ref{bs32} involving families $(L^t,E^t,b^t)$ for $t\in[0,\iy)$, starting with $E^0$ trivial, after a surgery at $t=T_i$, we can have $E^t$ nontrivial for~$t>T_i$.
\label{bs2rem1}
\end{rem}

\begin{dfn} Let $(M,J,g,\Om)$ be a Calabi--Yau $m$-fold, and $L,L'$
graded Lagrangians in $M$, with phase functions $\th_L,\th_{L'}$,
which intersect transversely at $p\in M$. By a kind of simultaneous
diagonalization, we may choose an isomorphism $T_pM\cong\C^m$ which
identifies $J\vert_p,g\vert_p,\om\vert_p$ on $T_pM$ with the
standard versions \eq{bs2eq2} on $\C^m$, and identifies $T_pL,T_pL'$
with the Lagrangian planes $\Pi_0,\Pi_{\bs\phi}$ in $\C^m$
respectively, where
\e
\Pi_0=\bigl\{(x_1,\ldots,x_m):x_j\in\R\bigr\},\;\>
\Pi_{\bs\phi}=\bigl\{({\rm e}^{i\phi_1}x_1,\ldots, {\rm
e}^{i\phi_m}x_m):x_j\in\R\bigr\},
\label{bs2eq12}
\e
for $\phi_1,\ldots,\phi_m\in(0,\pi)$. Then $\phi_1,\ldots,\phi_m$
are independent of choices up to order. Define the {\it degree\/}
$\mu_{L,L'}(p)\in\Z$ of $p$ by
\begin{equation*}
\mu_{L,L'}(p)=(\phi_1+\cdots+\phi_m+\th_L(p)-\th_{L'}(p))/\pi.
\end{equation*}
This an integer as $\th_{L'}(p)=\th_L(p)+\phi_1+\cdots+\phi_m\mod\pi\Z$. Exchanging $L,L'$ replaces $\phi_1,\ldots,\phi_m$ by $\pi-\phi_1,\ldots,\pi-\phi_m$, so that $\mu_{L,L'}(p)+\mu_{L',L}(p)=m$. Since $\phi_1,\ldots,\phi_m\in(0,\pi)$, we see that
\e
(\th_L(p)-\th_{L'}(p))/\pi<\mu_{L,L'}(p)<(\th_L(p)-\th_{L'}(p))/\pi+m.
\label{bs2eq13}
\e

\label{bs2def7}
\end{dfn}

Here is the basic idea of Lagrangian Floer cohomology. Let
$(L,E),(L',E')$ be Lagrangian branes in a Calabi--Yau $m$-fold $(M,J,g,\Om)$, and suppose $L,L'$ intersect transversely. The aim is to
define a $\La_\nov$-module $HF^*\bigl((L,E),(L',E')\bigr)$ called the {\it Lagrangian Floer cohomology}, which is the cohomology of a complex of $\La_\nov$-modules $\bigl(CF^*\bigl((L,E),(L',E')\bigr),\d\bigr)$ called the {\it Floer complex}.

\begin{figure}[htb]
\centerline{$\splinetolerance{.8pt}
\begin{xy}
0;<1mm,0mm>:
,(-20,0);(20,0)**\crv{(0,10)}
?(.95)="a"
?(.85)="b"
?(.75)="c"
?(.65)="d"
?(.55)="e"
?(.45)="f"
?(.35)="g"
?(.25)="h"
?(.15)="i"
?(.05)="j"
?(.5)="y"
,(-20,0);(-30,-6)**\crv{(-30,-5)}
,(20,0);(30,-6)**\crv{(30,-5)}
,(-20,0);(20,0)**\crv{(0,-10)}
?(.95)="k"
?(.85)="l"
?(.75)="m"
?(.65)="n"
?(.55)="o"
?(.45)="p"
?(.35)="q"
?(.25)="r"
?(.15)="s"
?(.05)="t"
?(.5)="z"
,(-20,0);(-30,6)**\crv{(-30,5)}
,(20,0);(30,6)**\crv{(30,5)}
,"a";"k"**@{.}
,"b";"l"**@{.}
,"c";"m"**@{.}
,"d";"n"**@{.}
,"e";"o"**@{.}
,"f";"p"**@{.}
,"g";"q"**@{.}
,"h";"r"**@{.}
,"i";"s"**@{.}
,"j";"t"**@{.}
,"y"*{<}
,"z"*{>}
,(-20,0)*{\bu}
,(-20,-3)*{p}
,(20,0)*{\bu}
,(20,-3)*{q}
,(0,0)*{\Si}
,(-32,4)*{L}
,(-32,-4)*{L'}
,(32,4)*{L}
,(32.5,-4)*{L'}
\end{xy}$}
\caption{Holomorphic disc $\Si$ with boundary in $L\cup L'$}
\label{bs2fig2}
\end{figure}
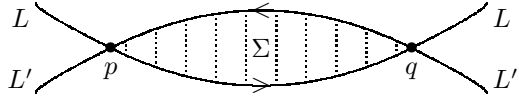

Define a free, graded $\La_\nov$-module $CF^*\bigl((L,E),(L',E')\bigr)$ by
\begin{equation*}
CF^k\bigl((L,E),(L',E')\bigr)=\bigop_{p\in L\cap L':\mu_{L,L'}(p)=k}
\Hom_\F\bigl(E\vert_p,E'\vert_p\bigr)\ot_\F\La_\nov.
\end{equation*}
Initially we define $\d:CF^k(L,L')\ra CF^{k+1}(L,L')$ by
\e
\d \al_p=\!\!\!\!\!\!\!\!\!\bigop_{\begin{subarray}{l} q\in L\cap L':\\
\mu_{L,L'}(q)=k+1\end{subarray}}\!\!\!\sum_{A>0} \bigl(\#_{\rm virt} \oM_{p,q}^A\bigr)\,P^A\cdot \mathop{PT}\limits_{\begin{subarray}{l} \text{$p\ra q$ in}\\ \pd\Si\cap L'\end{subarray}}(E')\ci \al_p\ci \mathop{PT}\limits_{\begin{subarray}{l} \text{$q\ra p$ in}\\ \pd\Si\cap L\end{subarray}}(E),
\label{bs2eq14}
\e
for $p\in L\cap L'$ with $\mu_{L,L'}(p)=k$ and $\al_p\in\Hom_\F\bigl(E\vert_p,E'\vert_p\bigr)\ot_\F\La_\nov$, where $\oM_{p,q}^A$ is the moduli space of stable $J$-holomorphic discs $\Si$ in $M$ with boundary in $L\cup L'$, corners at $p,q$ and area $A$, of the form shown in Figure \ref{bs2fig2}, where $\#_{\rm virt}\oM_{p,q}^A\in\Q$ is the `virtual number of points' in $\oM_{p,q}^A$, and the sum is weighted by composition with the {\it parallel transport maps\/} $PT_{\cdots}(E)\in \Hom_\F\bigl(E\vert_q,E\vert_p\bigr)$ and $PT_{\cdots}(E')\in \Hom_\F\bigl(E'\vert_p,E'\vert_q\bigr)$ in the $\F$-local systems $E,E'$ along the two segments of $\pd\Si$. These $PT_{\cdots}(E),PT_{\cdots}(E')$ are locally constant on~$\oM_{p,q}^A$.

Constructing an appropriate geometric structure (`Kuranishi space' or `polyfold') on $\oM_{p,q}^A$, and defining the virtual count $\#_{\rm virt}\oM_{p,q}^A$, raise many complicated issues which we will not go into.

For exact Lagrangians in an exact symplectic manifold, as in Seidel
\cite{Seid3}, the differential $\d$ in \eq{bs2eq14} has $\d^2=0$, so
$HF^*\bigl((L,E),(L',E')\bigr)$ is well-defined. However, in the non-exact case we may have $\d^2\ne 0$, because of contributions to the boundaries $\pd\oM_{p,q}^A$ from holomorphic discs with boundary in $L$ or in~$L'$.

To get round this, Fukaya, Oh, Ohta and Ono \cite[\S 3.6]{FOOO} introduce
the notion of a {\it bounding cochain\/} $b$ for $(L,E)$, an element
$b$ of the singular $(m-1)$-chains $C_{m-1}(L,\La_\nov^+)$ of $L$ with coefficients in $\La_\nov^+\subset\La_\nov$, satisfying an equation in $C_{m-2}(L,\La_\nov^+)$ which is (very roughly, and oversimplified) of the form
\e
\pd b+\sum_{k\ge 0}\sum_{A>0}P^A\cdot\bigl[\oM_{k+1}^A\t_{L^k} b^k\bigr]_{\rm virt}\cdot\Hol_{\pd\Si}(E)=0,
\label{bs2eq15}
\e
where $\oM_{k+1}^A$ is the moduli space (as a Kuranishi space or polyfold, of virtual dimension $m+k-2$) of isomorphism classes $[\Si,\vec z]$ where $\Si$ is a stable $J$-holomorphic disc of area $A>0$ in $M$ with boundary in $L$, and $\vec z=(z_0,z_1,\ldots,z_k)$ are cyclically ordered marked points in $\pd\Si$. Also $\oM_{k+1}^A\t_{L^k} b^k$ is the moduli space of such $[\Si,\vec z]$ in which $z_1,\ldots,z_k$ intersect the chain $b$ in $L$, and $\bigl[\oM_{k+1}^A\t_{L^k} b^k\bigr]{}_{\rm virt}$ is a virtual chain for this. The sum is weighted by the holonomy $\Hol_{\pd\Si}(E)\in\F^*$ of the rank one $\F$-local system $E\ra L$ around $\pd\Si\subset L$, which depends only on $[\pd\Si]\in H_1(L,\Z)$, and is locally constant on~$\oM_{k+1}^A$.

If a bounding cochain $b$ exists for $(L,E)$, we say that $(L,E)$ {\it has $HF^*$ unobstructed}, otherwise we say that $(L,E)$ {\it has $HF^*$ obstructed}. Implicitly we will always consider bounding cochains $b$ up to the appropriate notion of equivalence.

To oversimplify even further, suppose that the terms for $k\ge 1$ in \eq{bs2eq15} are zero, and $\pd\oM_1^A=\es$ for all $A>0$ when $k=0$, so that $\pd\bigl[\oM_1^A\bigr]{}_{\rm virt}=0$, and write $b=\sum_{A>0}P^A\cdot b_A$ for $b_A\in C_{m-1}(L,\Q)$. Then \eq{bs2eq15} becomes $\pd b_A=\bigl[\oM_1^A\bigr]{}_{\rm virt}\cdot\Hol_{\pd\Si}(E)$ for all $A>0$. So a bounding cochain $b$ exists if $\bigl[\bigl[\oM_1^A\bigr]{}_{\rm virt}\bigr]=0$ in $H_{m-2}(L,\Q)$ for all $A>0$. In particular, if $H_{m-2}(L,\Q)=0$ then a bounding cochain exists. 

In the general case, if $H_{m-2}(L,\Q)=0$ then \eq{bs2eq15} may be solved for $b=\sum_{A>0}P^A\cdot b_A$ by an inductive procedure in increasing $A$, yielding:

\begin{lem} Let\/ $(M,J,g,\Om)$ be a Calabi--Yau $m$-fold and\/
$(L,E)$ an embedded Lagrangian brane in $M$. If\/ $b_{m-2}(L)=0$ then $(L,E)$ has $HF^*$ unobstructed.

\label{bs2lem1}
\end{lem}

Suppose $b,b'$ are bounding cochains for $(L,E),(L',E')$. Then Fukaya et al.\ \cite{FOOO} define a modification $\d^{b,b'}$ of $\d$ in
\eq{bs2eq14} involving $b,b'$ and satisfying $(\d^{b,b'})^2=0$. The
{\it Lagrangian Floer cohomology\/} $HF^*\bigl((L,E,b),(L',E',b')\bigr)$
is the cohomology of $\bigl(CF^*\bigl((L,E),(L',E')\bigr),\d^{b,b'}\bigr)$, which may depend on $b,b'$. Here are some properties of Lagrangian Floer cohomology in the theory of Fukaya, Oh, Ohta and Ono~\cite{FOOO}:
\begin{itemize}
\setlength{\parsep}{0pt}
\setlength{\itemsep}{0pt}
\item[(a)] The Lagrangian Floer cohomology $HF^*\bigl((L,E,b),(L',E',b')\bigr)$ is independent of the choice of almost complex structure $J$ up to canonical isomorphism, although $\bigl(CF^*\bigl((L,E),(L',E')\bigr),\d^{b,b'}\bigr)$ does depend on~$J$.
\item[(b)] Let $(L^t,E^t):t\in[0,1]$ be a smooth family of Lagrangian branes, with the $L^t$ Hamiltonian isotopic and the $E^t$ locally constant in $t$, and $b^0$ a bounding cochain for $(L^0,E^0)$. Then Fukaya et al.\ \cite[Th.~B \& Th.~G]{FOOO} explain how to extend $b^0$ to a family of bounding cochains $b^t$ for $(L^t,E^t)$ for $t\in[0,1]$, by pushing forward $b^0$ along a Hamiltonian symplectomorphism of $(M,\om)$ taking $L^0$ to $L^t$. If $(L',E')$ is another Lagrangian brane with bounding cochain $b'$ then $HF^*\bigl((L^t,E^t,b^t),(L',E',b')\bigr)$ is independent of $t\in[0,1]$ up to canonical isomorphism. Thus $HF^*\bigl((L,E,b),\ab(L',E',b')\bigr)$ is an invariant of Lagrangian branes up to Hamiltonian isotopy.

This will be important in \S\ref{bs3}. If the $L^t$ for $t\in[0,1]$ are graded Lagrangians moving under Lagrangian MCF in a Calabi--Yau $m$-fold, then they are Hamiltonian isotopic. If $L^0$ extends to a triple $(L^0,E^0,b^0)$, then $L^t$ for $t\in[0,1]$ extends naturally to a family of Lagrangian branes $(L^t,E^t)$ with $E^t$ locally constant in $t$, and $b^0$ to a family of bounding cochains $b^t$ for $(L^t,E^t)$, so that Lagrangian MCF determines the flow of triples $(L^t,E^t,b^t)$, not just of Lagrangians $L^t$. We sometimes refer to this flow for $b^t$ as a kind of `parallel transport' on the space of bounding cochains for~$(L^t,E^t)$.
\end{itemize}

\begin{rem} We need $M$ to be (symplectic) Calabi--Yau and $L,L'$ to be graded to define the degree $\mu_{L,L'}(p)\in\Z$, which determines the grading of $CF^*\bigl((L,E),(L',E')\bigr)$ and $HF^*\bigl((L,E,b),(L',E',b')\bigr)$. If we took $M$ symplectic and $L,L'$ oriented, then $CF^*\bigl((L,E),(L',E')\bigr),HF^*\bigl((L,E,b),(L',E',b')\bigr)$ would only be graded over $\Z_2$ rather than~$\Z$.
\label{bs2rem2}
\end{rem}

Lagrangian Floer cohomology is only the beginning of a more general
theory of {\it Fukaya categories}, which may be still incomplete in
the general case. Let $(M,J,g,\Om)$ be a Calabi--Yau $m$-fold. The idea is to define the {\it Fukaya category\/} $\sF(M)$ of $M$, an $A_\iy$-{\it category\/} whose objects are triples $(L,E,b)$ of a Lagrangian brane $(L,E)$ in $M$ with $HF^*$ unobstructed, and a bounding cochain $b$ for $(L,E)$, such that the morphisms $\Hom\bigl((L_0,E_0,b_0),(L_1,E_1,b_1)\bigr)$ in $\sF(M)$ are the graded $\La_\nov$-modules $CF^*\bigl((L_0,E_0),(L_1,E_1)\bigr)$ from above, with $A_\iy$-operations
\e
\begin{split}
\mu^k:CF^{a_k}\bigl((L_{k-1},E_{k-1}),(L_k,E_k)\bigr)\t\cdots
\t CF^{a_1}\bigl((L_0,E_0),(L_1,E_1)\bigr)&\\
\longra CF^{a_1+\cdots+a_k+2-k}\bigl((L_0,E_0),(L_k,E_k)\bigr)&
\end{split}
\label{bs2eq16}
\e
for $k\ge 1$, with $\mu^1:CF^{a_1}\bigl((L_0,E_0),(L_1,E_1)\bigr)\ab\ra CF^{a_1+1}\bigl((L_0,E_0),(L_1,E_1)\bigr)$ the differential $\d^{b_0,b_1}$ in the Floer complex. The coefficients in the $\La_\nov$-multilinear map $\mu^k$ in \eq{bs2eq16} are obtained by `counting' $J$-holomorphic $(k\!+\!1)$-gons in $M$ with boundary in $L_0\cup\cdots\cup L_k$, weighted by parallel transport maps in~$E_0,\ldots,E_k$.

By a category theory construction, one then defines the {\it derived
Fukaya category}, a triangulated category. There are two versions, which we will write $D^b\sF(M)$ and $D^\pi\sF(M)$. For $D^b\sF(M)$, the objects are {\it twisted complexes}, as in  Seidel \cite[\S 3l]{Seid3}. Roughly speaking, a twisted complex consists of objects $(L_1,E_1,b_1),\ldots,(L_n,E_n,b_n)$ in $\sF(M)$ together with Floer cochains $b_{ij}\in CF^*\bigl((L_i,E_i),(L_j,E_j)\bigr)$ for $1\le i<j\le n$ satisfying an equation \cite[eq.~(3.19)]{Seid3} related to the bounding cochain equation. In particular, objects $(L,E,b)$ in $\sF(M)$ are also objects in~$D^b\sF(M)$.

The translation functor $[1]$ in the triangulated category $D^b\sF(M)$
acts on objects $(L,E,b)$ by reversing the orientation of $L$ and
changing the grading $\th_L$ to $\th_L+\pi$. The (graded) morphisms of objects $(L,E,b),(L',E',b')$ in $D^b\sF(M)$ are $\Hom^*\bigl((L,E,b),(L',E',b')\bigr)=HF^*\bigl((L,E,b),(L',E',b')\bigr)$. 

The second version $D^\pi\sF(M)$, called the {\it idempotent completion}, {\it Karoubi completion}, or {\it split closure\/} of $D^b\sF(M)$, is obtained by applying a further category theory construction to $D^b\sF(M)$, which adds direct summands (idempotents) of objects in $D^b\sF(M)$ as extra objects, as in Seidel~\cite[\S 4]{Seid3}.

Kontsevich's {\it Homological Mirror Symmetry Conjecture\/} \cite{Kont}, motivated by String Theory, says (very roughly) that if $M,\check M$ are `mirror' Calabi--Yau $m$-folds then there should be an equivalence of triangulated categories
\begin{equation*}
D^\pi\sF(M)\simeq D^b\mathop{\rm coh}(\check M).
\end{equation*}
This has driven much research in the area. 

For Mirror Symmetry, one must use $D^\pi\sF(M)$ rather than $D^b\sF(M)$, as the mirror category $D^b\coh(\check M)$ is automatically idempotent complete. In \S\ref{bs31} we will conjecture that in the situation we are interested in, our enlarged version of $D^b\sF(M)$ should be idempotent complete, so that~$D^\pi\sF(M)\simeq D^b\sF(M)$.

\subsection{$HF^*$ and $D^b\sF(M)$ for immersed Lagrangians}
\label{bs26}

For the programme of \S\ref{bs3}, it will be necessary to enlarge the derived Fukaya category $D^b\sF(M)$ of a Calabi--Yau $m$-fold $(M,J,g,\Om)$ to include {\it immersed\/} Lagrangians. As a first step in doing this, Akaho and the author \cite{AkJo} explain how to generalize the Lagrangian
Floer cohomology of Fukaya, Oh, Ohta and Ono \cite{FOOO} from
embedded Lagrangians to immersed Lagrangians with transverse
self-intersections. We now explain some of the main ideas
in~\cite{AkJo}.

Let $(M,J,g,\Om)$ be a Calabi--Yau $m$-fold, and $(L,E)$ a Lagrangian brane in $M$. As in \S\ref{bs25}, in the embedded case \cite{FOOO}, a bounding cochain for $(L,E)$ is a singular $(m\!-\!1)$-chain $b\in C_{m-1}(L,\La_\nov^+)$ (or equivalence class of such chains), satisfying an equation \eq{bs2eq15} involving virtual chains for moduli spaces $\oM_{k+1}^{\rm main}(J,\be)$ of $J$-holomorphic discs in $M$ with boundary in~$L$.

In the immersed case \cite{AkJo}, if $L$ has transverse
self-intersections, a bounding cochain $b$ for $(L,E)$ consists of
two pieces of data: a chain $b_{\rm ch}$ in $C_{m-1}(L,\La_\nov^+)$
as above, and also, for each point $p\in M$ at which two local
sheets $L_+,L_-$ of $L$ intersect transversely with
$\mu_{L_+,L_-}(p)=1$, an element
\e
b_p\in\Hom_\F\bigl(E_+\vert_p,E_-\vert_p\bigr)\ot_\F\La_\nov^{\ge 0},
\label{bs2eq17}
\e
where we write $E_\pm$ for the restriction of $E\ra L$ to the local sheets $L_\pm$. These $b_{\rm ch},b_p$ must satisfy equations involving virtual chains for moduli spaces of $J$-holomorphic discs in $M$ with boundary in $L$, but now these $J$-holomorphic discs can be polygons with `corners' at self-intersection points of~$L$.

For example, suppose $(L_1,E_1),(L_2,E_2)$ are embedded, transversely
intersecting Lagrangian branes in $M$. Then $(L,E)=(L_1,E_1)\cup (L_2,E_2)$ is an immersed Lagrangian brane in $M$. A bounding cochain $b$ for $(L,E)$ could consist of $b_{\rm ch}=b_1\op b_2$, where $b_i\in C_{m-1}(L_i,\La_\nov^+)$ for $i=1,2$ are embedded bounding cochains for $(L_1,E_1),(L_2,E_2)$, together with elements $b_p$ in \eq{bs2eq17} for $p\in L_1\cap L_2$ with $\mu_{L_1,L_2}(p)=1$ or $\mu_{L_2,L_1}(p)=1$ which encode how the objects $(L_1,E_1,b_1),(L_2,E_2,b_2)$ in $D^b\sF(M)$ are glued together to make $(L,E,b)$. For instance, if we have a distinguished triangle in $D^b\sF(M)$
\begin{equation*}
\smash{\xymatrix@C=33pt{ (L_1,E_1,b_1) \ar[r] & (L,E,b) \ar[r] & (L_2,E_2,b_2) \ar[r]^(0.45){\be} & (L_1,E_1,b_1)[1], }}
\end{equation*}
then the $b_p\in\Hom_\F\bigl(E_2\vert_p,E_1\vert_p\bigr)\ot_\F\La_\nov^{\ge 0}$ for $p\in L_1\cap L_2$ with $\mu_{L_2,L_1}(p)=1$
form a chain in $CF^1\bigl((L_2,E_2),(L_1,E_1)\bigr)$ representing $\be$, and $b_p=0$ otherwise. 

Note that $\be$ is represented by $(b_p)$ with $b_p\in\Hom_\F\bigl(E_2\vert_p,E_1\vert_p\bigr)\ot_\F\La_\nov$, but to define a bounding cochain we require that $b_p\in\Hom_\F\bigl(E_2\vert_p,E_1\vert_p\bigr)\ot_\F\La_\nov^{\ge 0}$. This can be achieved by multiplying $\be,b_p$ by $P^\la$ for $\la\gg 0$, which does not change $(L,E,b)$ up to isomorphism in~$D^b\sF(M)$.

\begin{figure}[htb]
\centerline{$\splinetolerance{.8pt}
\begin{xy}
0;<1mm,0mm>:
,(10,0);(-20,0)**\crv{(10,5)&(0,10)}
?(.95)="a"
?(.88)="b"
?(.8)="c"
?(.72)="d"
?(.63)="e"
?(.53)="f"
?(.42)="g"
?(.27)="i"
?(.6)="y"
,(-20,0);(-30,-6)**\crv{(-30,-5)}
,(10,0);(-20,0)**\crv{(10,-5)&(0,-10)}
?(.95)="k"
?(.88)="l"
?(.8)="m"
?(.72)="n"
?(.63)="o"
?(.53)="p"
?(.42)="q"
?(.27)="s"
?(.6)="z"
,(-20,0);(-30,6)**\crv{(-30,5)}
,"a";"k"**@{.}
,"b";"l"**@{.}
,"c";"m"**@{.}
,"d";"n"**@{.}
,"e";"o"**@{.}
,"f";"p"**@{.}
,"g";"q"**@{.}
,"i";"s"**@{.}
,"y"*{<}
,"z"*{>}
,(-20,0)*{\bu}
,(-20,-3)*{q}
,(-47,0)*{\mu_{L_+,L_-}(q)=2}
,(0,0)*{\Si}
,(12,0)*{L}
,(-32.5,4)*{L_-}
,(-32.5,-4)*{L_+}
\end{xy}$}
\caption{$J$-holomorphic `teardrop' making immersed $HF^*$ obstructed}
\label{bs2fig3}
\end{figure}
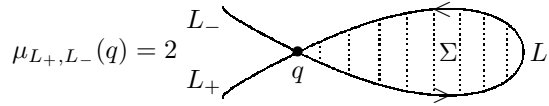

The new cause of obstructions to $HF^*$ for immersed Lagrangians $L$
with transverse self-intersections is `teardrop-shaped'
$J$-holomorphic discs $\Si$ of the form shown in Figure
\ref{bs2fig3}, with one corner at $q\in M$, and with
$\mu_{L_+,L_-}(q)=2$, where $L_\pm$ are the local sheets of $L$
intersecting at $q$. As $\mu_{L_+,L_-}(q)=2$ the moduli space of
such discs has virtual dimension 0. Such $\Si$ only obstruct $HF^*$
if they have `small area' (that is, $\area(\Si)$ is
smaller than the areas of other relevant curves with boundary in
$L$). Thus we deduce an analogue of Lemma~\ref{bs2lem1}:

\begin{lem} Suppose $(M,J,g,\Om)$ is a Calabi--Yau $m$-fold and\/
$(L,E)$ is an immersed Lagrangian brane in $M$ with only
transverse self-intersections. If\/ $b_{m-2}(L)=0$ and\/ $L$ has no
self-intersection points $p$ with\/ $\mu_{L^+,L^-}(p)=2$ or $m-2,$
where $L_\pm$ are the local sheets of\/ $L$ at\/ $p,$ then $(L,E)$ has
$HF^*$ unobstructed.
\label{bs2lem2}
\end{lem}

In \S\ref{bs25} we explained that if $(L^t,E^t):t\in[0,1]$ is a smooth family of embedded Lagrangian branes with the $L^t$ Hamiltonian isotopic and the $E^t$ locally constant in $t$, and $b^0$ is a bounding cochain for $(L^0,E^0)$, then $b^0$ extends to bounding cochains $b^t:t\in[0,1]$ for $(L^t,E^t)$ by a kind of `parallel transport', and the isomorphism class of $(L^t,E^t,b^t)$ in $D^b\sF(M)$ is independent of~$t\in[0,1]$.

In the immersed case, things are more complicated. Firstly, there
are two notions of Hamiltonian isotopy. Let $\io^t:L\ra M$ for
$t\in[0,1]$ be a smooth family of compact, immersed Lagrangians in
$M$, where we also write $\io^t:L\ra M$ as $L^t$. We call the family
{\it globally Hamiltonian isotopic\/} if $\frac{\d}{\d t}\io^t$ for
$t\in[0,1]$ is Hamiltonian flow by $H^t\ci\io^t$ for some smooth
$H^t:M\ra\R$.  We call the family {\it locally Hamiltonian
isotopic\/} if $\frac{\d}{\d t}\io^t$ for $t\in[0,1]$ is Hamiltonian
flow by some smooth $H^t:L\ra\R$, where there may exist $p_+,p_-\in
L$ with $\io^t(p_+)=\io^t(p_-)$ but $H^t(p_+)\ne H^t(p_-)$, so that $H^t$ does not descend from $L$ to~$M$.

There is a notion of `parallel transport' for bounding cochains $b^t$ along such local Hamiltonian isotopies, but it does not work all the time. Suppose for simplicity that $L^t$ has only transverse self-intersections for all $t\in[0,1]$. Then the self-intersection points of $L^t$ in $M$ depend smoothly on $t\in[0,1]$, so we can write $p^t=\io^t(p^t_+)=\io^t(p^t_-)$ for the intersection of local sheets $L^t_+\ni p^t_+$, $L^t_-\ni p^t_-$ of $L^t$ for $t\in[0,1]$, where $p^t_\pm,L^t_\pm$ depend smoothly on $t$. Then $\mu_{L^t_+,L^t_-}(p^t)$ is independent of~$t$.

Let $b^t$ be a bounding cochain for $(L^t,E^t)$ depending
smoothly on $t$, with $(L^t,E^t,b^t)\cong(L^0,E^0,b^0)$ in $D^b\sF(M)$. Then $b^t$ evolves in time by a kind of `parallel transport'. Let
$p^t,p^t_\pm,L^t_\pm$ be as above with $\mu_{L^t_+,L^t_-}(p^t)=1$. As above, $b^t$ includes an element $b^t_{p^t}\in\Hom_\F\bigl(E_+^t\vert_{p^t},E_-^t\vert_{p^t}\bigr)\ot_\F\La_\nov^{\ge 0}$. Since the local systems $E^t$ are locally constant in $t$, we can identify the fibres $E_+^t\vert_{p^t}$ for $t\in[0,1]$, and the fibres $E_-^t\vert_{p^t}$ for $t\in[0,1]$, and so regard $\Hom_\F\bigl(E_+^t\vert_{p^t},E_-^t\vert_{p^t}\bigr)\ot_\F\La_\nov^{\ge 0}$ as being independent of $t$. Then $b^t_{p^t}$ is not constant, but
evolves by
\e
\frac{\d }{\d t}b^t_{p^t}=\bigl(H^t(p^t_+)-H^t(p^t_-)\bigr)\cdot\log P\cdot b^t_{p^t}.
\label{bs2eq18}
\e
Integrating this over $[0,t]$ gives
\e
b^t_{p^t}=P^{\ts\int_0^t(H^t(p^s_+)-H^t(p^s_-)\d s}\cdot b^0_{p^0}.
\label{bs2eq19}
\e

Suppose $b^0_{p^0}\ne 0$, and write $b^0_{p^0}=\sum_{i=0}^\iy a_iP^{\la_i}$ with $a_i\in\Hom_\F(E_+^0\vert_{p^0},E_-^0\vert_{p^0})$, $a_0\ne 0$, and $0\le\la_0<\la_1<\la_2<\cdots$.
Then
\e
b^t_{p^t}=\sum_{i=0}^\iy
a_iP^{\la_i+\ts\int_0^t\bigl(H^t(p^s_+)-H^t(p^s_-)\bigr)\d
s}.
\label{bs2eq20}
\e
Thus $b^t_{p^t}\in\Hom_\F\bigl(E_+^t\vert_{p^t},E_-^t\vert_{p^t}\bigr)\ot_\F\La_\nov^{\ge 0}\subset\Hom_\F\bigl(E_+^t\vert_{p^t},E_-^t\vert_{p^t}\bigr)\ot_\F\La_\nov$, required for $b^t$ to
be a bounding cochain by \eq{bs2eq17}, if and only if
\e
\la_0+\int_0^t\bigl(H^t(p^s_+)-H^t(p^s_-)\bigr)\d s\ge
0.
\label{bs2eq21}
\e

Hence we have the following situation, which will be important in \S\ref{bs34}. Let $(L^t,E^t),$ $t\in[0,1]$ be a local Hamiltonian isotopy of Lagrangian branes in $M$, and $b^0$ a bounding cochain for $(L^0,E^0)$. We may extend $b^0$ to a family of bounding cochains $b^t:t\in[0,T]$ for $(L^t,E^t),t\in[0,T]$ for some $T\in[0,1]$, so that $(L^t,E^t,b^t)\cong(L^0,E^0,b^0)$ in $D^b\sF(M)$. But at time $t=T$ we may cross a `wall' when the l.h.s.\ of \eq{bs2eq21} becomes zero, and we cannot define $b^t$ for $t>T$. Either $(L^t,E^t)$ for $t>T$ may have $HF^*$ obstructed, or a bounding cochain $\ti b^t$ may exist but $(L^t,E^t,\ti b^t)\not\cong(L^0,E^0,b^0)$ in $D^b\sF(M)$.

The {\it Lagrangian $h$-principle}, due to Gromov \cite[p.~60-61]{Grom}
and Lees \cite{Lees}, says that two Lagrangians $L,L'$ are locally
Hamiltonian isotopic in $(M,\om)$ if and only if they are homotopic
in a weak sense, which can be well understood using homotopy theory,
and is weaker than isomorphism in $D^b\sF(M)$. So we should expect local Hamiltonian isotopies to connect Lagrangians with $HF^*$ unobstructed and with $HF^*$ obstructed, or to connect non-isomorphic Lagrangians in~$D^b\sF(M)$.

\begin{rem} As in \S\ref{bs25}, in the embedded case, the Fukaya category $\sF(M)$ has objects $(L,E,b)$ for $(L,E)$ an embedded Lagrangian brane and $b$ a bounding cochain, but the derived Fukaya category $D^b\sF(M)$ has objects {\it twisted complexes}, consisting of objects $(L_1,E_1,b_1),\ldots,(L_n,E_n,b_n)$ in $\sF(M)$ together with $b_{ij}\in CF^*\bigl((L_i,E_i),(L_j,E_j)\bigr)$ for $1\le i<j\le n$ satisfying an equation.

In the immersed case (at least, if we do not require objects of our immersed Fukaya category to have transverse self-intersections), we can regard such a twisted complex as a single object $(L,E,b)$ in $\sF(M)$, where $L$ is the disjoint union $L_1\amalg\cdots\amalg L_n$, considered as a single immersed Lagrangian, $E\vert_{L_i}=E_i$, and $b$ is a bounding cochain for $(L,E)$ built from $b_1,\ldots,b_n$ and $b_{ij}$ for $i<j$. Thus there is no need to add twisted complexes, and we can suppose all objects of $D^b\sF(M)$ are of the form~$(L,E,b)$.

The idempotent completion $D^\pi\sF(M)$ of $D^b\sF(M)$ as in \S\ref{bs25} could still include objects which are direct summands of some $(L,E,b)$, but do not have a good geometric interpretation. However, in \S\ref{bs31} we will conjecture that in the situation we are interested in, $D^b\sF(M)$ is already idempotent complete, so that we can take all objects of $D^\pi\sF(M)$ to be of the form~$(L,E,b)$.
\label{bs2rem3}
\end{rem}

\section{The conjectures}
\label{bs3}

I now explain a conjectural picture linking Bridgeland stability on
the derived Fukaya category $D^b\sF(M)$ of a Calabi--Yau manifold
$M$, special Lagrangians, Lagrangian mean curvature flow, and
obstructions to Lagrangian Floer cohomology. I had help from many
people in forming this picture, and drew inspiration from
\cite{Brid1,FOOO,GSSZ,Neve1,Neve3,Thom,ThYa}, and other places. Any mistakes are my own.

I will state some Conjectures, and also `Principles', which are too
vague to be called conjectures, but describe how I think the
mathematics ought to work. This material is intended to motivate
future research. Note that even the Conjectures are imprecise, and
may well be false in their current form.

So, for ambitious readers: few points will be awarded for {\it
disproving\/} the conjectures below, if there is some simple way to
rephrase them, retaining their spirit, but excluding the
counterexample you have in mind. Your mission, should you choose to
accept it, is to find the {\it correct\/} version of the
conjectures, and prove them; or else to show that the whole picture
is fundamentally flawed.

Be warned that I expect the difficulty of proving Conjectures
\ref{bs3conj1} and \ref{bs3conj6} increases sharply with dimension,
and even in dimension 3 is probably comparable in difficulty to the
three-dimensional Poincar\'e Conjecture, as proved by Perelman and
others (see \cite{Pere1,Pere2,Pere3} and Morgan and Tian \cite{MoTi}). The two-dimensional case may be feasible, though challenging. However, verifying that smaller parts of the picture work as expected could provide a lot of
interesting research projects.

\subsection{Bridgeland stability on $D^b\sF(M)$ for $M$
Calabi--Yau}
\label{bs31}

Let $(M,J,g,\Om)$ be a Calabi--Yau $m$-fold, with K\"ahler form $\om$, so that $(M,\om)$ is a symplectic Calabi--Yau manifold. As in \S\ref{bs25}, we will consider the {\it derived Fukaya category\/} $D^b\sF(M)$ of $M$, in the sense of Fukaya, Oh, Ohta and Ono \cite{Fuka1,FOOO}. Objects of $D^b\sF(M)$ include triples $(L,E,b)$, where $L$ is a compact, spin, graded Lagrangian in $M$ and $E\ra L$ a rank one $\F$-local system such that $(L,E)$ has $HF^*$ unobstructed, and $b$ is a bounding cochain for~$(L,E)$.

Note in particular that {\it not every compact, graded Lagrangian
$L$ or brane $(L,E)$ yields an object of\/} $D^b\sF(M)$, but only those $(L,E)$ with $HF^*$ unobstructed. One of our themes will be that we expect
Lagrangians $L$ with $HF^*$ unobstructed to be better-behaved from
the point of view of Lagrangian MCF.

We hope to use special Lagrangians and Lagrangian MCF in
$(M,J,g,\Om)$ to define an additional structure on the triangulated
category $D^b\sF(M)$, a stability condition in the sense of
Bridgeland \cite{Brid1} (see also Huybrechts~\cite{Huyb}):

\begin{dfn} Let $\cT$ be a triangulated category. A ({\it
Bridgeland\/}) {\it stability condition\/} $(Z,\cP)$ on $\cT$
consists of a group homomorphism $Z:K_0(\cT)\ra\C$ called the {\it
central charge}, and full additive subcategories
$\cP(\phi)\subset\cT$ for each $\phi\in\R$, satisfying the following
properties:
\begin{itemize}
\setlength{\parsep}{0pt}
\setlength{\itemsep}{0pt}
\item[(i)] If $A\in\cP(\phi)$ then $Z([A])=m(A)e^{i\pi\phi}$
for some $m(A)>0$.
\item[(ii)] For all $\phi\in\R$, $\cP(\phi+1)=\cP(\phi)[1]$.
\item[(iii)] If $\phi_1>\phi_2$ and $A_j\in\cP(\phi_j)$ then
$\Hom_\cT(A_1,A_2)=0$.
\item[(iv)] For each nonzero object $F\in\cT$ there is a finite
sequence of real numbers $\phi_1>\phi_2>\cdots>\phi_n$ and a
diagram in $\cT$
\end{itemize}
\e
\begin{gathered}
\xymatrix@!0@C=31.5pt@R=28pt{ 0=F_0 \ar[rr] && F_1 \ar[rr] \ar[dl]
&& F_2 \ar[r] \ar[dl] & \cdots
\ar[r] & F_{n-1} \ar[rr] && F_n=F, \ar[dl] \\
& A_1 \ar[ul]^{[1]}  && A_2 \ar[ul]^{[1]} &&&& A_n \ar[ul]^{[1]} }
\end{gathered}
\label{bs3eq1}
\e
\begin{itemize}
\setlength{\parsep}{0pt}
\setlength{\itemsep}{0pt}
\item[]where the triangles are distinguished and $A_j\in\cP(\phi_j)$
for $j=1,\ldots,n$.
\end{itemize}
Objects in $\cP(\phi)$ for some $\phi\in\R$ are called {\it
semistable}.
\label{bs3def1}
\end{dfn}

The following conjecture extending Thomas \cite{Thom} (perhaps excluding (c),(c$)'$?) is folklore, known for years in some form to many in the Geometry and String Theory communities, and is mentioned briefly in Bridgeland~\cite[\S 1.4]{Brid1}.

\begin{conj} Let\/ $(M,J,g,\Om)$ be a Calabi--Yau $m$-fold, either compact or suitably convex at infinity, and\/ $D^b\sF(M)$ the derived Fukaya category of\/ $M$ in the sense of\/ {\rm\cite{Fuka1,FOOO}}. Then there exists a natural Bridgeland stability condition $(Z,\cP)$ on $D^b\sF(M)$ such that:
\begin{itemize}
\setlength{\parsep}{0pt}
\setlength{\itemsep}{0pt}
\item[{\bf(a)}] The central charge $Z$ is the composition of
the natural maps
\e
\smash{\xymatrix@C=70pt{ K_0(D^b\sF(M)) \ar[r]^{(L,E,b)\mapsto
[L]} & H_m(M;\Z)
\ar[r]^(0.55){[L]\mapsto[\Om]\cdot[L]=\int_L\Om} & \C. }}
\label{bs3eq2}
\e
\item[{\bf(b)}] If\/ $(L,E,b)\in D^b\sF(M)$ with\/ $L$ special
Lagrangian of phase $e^{i\pi\phi},$ so that\/ $L$ has constant
phase function $\th_L=\pi\phi,$ then\/ $(L,E,b)\in\cP(\phi)$.
\item[{\bf(c)}] {\bf (Dubious, probably false as stated.)}
Suppose we enlarge the definition of\/ $D^b\sF(M)$ so that it
contains `as many Lagrangians $L$ as possible for which\/ $HF^*$
can be defined', including immersed Lagrangians as in\/ {\rm\S\ref{bs26},} and some classes of singular Lagrangians. Then every isomorphism class of objects in $\cP(\phi)$ for any $\phi\in\R$ contains a unique representative $(L,E,b)$ with\/ $L$ a (possibly immersed or
singular) special Lagrangian of phase~$e^{i\pi\phi}$.
\end{itemize}

Part\/ {\bf(c)} requires the inclusion of badly singular Lagrangians
in $D^b\sF(M),$ which may not be feasible. Here is an alternative
which may work with\/ $D^b\sF(M)$ containing only more mildly
singular Lagrangians:
\begin{itemize}
\setlength{\parsep}{0pt}
\setlength{\itemsep}{0pt}
\item[{\bf(c$\bs ){}'$}] {\bf (Still dubious.)} Suppose we
enlarge\/ $D^b\sF(M)$ so that it contains `sufficiently many
Lagrangians $L$ for which\/ $HF^*$ can be defined', including
immersed and some singular Lagrangians. Then for any $\ep>0$
and\/ $\phi\in\R,$ every isomorphism class of objects in
$\cP(\phi)$ contains a representative $(L,E,b)$ whose phase
function $\th_L$ maps $\th_L:L\ra(\pi\phi-\ep,\pi\phi+\ep)$.
\end{itemize}

\label{bs3conj1}
\end{conj}

\begin{rem}{\bf(i)} The enlargement of $D^b\sF(M)$ envisaged in
(c),(c$)'$ adds more objects to $D^b\sF(M)$, but it need not change
$D^b\sF(M)$ up to equivalence.

An example of the kind of enlargement the author has in mind is
including {\it immersed\/} Lagrangians in $D^b\sF(M)$, as in \S\ref{bs26}. We have embedded and immersed derived Fukaya categories $D^b\sF(M)_{\rm em}\subset D^b\sF(M)_{\rm im}$, but if every immersed Lagrangian $(L,E,b)$ in $D^b\sF(M)_{\rm im}$ is equivalent to a twisted complex of embedded Lagrangians, then~$D^b\sF(M)_{\rm em}\simeq D^b\sF(M)_{\rm im}$.

For many applications in symplectic topology, one only really cares
about $D^b\sF(M)$ up to equivalence, so adding extra geometric
objects to $D^b\sF(M)$ in this way is unnecessary. But for
Conjecture \ref{bs3conj1}(c),(c$)'$, it is vital --- if an isomorphism
class in $\cP(\phi)$ contains a unique special Lagrangian
representative $(L,E,b)$, and $L$ happens to be immersed, then
restricting to embedded Lagrangians would make Conjecture
\ref{bs3conj1}(c) false. (Examples of such isomorphism classes are easy to find, e.g.\ take $L=L_1\cup L_2$ for $L_1,L_2$ embedded special Lagrangians of phase $e^{i\pi\phi}$ with $L_1\cap L_2\ne\es$.) Similarly, we will see that the programme of long-time existence for Lagrangian MCF we outline below must take
place in an enlarged category of Lagrangians to have any chance of
working.

\smallskip

\noindent{\bf(ii)} The uniqueness of $(L,E,b)$ in its isomorphism
class in Conjecture \ref{bs3conj1}(c), provided it exists, should be
proved as in Thomas and Yau \cite[Th.~4.3]{ThYa}.

Note however that Thomas and Yau's method does not exclude the possibility that $L'\ra L$ and $L''\ra L$ are non-isomorphic $k$-fold multiple covers of a non-simply-connected special Lagrangian $L$ in $M$ for $k>1$, with $(L',E',b')\cong (L'',E'',b'')$ in $D^b\sF(M)$. A good uniqueness statement in Conjecture \ref{bs3conj1}(c) may be that the {\it special Lagrangian integral current\/} in Geometric Measure Theory induced by $L$ is unique, so that in the case above the special Lagrangian integral currents of both $L',L''$ would be $kL$.
\smallskip

\noindent{\bf(iii)} There may be a way to construct the expected Bridgeland stability conditions on $D^b\sF(M)$ in examples (though initially without proving that semistable objects are represented by special Lagrangians) using Mirror Symmetry.

Kontsevich's {\it Homological Mirror Symmetry Conjecture\/} \cite{Kont} roughly says that Calabi--Yau $m$-folds should exist in `mirror pairs' $M,\check M$ for which there should be an equivalence of triangulated categories
\e
D^\pi\sF(M)\simeq D^b\coh(\check M),
\label{bs3eq3}
\e
where $D^b\coh(\check M)$ is the derived category of coherent sheaves on $\check M$. (Really $\check M$ should be defined over the Novikov ring~$\La_\nov$.)

Kontsevich \cite{Kont} proved \eq{bs3eq3} when $M$ is an elliptic curve (a Calabi--Yau 1-fold). Seidel \cite{Seid2} proved it for $M$ a quartic surface in $\CP^3$ (a Calabi--Yau 2-fold), and Sheridan \cite{Sher1} proved it for $M$ a smooth Calabi--Yau $m$-fold hypersurface in $\CP^{m+1}$ for $m\ge 3$. If \eq{bs3eq3} holds then stability conditions on $D^\pi\sF(M)$ are equivalent to stability conditions on $D^b\coh(\check M)$. But derived categories of coherent sheaves are generally better understood than derived Fukaya categories.

Bridgeland stability conditions on $D^b\coh(M)$ are defined by Bridgeland \cite[Ex.~5.4]{Brid1} for $M$ a Calabi--Yau 1-fold and \cite{Brid2} for $M$ an algebraic $K3$ surface (a Calabi--Yau 2-fold). Assuming a conjecture on `Bogomolov--Gieseker type inequalities', Bayer, Macr\`\i\ and Toda \cite{BMT} construct Bridgeland stability conditions on $D^b\coh(M)$ for $M$ a Calabi--Yau 3-fold; the conjecture is proved by Macioca and Piyaratne \cite{MaPi1,MaPi2} when $M$ is an abelian 3-fold.

Combining the two, one may be able to construct examples of Bridgeland stability conditions on $D^\pi\sF(M)$ for $M$ a Calabi--Yau 1-fold, 2-fold or 3-fold.
\label{bs3rem1}
\end{rem}

The next definition and conjecture give an alternative formulation of stability which is much closer to Thomas' definition~\cite[Def.~5.1]{Thom}:

\begin{dfn} Let $(M,J,g,\Om)$ be a Calabi--Yau $m$-fold, either compact or suitably convex at infinity, and $D^b\sF(M)$ the derived Fukaya category of $M$, enlarged as in Conjecture \ref{bs3conj1} to include immersed Lagrangians, and maybe also some classes of singular Lagrangians. As in Remark \ref{bs2rem3}, we may take all objects in $D^b\sF(M)$ to be of the form $(L,E,b)$, we do not need twisted complexes. 

Suppose $\al\in\R$ is such that $[\Om]\cdot[L]\notin e^{i\pi\al}\cdot(0,\iy)$ for all $(L,E,b)$ in $D^b\sF(M)$, where $[\Om]\in H^m(M;\C)$ and $[L]\in H_m(M;\Z)$. As there are only countably many such homology classes $[L]$, this holds for generic $\al\in\R$. Write $\cA_\al$ for the full subcategory of $D^b\sF(M)$ with objects $(L,E,b)$ such that the phase function $\th_L$ of $L$ maps $L\ra(\pi\al,\pi(\al+1))$. Write $\oA_\al$ for the full subcategory of objects in $D^b\sF(M)$ isomorphic to an object of $\cA_\al$, so that $\cA_\al,\oA_\al$ are equivalent categories with~$\cA_\al\subset\oA_\al\subset D^b\sF(M)$.

We have $\cA_\al[1]\!=\!\cA_{\al+1}$ and $\oA_\al[1]\!=\!\oA_{\al+1}$.
The condition on $\al$ is to avoid taking phases in a half-open interval $(\pi\al,\pi(\al+1)]$, which could cause problems. If $(L,E,b)\in\cA_\al$, then $L$ is almost calibrated (has phase variation less than~$\pi$). 

Using the almost calibrated condition, we see that every $(L,E,b)\in\cA_\al$ has a unique {\it global phase\/} $\phi(L)\in(\pi\al,\pi(\al+1))$ with $\int_L\Om=Re^{i\phi(L)}$ for $R>0$, as in Thomas \cite[\S 3]{Thom}. If $(L',E',b')\in\oA_\al$ then $(L',E',b')\cong(L,E,b)$ in $D^b\sF(M)$ for some $(L,E,b)\in\cA_\al$, and $\int_{L'}\Om=\int_L\Om=Re^{i\phi(L)}$, where $\phi(L)$ is independent of the choice of $(L,E,b)$. Thus we may  define $\phi(L')=\phi(L)$ for~$(L',E',b')\in\oA_\al$.

In a similar way to Thomas \cite[Def.~5.1]{Thom}, we say that a nonzero object $(L,E,b)$ in $\cA_\al$ or $\oA_\al$ is {\it stable\/} (or {\it semistable\/}) if there is no distinguished triangle 
\e
\smash{\xymatrix@C=25pt{ (L_1,E_1,b_1) \ar[r] & (L,E,b) \ar[r] & (L_2,E_2,b_2)
\ar[r] & (L_1,E_1,b_1)[1] }}
\label{bs3eq4}
\e
in $D^b\sF(M)$ with $(L_1,E_1,b_1),(L_2,E_2,b_2)$ nonzero objects in $\cA_\al$ or $\oA_\al$ such that $\phi(L_1)\ge\phi(L_2)$ (or~$\phi(L_1)>\phi(L_2)$).
\label{bs3def2}
\end{dfn}

\begin{conj} In Definition\/ {\rm\ref{bs3def2},} $\oA_\al$ is the heart of a bounded t-structure on $D^b\sF(M),$ and so $\cA_\al,\oA_\al$ are abelian categories, and\/ \eq{bs3eq4} becomes a short exact sequence in $\cA_\al$ or $\oA_\al$. Furthermore, the Bridgeland stability condition $(Z,\cP)$ on $D^b\sF(M)$ in Conjecture\/ {\rm\ref{bs3conj1}} may be described as follows: $Z$ is defined by {\rm\eq{bs3eq2},} and\/ $\cP(\al)=\es,$ and for each\/ $\be\in(\al,\al+1),$ $\cP(\be)$ is the full subcategory of semistable objects $(L,E,b)$ in $\oA_\al$ with\/~$\phi(L)=\pi\be$.
\label{bs3conj2}
\end{conj}

Note that (semi)stability in Definition \ref{bs3def2} is equivalent to slope (semi)\-stab\-il\-ity on the (conjecturally abelian) categories $\cA_\al,\oA_\al$, with slope function
\begin{equation*}
\mu(L,E,b)=\frac{-\cos\pi\al \int_{L_{\vphantom{l}}}\Re\Om-\sin\pi\al \int_L\Im\Om}{-\sin\pi\al \int_L\Re\Om+\cos\pi\al \int_L\Im\Om}\,,
\end{equation*}
since $\phi(L)=\tan^{-1}(\mu(L,E,b))+\pi\al+\frac{\pi}{2}$. Thomas' analogue of \eq{bs3eq4} is to require $L_1,L_2$ to intersect transversely at one point $p$, and $L$ to be Hamiltonian isotopic to the Lagrangian connect sum $L_1\# L_2$ at $p$. Equation \eq{bs3eq4} is more general, e.g.\ it does not imply that $L$ is diffeomorphic to $L_1\# L_2$. It would be nice to state the relationship between $L$ and $L_1,L_2$ geometrically rather than categorically.

As in \S\ref{bs25}, there are two versions $D^b\sF(M)\subseteq D^\pi\sF(M)$ of the derived Fukaya category, where $D^b\sF(M)$ has objects twisted complexes in $\sF(M)$, and $D^\pi\sF(M)$ has objects direct summands of objects in $D^b\sF(M)$. By Remark \ref{bs2rem3}, for immersed Lagrangians we do not need to add twisted complexes, so we can take all objects in $D^b\sF(M)$ to be of the form~$(L,E,b)$.

We wrote Conjecture \ref{bs3conj1} using $D^b\sF(M)$, since the extra objects in $D^\pi\sF(M)$ are not geometric, and our programme does not make sense for them. For example, the map $K_0(D^b\sF(M))\ra H_m(M;\Z)$ in \eq{bs3eq2} is not defined for $D^\pi\sF(M)$, as we cannot associate a homology class to a direct summand of~$(L,E,b)$.

However, if $D^b\sF(M)$ has a Bridgeland stability condition, then it has a bounded t-structure, and so by Huybrechts \cite[Rem.~1.15]{Huyb} it is idempotent complete. Thus Conjecture \ref{bs3conj1} or Conjecture \ref{bs3conj2} imply:

\begin{conj} In the situation of Conjecture\/ {\rm\ref{bs3conj1},} the enlarged version of\/ $D^b\sF(M)$ with objects $(L,E,b)$ for $L$ a possibly singular, compact, immersed, graded Lagrangian is idempotent complete. Hence $D^\pi\sF(M)\simeq D^b\sF(M),$ and we can take all objects of\/ $D^\pi\sF(M)$ to be geometric, of the form~$(L,E,b)$.
\label{bs3conj3}
\end{conj}

\begin{rem} A partial verification of Conjecture \ref{bs3conj3} in the case $M=T^2$ is provided by Haug \cite{Haug}. He defines a version of the derived Fukaya category $D^b\sF(T^2)$ in which the objects are twisted complexes built out of pairs $(L,E)$ for $L$ a compact, spin, graded, embedded Lagrangian in $T^2$, and $E\ra L$ a local system, and proves that $D^b\sF(T^2)$ is idempotent complete. 

Haug remarks \cite[\S 1]{Haug} that for $T^2$, including local systems $E\ra L$ has the effect of making $D^b\sF(T^2)$ idempotent complete, and that $D^b\sF(T^2)$ would not be idempotent complete if we took objects to be twisted complexes of Lagrangians $L$ rather than pairs $(L,E)$. This shows that {\it including local systems $E\ra L$ in objects $(L,E,b)$ is necessary for our programme}, since otherwise Conjecture \ref{bs3conj3} and hence Conjecture \ref{bs3conj1} would be false even for $M=T^2$. We will see in \S\ref{bs34} how nontrivial local systems are needed for some kinds of surgeries.

Haug's definition of $D^b\sF(T^2)$ is not quite the same as ours. He does not include bounding cochains $b$ in his objects $(L,E)$ (the simplicity of dimension 1 permits this). He fixes $\F=\C$. His local systems $E\ra L$ \cite[\S 3.1.1]{Haug} are not $\F$-local systems, as in \S\ref{bs25}, but $\La_\nov$-local systems of arbitrary finite rank, such that (roughly) the eigenvalues of $\Hol(\nabla_E)$ lie in $\F^*\subset\La_\nov^*$ to leading order.

I expect this should be related to our definition of $D^b\sF(T^2)$ as follows. In dimension 1, the combination of a rank one $\F$-local system $E\ra L$ and a bounding cochain $b$ is essentially equivalent to a rank one $\La_\nov$-local system $E_\nov\ra L$ satisfying Haug's condition, where the holonomies satisfy $\Hol(\nabla_{E_\nov})[\ga]=\Hol(\nabla_E)[\ga]\cdot e^{\int_\ga b}$ for $[\ga]\in\pi_1(L)$. Also, I expect that for $T^2$, considering rank one local systems $E\ra L$ on immersed Lagrangians has a similar effect to considering higher rank local systems $E\ra L$ on embedded Lagrangians.
\label{bs3rem2}
\end{rem}

\subsection{Approaching Conjecture \ref{bs3conj1} using Lagrangian
MCF}
\label{bs32}

Here is our suggestion for a programme to prove Conjecture
\ref{bs3conj1} using Lagrangian MCF, building on Thomas and Yau \cite{ThYa}. We will state a conjecture about it in \S\ref{bs39}, after discussing issues that arise in the programme in~\S\ref{bs33}--\S\ref{bs38}.
\smallskip

\noindent{\bf Programme for (partially?) proving Conjecture
\ref{bs3conj1} using LMCF.} {\it Let\/ $(M,J,g,\Om)$ be a Calabi--Yau $m$-fold, either compact or suitably convex at infinity, and suppose as in Conjecture\/ {\rm\ref{bs3conj1}} that we have extended the definition of\/
$D^b\sF(M)$ to include immersed Lagrangians, as in {\rm\cite{AkJo},}
and some classes of singular Lagrangians.

Define $Z:K_0(D^b\sF(M))\ra\C$ by {\rm\eq{bs3eq2},} and define
$\cP(\phi)$ for $\phi\in\R$ to be the full subcategory of objects
$A$ in $D^b\sF(M)$ isomorphic to $(L,E,b)$ for $L$ a (possibly
singular) special Lagrangian of phase $e^{i\pi\phi}$ with\/
$\th_L=\pi\phi,$ as in Conjecture\/ {\rm\ref{bs3conj1}(c),} or
alternatively those objects $A$ in $D^b\sF(M)$ which for any $\ep>0$
are isomorphic to some $(L,E,b)$ with phase function
$\th_L:L\ra(\pi\phi-\ep,\pi\phi+\ep),$ as in
Conjecture\/~{\rm\ref{bs3conj1}(c$)'$}.

We must prove $(Z,\cP)$ is a Bridgeland stability condition on
$D^b\sF(M)$. We discuss only the problem of verifying Definition\/
{\rm\ref{bs3def1}(iv)} for objects $F=(L,E,b),$ where $(L,E)$ is a nonsingular, immersed Lagrangian brane with\/ $HF^*$ unobstructed. For such $(L,E,b),$ we must construct a diagram
\e
\begin{gathered}
\xymatrix@!0@C=32.5pt@R=28pt{ 0=F_0 \ar[rr] && F_1 \ar[rr] \ar[dl] &&
F_2 \ar[r] \ar[dl] & \cdots
\ar[r] & F_{n-1} \ar[rr] && F_n\!=\!(L,E,b), \ar[dl] \\
& (L_1,E_1,b_1) \ar[ul]^(0.6){[1]}  && (L_2,E_2,b_2) \ar[ul]^(0.6){[1]} &&&&
(L_n,E_n,b_n) \ar[ul]^(0.6){[1]} }
\end{gathered}
\label{bs3eq5}
\e
in $D^b\sF(M),$ where $L_1,\ldots,L_n$ are either unique (possibly
singular) special Lagrangians with\/ $\th_{L_j}=\pi\phi_j$ for\/
$\phi_1>\phi_2>\cdots>\phi_n,$ or else (possibly singular)
Lagrangians with\/ $\th_{L_j}:L_j\ra(\pi\phi_j-\ep,\pi\phi_j+\ep)$ for
arbitrarily small\/~$\ep>0$.

We aim to construct a unique family
$\bigl\{(L^t,E^t,b^t):t\in[0,\iy)\bigr\}$ satisfying:
\begin{itemize}
\setlength{\parsep}{0pt}
\setlength{\itemsep}{0pt}
\item[{\bf(a)}] $(L^0,E^0,b^0)=(L,E,b)$.
\item[{\bf(b)}] There is a (hopefully finite) series of
\begin{bfseries}singular times\end{bfseries} $0<T_1<T_2<\cdots,$
such that if\/ $t\in [0,\iy)\sm\{T_1,T_2,\ldots\}$ then
$(L^t,E^t,b^t)$ is an object in $D^b\sF(M)$ isomorphic to $(L,E,b),$
with\/ $L^t$ a (possibly immersed or singular) compact, graded
Lagrangian in $M,$ with\/ $HF^*$ unobstructed.
\item[{\bf(c)}] The family
$\bigl\{L^t:t\in[0,\iy)\sm\{T_1,T_2,\ldots\}\bigr\}$ satisfies
Lagrangian mean curvature flow, and $\bigl\{E^t:t\in[0,\iy)\sm\{T_1,T_2,\ldots\}\bigr\}$ are locally constant in $t$. (As a shorthand, we will say that \begin{bfseries}the family of Lagrangian branes $\bigl\{(L^t,E^t):t\in[0,\iy)\sm\{T_1,T_2,\ldots\}\bigr\}$ satisfies Lagrangian MCF\end{bfseries}.) The bounding cochains $b^t$ also change by a kind of `parallel transport' for
$t\in[0,\iy)\sm\{T_1,T_2,\ldots\}$ as in\/ {\rm\S\ref{bs25}--\S\ref{bs26},} to ensure that the isomorphism class of\/ $(L^t,E^t,b^t)$ in $D^b\sF(M)$ remains constant.
\item[{\bf(d)}] Let\/ $T_i$ for $i=1,2,\ldots$ be a singular
time and $\ep>0$ be small, so that\/
$\bigl\{L^t:t\in(T_i-\ep,T_i)\bigr\}$ and\/
$\bigl\{L^t:t\in(T_i,T_i+\ep)\bigr\}$ satisfy Lagrangian MCF. As
$t\ra T_i$ in $(T_i-\ep,T_i),$ the flow usually undergoes a finite
time singularity of Lagrangian MCF. But see {\rm\S\ref{bs34}} for a case in which the limit is smooth as $t\ra T_i$ in $(T_i-\ep,T_i),$ and singular as $t\ra T_i$ in~$(T_i,T_i+\ep)$.

We do not require $(L^{T_i},E^{T_i},b^{T_i})$ to be an object in
$D^b\sF(M),$ as the singularities of\/ $L^{T_i}$ may be too bad,
and if so, $b^{T_i}$ is meaningless.

The topologies of\/ $L^t$ for\/ $t\in(T_i-\ep,T_i),$ and\/
$L^{T_i},$ and\/ $L^t$ for\/ $t\in(T_i,T_i+\ep),$ may all be
different, so we may think of the (possibly singular) manifolds
$L^t$ as undergoing a surgery at time $t=T_i$. Nonetheless, the
family $\bigl\{L^t:t\in(T_i-\ep,T_i+\ep)\bigr\}$ is in a
suitable sense continuous, for instance, as graded Lagrangian integral
currents in $M$ in Geometric Measure Theory.
\item[{\bf(e)}] For the case of Conjecture\/
{\rm\ref{bs3conj1}(c),} we have
$\lim_{t\ra\iy}L^t=L_1\cup\cdots\cup L_n,$ where $L_j$ is a
(possibly badly singular) special Lagrangian with phase
$e^{i\pi\phi_j}$ and phase function $\th_{L_j}=\pi\phi_j,$ for
$\phi_1>\phi_2>\cdots>\phi_n$. The local systems $E_1,\ldots,E_n,$ bounding cochains $b_1,\ldots,b_n$ and morphisms in \eq{bs3eq5} are obtained from $\lim_{r\ra\iy}E^t$ and\/~$\lim_{t\ra\iy}b^t$.

For the case of Conjecture\/ {\rm\ref{bs3conj1}(c$)'$,} if\/
$t\gg 0$ then there is a decomposition
$L^t=L^t_1\amalg\cdots\amalg L^t_n,$ such that\/ $\th_{L^t}$
maps $L^t_j\ra (\pi\phi_j-\ep^t,\pi\phi_j+\ep^t)$ for
$j=1,\ldots,n,$ where $\ep^t>0$ with $\ep^t\ra 0$ as $t\ra\iy$.
\end{itemize}}

\begin{rem}{\bf(i)} In dimension $m>1$, Lagrangian MCF $\{L^t:t\in[0,\iy)\}$ starting from a compact, embedded Lagrangian $L^0$ can flow to immersed
Lagrangians $L^t$ in finite time, as sketched in Figure \ref{bs3fig1}, or vice versa. (When $m=1$, embedded curves remain embedded.)

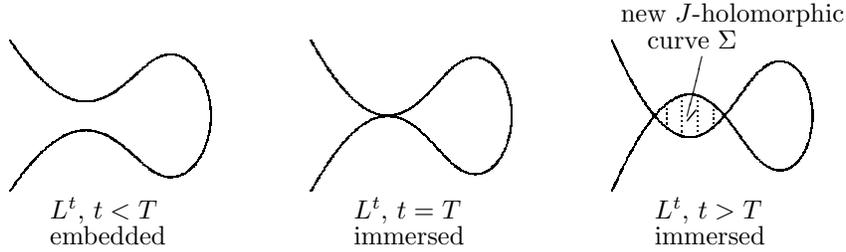
\begin{figure}[htb]
\centerline{$\splinetolerance{.8pt}
\begin{xy}
0;<1mm,0mm>:
,(-40,10);(-40,-10)**\crv{(-30,-5)&(-20,15)&(-10,0)&(-20,-15)&(-30,5)}
,(0,10);(0,-10)**\crv{(10,-8.15)&(20,15)&(30,0)&(20,-15)&(10,8.15)}
,(40,10);(40,-10)**\crv{(50,-13)&(60,15)&(70,0)&(60,-15)&(50,13)}
?(.07)="a"
?(.093)="b"
?(.12)="c"
?(.15)="d"
?(.93)="e"
?(.907)="f"
?(.88)="g"
?(.85)="h"
,"a";"e"**@{.}
,"b";"f"**@{.}
,"c";"g"**@{.}
,"d";"h"**@{.}
,(-27,-14)*{\begin{subarray}{l}\ts\text{$L^t,$ $t<T$} \\ \ts \text{embedded}\end{subarray}}
,(13,-14)*{\begin{subarray}{l}\ts\text{$L^t$, $t=T$} \\ \ts \text{immersed}\end{subarray}}
,(53,-14)*{\begin{subarray}{l}\ts\text{$L^t,$ $t>T$} \\ \ts \text{immersed}\end{subarray}}
,(50,0);(52,8)**\crv{}
,(50,0);(51.2,1.5)**\crv{}
,(50,0);(49.7,1.6)**\crv{}
,(56,12)*{\begin{subarray}{l}\ts\text{new $J$-holomorphic} \\
\ts \text{\quad curve $\Si$}\end{subarray}}\end{xy}$}
\caption{LMCF flowing from embedded to immersed Lagrangians}
\label{bs3fig1}
\end{figure}

Therefore, to carry out the programme above, we must include {\it
immersed\/} Lagrangians in $D^b\sF(M)$, since otherwise in the
situation of Figure \ref{bs3fig1} we could not continue the
programme past $t=T$. This inclusion was discussed in \S\ref{bs26},
using the extension of \cite{FOOO} to immersed Lagrangians in Akaho
and Joyce~\cite{AkJo}.

Observe that for Lagrangian MCF $\{L^t:t\in[0,T)\}$ of immersed,
graded Lagrangians $L^t$ in a Calabi--Yau $m$-fold, the $L^t$ for
$t\in[0,T)$ are all {\it locally Hamiltonian isotopic\/} in the
sense of \S\ref{bs26}, but not necessarily globally Hamiltonian
isotopic, as in Figure~\ref{bs3fig1}.

Thus, for immersed Lagrangian MCF we must deal with the possibility
that even without finite time singularities, the flow may take us
from Lagrangians with unobstructed $HF^*$ to Lagrangians with
obstructed $HF^*$, or change the isomorphism class in $D^b\sF(M)$,
since we explained in \S\ref{bs26} that local Hamiltonian isotopies
can do this. We discuss this further in~\S\ref{bs34}.
\smallskip

\noindent{\bf(ii)} Notice the strong similarity of the programme
above with the proof of the three-dimensional Poincar\'e Conjecture
by Perelman, Hamilton and others, as in \cite{MoTi,Pere1,Pere2,Pere3}.
There one starts with a Riemannian 3-manifold $(M,g)$ (the analogue
of Lagrangians), and applies rescaled Ricci flow, encountering
finite time singularities at times $0<T_1<T_2<\cdots$ when one does
surgery, until as $t\ra\iy$ the flow converges to a disjoint union
of constant curvature Riemannian 3-manifolds (the analogue of
special Lagrangians).

In dimension $m=3$, I expect the programme above to be of comparable
difficulty to the Poincar\'e Conjecture. As the dimension increases,
so should the difficulty, as there will be more kinds of finite-time
singularities to worry about.

\smallskip

\noindent{\bf(iii)} There is a literature on mean curvature flow with surgeries for (suitably convex) hypersurfaces in $\R^n$ which resembles our proposal --- see for instance Huisken and Sinestrari \cite{HuSi} and Brendle and Huisken \cite{BrHu}.
\smallskip

\noindent{\bf(iv)} As for isolated conical singularities of special
Lagrangians \cite[\S 3]{Joyc2}, one could try to define an `index'
$\ind(\tau)$ for different `types' $\tau$ of finite time
singularities of Lagrangian MCF, which measures the codimension in
the infinite-dimensional family $\scr L$ of Lagrangians $L$ in $M$
in which singularities of type $\tau$ occur in Lagrangian MCF
starting from $L$. So for instance, Lagrangian MCF starting from a
generic Lagrangian $L$ could only develop singularities
with~$\ind(\tau)=0$.

We could modify the programme above by taking $L^0$ to be a generic
Hamiltonian perturbation of $L$ in (a), rather than $L^0=L$. Then
the Lagrangian MCF singularities occurring at the singular times
$T_1,T_2,\ldots$ would have to have index 0. This might have the
effect of limiting the kinds of singular Lagrangians that must be
included in $D^b\sF(M)$ to make the programme work.

For similar ideas in MCF of hypersurfaces in $\R^n$, see Angenent and Vel\'azquez \cite{AnVe} who construct examples of non-generic finite time singularities of MCF, and Colding and Minicozzi \cite{CoMi}, who classify the possible finite time singularities of MCF starting from a generic, compact, embedded surface $\Si^2$ in $\R^3$.
\smallskip

\noindent{\bf(v)} Taking limits $\lim_{t\ra\iy}L^t$ in (e) above is
likely to introduce different, and worse, singularities than those
in the finite time singularities $L^{T_1},L^{T_2},\ldots.$ Also, I
expect $\lim_{t\ra\iy}L^t$ to be unchanged by Hamiltonian
perturbations of $L^0$, so taking $L^0$ generic as in {\bf(iv)} will not
help.

It seems likely that the possible singularities occurring in
$\lim_{t\ra\iy}L^t$ may be too severe to incorporate as objects in
$D^b\sF(M)$. Thus, although Conjecture \ref{bs3conj1}(c) is more
attractive, Conjecture \ref{bs3conj1}(c$)'$ is more plausible.
\smallskip

\noindent{\bf(vi)} Since $\{L^t:t\in[0,\iy)\}$ above satisfies Lagrangian MCF, one might expect that $\{L^t:t\in[0,\iy)\}$ depends only on $L^0=L$, and is independent of $E,b$ in $(L,E,b)$. However, in \S\ref{bs34} we will describe a surgery `opening a neck' depending on $E,b$, so $\{L^t:t\in[0,\iy)\}$ does depend on all of $L,E,b$, not just on~$L$.
\smallskip

\noindent{\bf(vii)} Behrndt \cite{Behr1} defines a modification of
Lagrangian MCF which works in {\it almost Calabi--Yau manifolds\/}
$(M,J,g,\Om)$, that is, a complex $m$-manifold $(M,J)$ with K\"ahler
metric $g$ and nonvanishing holomorphic $(m,0)$-form $\Om$ which
need not satisfy \eq{bs2eq1}, so that $g$ need not be Ricci-flat. I
expect the whole of this paper also to work for modified Lagrangian MCF in almost Calabi--Yau $m$-folds.

\label{bs3rem3}
\end{rem}

\subsection{On finite time singularities of Lagrangian MCF}
\label{bs33}

Finite time singularities of Lagrangian MCF were discussed in \S\ref{bs23}. For {\it graded\/} Lagrangian MCF, Theorem \ref{bs2thm4} says that any finite time singularity must be of type II, and Theorem \ref{bs2thm5} that any finite time singularity must admit a `type II blow up' modelled on a nontrivial eternal solution of Lagrangian MCF in $\C^m$. As in the end of \S\ref{bs23}, two natural classes of eternal solutions are provided by SL $m$-folds in $\C^m$, and Lagrangian MCF translators.

Motivated by this, the next `principle' gives heuristic pictures of how the author expects two different classes of finite time singularities to work.

\begin{princ} Let\/ $(M,J,g,\Om)$ be a compact Calabi--Yau $m$-fold
and\/ $\{L^t:t\in[0,T)\}$ a family of compact, immersed, graded
Lagrangians in $M$ satisfying Lagrangian MCF, with a finite time
singularity at\/ $t=T,$ and a singular point at\/ $x\in M$. Here are
broad descriptions of two classes of such singularities:
\begin{itemize}
\setlength{\parsep}{0pt}
\setlength{\itemsep}{0pt}
\item[{\bf(a)}] Let\/ $U$ be a small open neighbourhood of\/ $x$ in
$M,$ which we identify with a small open neighbourhood of\/ $0$
in $\C^m=T_xM,$ and\/ $\ep>0$ be small. Then $L^t\cap U$
approximates a closed, exact SL\/ $m$-fold in $\C^m$
for~$t\in(T-\ep,T)$.

Since SL\/ $m$-folds are stationary points of LMCF, to `first
order'\/ $L^t\cap U$ is constant in $t,$ but to `second order'
$L^t\cap U$ wanders slowly in the moduli space of closed, exact
SL\/ $m$-folds in $\C^m,$ until at time\/ $t=T$ it hits a
singular SL\/ $m$-fold. This `wandering' is driven by `outside
influences' from the whole of\/ $L^t,$ not just from\/ $L^t\cap
U$.

For example, if\/ $N$ is an exact asymptotically conical SL\/
$m$-fold in $\C^m,$ we could have $L^t\cap U \approx f(t)\cdot
N$ for $t\in(T-\ep,T),$ where $f:(T-\ep,T)\ra(0,\iy)$ is smooth
with\/ $f(t)\ra 0$ as~$t\ra T$.
\item[{\bf(b)}] Let\/ $U,\ep$ be as in {\bf(a)}. Then $L^t\cap U$
approximates a closed, exact LMCF translator in
$\C^m=T_xM$ for $t\in(T-\ep,T)$. To `first order'\/ $L^t\cap U$ moves by translation in $\C^m=T_xM,$ since it approximates a translating soliton. But to second order it also wanders slowly in the moduli space of
closed, exact LMCF translators in $\C^m,$ driven by `outside influences' from the whole of\/ $L^t,$ until at time $t=T$ it hits a singular soliton.

For example, if\/ $N$ is an exact LMCF translator in $\C^m$ with translating vector $v\in\C^m,$ we could have $L^t\cap U \approx f(t)\cdot N+g(t)\cdot v$ for $t\in(T-\ep,T),$ where $f,g:(T-\ep,T)\ra(0,\iy)$ are smooth with\/ $f(t)\ra 0$ as~$t\ra
T$.
\end{itemize}
\label{bs3princ1}
\end{princ}

\begin{rem}{\bf(i)} We will describe examples of behaviours (a),(b) in \S\ref{bs35} and \S\ref{bs38}. Section \ref{bs37} discusses a class of singularities not of type (a) or~(b).

Note that in (a),(b) we do {\it not\/} simply mean that the
singularity has a type II blow up $\{\ti L^s:s\in\R\}$ in Theorem
\ref{bs2thm5} with $\ti L^s$ special Lagrangian or an LMCF
translator. In general type II blow ups describe {\it only a small
part\/} of the singularity, and may give little idea of the global
geometry and topology near the singular point. The point of (a),(b)
is that in these cases we have a more complete picture of the
singularity than a general type II blow up gives.
\smallskip

\noindent{\bf(ii)} As in \S\ref{bs23}, Lagrangian MCF shrinkers do not occur
in the graded case. The other major class of Lagrangian MCF solitons, Lagrangian MCF expanders (as in \S\ref{bs23}) are not relevant to the formation of singularities of the flow (that is, to describing the flow
immediately before the singular time $t=T_i$). However, we can use
Lagrangian MCF expanders to model the flow immediately {\it after\/}
a surgery at a singular time $t=T_i$, and we do this in~\S\ref{bs34}.
\label{bs3rem4}
\end{rem}

If we believe Principle \ref{bs3princ1}, stretching credulity
a little further gives:

\begin{princ} Any type of (sufficiently well-behaved) singularity
of SL\/ $m$-folds, which can appear as a limit of nonsingular,
locally exact SL $m$-folds, may provide a local model for finite
time singularities of Lagrangian MCF.

Similarly, any (sufficiently well-behaved) singular Lagrangian in
$\C^m$ which can appear as a limit of nonsingular, exact Lagrangian
MCF translators in $\C^m,$ may provide a local model for finite time
singularities of Lagrangian MCF.
\label{bs3princ2}
\end{princ}

This suggests a class of research problems:

\begin{prob}{\bf(a)} Choose from the literature your favourite
family of explicit, nonsingular, exact SL\/ $m$-folds $N_s$ in
$\C^m$ which converge to an explicit singular SL\/ $m$-fold\/ $N_0$
as $s\ra 0$. For example, let $N$ be an exact AC SL\/ $m$-fold in
$\C^m$ with cone $C,$ and take $N_s=s\cdot N$ for $s>0$
and\/~$N_0=C$.

Construct examples $\{L^t:t\in [0,T]\}$ of Lagrangian MCF in $\C^m$
or in a Calabi--Yau $m$-fold\/ $(M,J,g,\Om)$ with finite time
singularities at\/ $t=T$ for which\/ $L^T$ has a singularity at\/ $x\in\C^m$ modelled on\/ $N_0,$ and\/ $L^t$ near\/ $x$ for\/ $t\in(T-\ep,T)$ approximates\/ $N_{s(t)},$ where\/ $s(t)\ra 0$ as\/ $t\ra T,$ as in Principle\/~{\rm\ref{bs3princ1}(a)}.

\smallskip

\noindent{\bf(b)} If you can do {\bf(a)\rm,} determine whether
Lagrangian MCF starting from a small generic Hamiltonian
perturbation of\/ $L^0$ also develops finite time singularities of
the same type. In this case, we call this type a
\begin{bfseries}generic singularity\end{bfseries} of Lagrangian MCF.
If it is not generic, compute the expected codimension amongst
Hamiltonian perturbations of\/ $L^0$ in which singularities of this
type occur.
\smallskip

\noindent{\bf(c)} Repeat\/ {\bf(a)\rm,\bf(b)} for LMCF translators rather than SL\/ $m$-folds.
\label{bs3prob1}
\end{prob}

\subsection{Flowing from unobstructed to obstructed immersed Lagrangians}
\label{bs34}

In Remark \ref{bs3rem3}(i) we noted that Lagrangian MCF may take an immersed Lagrangian brane $(L^t,E^t)$ with $HF^*$ unobstructed to one $(L^{t'},E^{t'})$ for $t'>t$ with $HF^*$ obstructed, without finite time singularities. This is a problem for the programme of \S\ref{bs32}, as we need $(L^t,E^t)$ to have $HF^*$ unobstructed for all $t$. We now discuss this problem in more detail, and explain how to solve it.

Let $(M,J,g,\Om)$ be a Calabi--Yau $m$-fold and $\{(L^t,E^t):t\in [0,T)\}$ a family of Lagrangian branes satisfying Lagrangian MCF. Suppose, for simplicity, that all the $L^t$ have transverse self-intersections. Then the self-intersection points of $L^t$ in $M$ depend smoothly on $t\in[0,T)$, so we can write $p^t$ for the intersection of local sheets $L^t_+,L^t_-$ of $L^t$ for $t\in[0,T)$, where $p^t,L^t_\pm$ depend smoothly on $t$. Then $\mu_{L^t_+,L^t_-}(p^t)$ is independent of~$t$.

Suppose that $b^t$ is a bounding cochain for $L^t$ depending
smoothly on $t$, with $(L^t,E^t,b^t)\cong(L^0,E^0,b^0)$ in $D^b\sF(M)$. Then $b^t$ evolves in time by a kind of `parallel transport'. Let
$p^t,L^t_\pm$ be as above with $\mu_{L^t_+,L^t_-}(p^t)=1$. Then as
in \S\ref{bs26}, $b^t$ includes an element $b^t_{p^t}\in\Hom_\F\bigl(E_+^t\vert_{p^t},E_-^t\vert_{p^t}\bigr)\ot_\F\La_\nov^{\ge 0}$. The analysis of \eq{bs2eq18}--\eq{bs2eq21} holds, with $H^t=-\th_{L^t}$. Thus, writing $b^0_{p^0}=\sum_{i=0}^\iy
a_iP^{\la_i}$ with $a_0\ne 0$ and $0\le\la_0<\la_1<\la_2<\cdots$, we have
\begin{equation*}
b^t_{p^t}=\sum_{i=0}^\iy
a_iP^{\la_i+\ts\int_0^t\bigl(\th_{L^s_-}(p^s)-\th_{L^s_+}(p^s)\bigr)\d
s},
\end{equation*}
and $b^t_{p^t}\in\Hom_\F\bigl(E_+^t\vert_{p^t},E_-^t\vert_{p^t}\bigr)\ot_\F\La_\nov^{\ge 0}\subset\Hom_\F\bigl(E_+^t\vert_{p^t},E_-^t\vert_{p^t}\bigr)\ot_\F\La_\nov$, required for $b^t$ to be a bounding cochain, if and only if
\e
\la_0+\int_0^t\bigl(\th_{L^s_-}(p^s)-\th_{L^s_+}(p^s)\bigr)\d s\ge
0.
\label{bs3eq6}
\e

We can now explain how Lagrangian MCF can flow from $HF^*$
unobstructed to $HF^*$ obstructed: as $t$ increases, we can cross a
`wall' at $t=T_1$ when the l.h.s.\ of \eq{bs3eq6} becomes
negative, so that $b^t_{p^t}\notin\Hom_\F\bigl(E_+^t\vert_{p^t},E_-^t\vert_{p^t}\bigr)\ot_\F\La_\nov^{\ge 0}$ for $t>T_1$.
Then $b^t$ is not a bounding cochain, and $(L^t,E^t)$ may have $HF^*$
obstructed.

\begin{figure}[htb]
\centerline{$\splinetolerance{.8pt}
\begin{xy}
0;<1mm,0mm>:
,(10,0);(-20,0)**\crv{(10,5)&(0,10)}
?(.95)="aa"
?(.88)="bb"
?(.8)="cc"
?(.72)="dd"
?(.63)="ee"
?(.53)="ff"
?(.42)="gg"
?(.27)="ii"
?(.6)="yy"
,(10,0);(-20,0)**\crv{(10,-5)&(0,-10)}
?(.95)="kk"
?(.88)="ll"
?(.8)="mm"
?(.72)="nn"
?(.63)="oo"
?(.53)="pp"
?(.42)="qq"
?(.27)="ss"
?(.6)="zz"
,"aa";"kk"**@{.}
,"bb";"ll"**@{.}
,"cc";"mm"**@{.}
,"dd";"nn"**@{.}
,"ee";"oo"**@{.}
,"ff";"pp"**@{.}
,"gg";"qq"**@{.}
,"ii";"ss"**@{.}
,"yy"*{<}
,"zz"*{>}
,(-60,0);(-20,0)**\crv{(-40,10)}
?(.95)="a"
?(.85)="b"
?(.75)="c"
?(.65)="d"
?(.55)="e"
?(.45)="f"
?(.35)="g"
?(.25)="h"
?(.15)="i"
?(.05)="j"
?(.5)="y"
,(-60,0);(-70,-6)**\crv{(-70,-5)}
,(-60,0);(-20,0)**\crv{(-40,-10)}
?(.95)="k"
?(.85)="l"
?(.75)="m"
?(.65)="n"
?(.55)="o"
?(.45)="p"
?(.35)="q"
?(.25)="r"
?(.15)="s"
?(.05)="t"
?(.5)="z"
,(-60,0);(-70,6)**\crv{(-70,5)}
,"a";"k"**@{.}
,"b";"l"**@{.}
,"c";"m"**@{.}
,"d";"n"**@{.}
,"e";"o"**@{.}
,"f";"p"**@{.}
,"g";"q"**@{.}
,"h";"r"**@{.}
,"i";"s"**@{.}
,"j";"t"**@{.}
,"y"*{<}
,"z"*{>}
,(-60,0)*{\bu}
,(-60,-3)*{p^t}
,(-20,0)*{\bu}
,(-20,-3)*{q^t}
,(-40,0)*{\Si^t_1}
,(-72,4)*{L^t_+}
,(-72,-4)*{L^t_-}
,(-22,8)*{\mu_{L_+^t,L_-^t}(q^t)\!=\!2}
,(-57,8)*{\mu_{L_+^t,L_-^t}(p^t)\!=\!1}
,(0,0)*{\Si^t_2}
,(13,0)*{L^t}
,(-40,8)*{L^t_-}
,(-40,-8)*{L^t_+}
,(-5,8)*{L^t_+}
,(-5,-8)*{L^t_-}
\end{xy}$}
\caption{Crossing between $HF^*$ unobstructed when $\area(\Si_1^t)
<\area(\Si_2^t)$ and $HF^*$ obstructed when $\area(\Si_1^t)>\area(\Si_2^t)$}
\label{bs3fig2}
\end{figure}
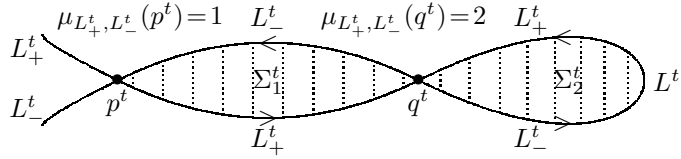

To make this more explicit, let us simplify further, and suppose
that $L^t$ has only two self-intersection points $p^t,q^t$ with
$\mu_{L^t_+,L^t_-}(p^t)=1$ and $\mu_{L^t_+,L^t_-}(q^t)=2$, and there
are only two $J$-holomorphic curves $\Si^t_1,\Si^t_2$ with boundary
in $L^t$ which are relevant to obstructions to $HF^*$, which are as
shown in Figure \ref{bs3fig2}, so that $\Si^t_1$ has two corners at
$p^t,q^t$ and $\Si^t_2$ one corner at $q^t$. Note that $\Si^t_2$ is
the type of curve in Figure \ref{bs2fig3} that can cause
obstructions to immersed~$HF^*$.

Then $(L^t,E^t)$ has $HF^*$ unobstructed if and only if
$\area(\Si^t_2)\ge\area(\Si^t_1)$, and if so, the bounding cochain
$b^t$ has
\e
b^t_{p^t}=a_0 P^{\,\area(\Si^t_2)-\area(\Si^t_1)}+\text{higher order
terms,}
\label{bs3eq7}
\e
where $0\ne a_0\in\Hom_\F\bigl(E_+^t\vert_{p^t},E_-^t\vert_{p^t}\bigr)$. We can think of $\Si^t_2-\Si^t_1$ as a `virtual $J$-holomorphic curve' with `virtual area' $\area(\Si^t_2)-\area(\Si^t_1)$ and one corner at $p^t$, which obstructs $HF^*$ if this virtual area is negative.

Under Lagrangian MCF we have
\ea
\frac{\d}{\d t}&\bigl(\area(\Si^t_2)-\area(\Si^t_1)\bigr)=
-\int_{\pd\Si^t_2}\d\th_{L^t}+\int_{\pd\Si^t_1}\d\th_{L^t}=
-\bigl[\th_{L_+^t}(q^t)-\th_{L_-^t}(q^t)\bigr]
\nonumber\\
&+\bigl[\th_{L_+^t}(q^t)-\th_{L_-^t}(q^t)+\th_{L_-^t}(p^t)
-\th_{L_+^t}(p^t)\bigr]=\th_{L_-^t}(p^t)-\th_{L_+^t}(p^t).
\label{bs3eq8}
\ea

Suppose now that the family $\{(L^t,E^t):t\in [0,T)\}$ passes from $HF^*$
unobstructed when $t<T_1$ to $HF^*$ obstructed when $t>T_1$. Then
$\area(\Si^t_2)-\area(\Si^t_1)$ crosses zero at $t=T_1$ going from
positive to negative, so \eq{bs3eq8} shows that
\e
\th_{L_-^{T_1}}(p^{T_1})-\th_{L_+^{T_1}}(p^{T_1})\le 0.
\label{bs3eq9}
\e

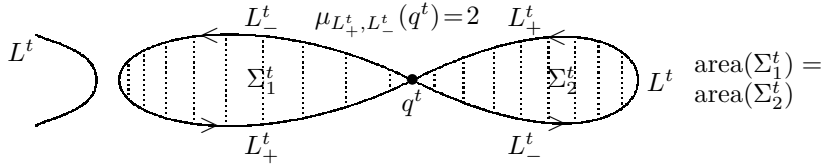
\begin{figure}[htb]
\centerline{$\splinetolerance{.8pt}
\begin{xy}
0;<1mm,0mm>:
,(10,0);(-20,0)**\crv{(10,5)&(0,10)}
?(.95)="aa"
?(.88)="bb"
?(.8)="cc"
?(.72)="dd"
?(.63)="ee"
?(.53)="ff"
?(.42)="gg"
?(.27)="ii"
?(.6)="yy"
,(10,0);(-20,0)**\crv{(10,-5)&(0,-10)}
?(.95)="kk"
?(.88)="ll"
?(.8)="mm"
?(.72)="nn"
?(.63)="oo"
?(.53)="pp"
?(.42)="qq"
?(.27)="ss"
?(.6)="zz"
,"aa";"kk"**@{.}
,"bb";"ll"**@{.}
,"cc";"mm"**@{.}
,"dd";"nn"**@{.}
,"ee";"oo"**@{.}
,"ff";"pp"**@{.}
,"gg";"qq"**@{.}
,"ii";"ss"**@{.}
,"yy"*{<}
,"zz"*{>}
,(-59,0);(-20,0)**\crv{(-59,6)&(-40,10)}
?(.95)="a"
?(.85)="b"
?(.75)="c"
?(.65)="d"
?(.55)="e"
?(.45)="f"
?(.35)="g"
?(.25)="h"
?(.15)="i"
?(.05)="j"
?(.5)="y"
,(-59,0);(-20,0)**\crv{(-59,-6)&(-40,-10)}
?(.95)="k"
?(.85)="l"
?(.75)="m"
?(.65)="n"
?(.55)="o"
?(.45)="p"
?(.35)="q"
?(.25)="r"
?(.15)="s"
?(.05)="t"
?(.5)="z"
,(-62,0);(-70,6)**\crv{(-62,4)&(-68,5)}
,(-62,0);(-70,-6)**\crv{(-62,-4)&(-68,-5)}
,"a";"k"**@{.}
,"b";"l"**@{.}
,"c";"m"**@{.}
,"d";"n"**@{.}
,"e";"o"**@{.}
,"f";"p"**@{.}
,"g";"q"**@{.}
,"h";"r"**@{.}
,"i";"s"**@{.}
,"j";"t"**@{.}
,"y"*{<}
,"z"*{>}
,(-20,0)*{\bu}
,(-20,-3)*{q^t}
,(-40,0)*{\Si^t_1}
,(-72,4)*{L^t}
,(-20,0)*{\bu}
,(-22,8)*{\mu_{L_+^t,L_-^t}(q^t)\!=\!2}
,(0,0)*{\Si^t_2}
,(13,0)*{L^t}
,(-40,8)*{L^t_-}
,(-40,-9)*{L^t_+}
,(-5,8)*{L^t_+}
,(-5,-9)*{L^t_-}
,(26,0)*{\begin{subarray}{l}\ts \area(\Si^t_1)=\\
\ts\area(\Si^t_2)\end{subarray}}
\end{xy}$}
\caption{$L^t$ for $t>T_1$, after Lagrangian connected sum surgery at $p^t$}
\label{bs3fig3}
\end{figure}

We claim that in the programme of \S\ref{bs32}, the correct thing to
do is to change $L^t$ for $t\ge T_1$ by doing a surgery at $p^t$
when $t=T_1$, a Lagrangian connected sum of the two sheets
$L^t_+,L^t_-$ at $p^t$, so that $L^t$ for $T_1<t<T_1+\ep$ looks
roughly like Figure \ref{bs3fig3}. We will call this surgery
`opening a neck'. The self-intersection $p^t$ is
now gone, and there are two $J$-holomorphic discs $\Si^t_1,\Si^t_2$
with one corner at $q^t$. Since we do the surgery when
$\area(\Si^t_1)=\area(\Si^t_2)$, we have $\area(\Si^t_1)=
\area(\Si^t_2)$ for all $t>T_1$, though $\Si^t_1,\Si^t_2$ are in
different relative homology classes. As their areas are equal, the
obstructions from $\Si^t_1,\Si^t_2$ cancel for suitable $E^t$, and $(L^t,E^t)$ for $t>T_1$
has $HF^*$ unobstructed.

We have $\mu_{L_+^{T_1},L_-^{T_1}}(p^{T_1})=1$ and
$\th_{L_+^{T_1}}(p^{T_1})\ge\th_{L_-^{T_1}}(p^{T_1})$ by
\eq{bs3eq9}. Suppose strict inequality holds,
$\th_{L_+^{T_1}}(p^{T_1})>\th_{L_-^{T_1}}(p^{T_1})$. Then from
Definition \ref{bs2def7}, we see that there is an identification
$T_{p^{T_1}}M\cong\C^m$ identifying
$J\vert_{p^{T_1}},g\vert_{p^{T_1}}$ with the standard versions on
$\C^m$, and identifying $T_{p^{T_1}}L_+^{T_1},T_{p^{T_1}}L_-^{T_1}$
with the Lagrangian planes $\Pi_0,\Pi_{\bs\phi}$ in \eq{bs2eq12} for
$\phi_1,\ldots,\phi_m\in(0,\pi)$ with $0<\phi_1+\cdots+\phi_m<\pi$,
where $\phi_1+\cdots+\phi_m<\pi$ comes from
$\mu_{L_+^{T_1},L_-^{T_1}} (p^{T_1})=1$
and~$\th_{L_+^{T_1}}(p^{T_1})>\th_{L_-^{T_1}}(p^{T_1})$.

Thus, by Example \ref{bs2ex4} there is a unique, exact
Joyce--Lee--Tsui Lagrangian MCF expander $L_{\bs\phi}^1$ with
$\al=1$ in $T_{p^{T_1}}M$ asymptotic to $T_{p^{T_1}}L_+^{T_1}\cup
T_{p^{T_1}}L_-^{T_1}$, and Theorem \ref{bs2thm6} shows that
$L_{\bs\phi}^1$ is the only LMCF expander with $\al=1$ in
$T_{p^{T_1}}M$ asymptotic to $T_{p^{T_1}}L_+^{T_1}\cup
T_{p^{T_1}}L_-^{T_1}$. Note that $\sqrt{2(t-T_1)}\cdot
L_{\bs\phi}^1$ for $t>T_1$ satisfy Lagrangian MCF in
$T_{p^{T_1}}M\cong\C^m$. We now aim to define the $L^t$ for
$T_1<t<T_1+\ep$ by gluing in $\sqrt{2(t-T_1)}\cdot L_{\bs\phi}^1$
into $L^{T_1}$ near~$p^{T_1}$.

To define the local systems $E^t$ for $t>T_1$, note that $b^{T_1}_{p^{T_1}}=a_0 +\cdots$ by \eq{bs3eq7}, where $0\ne a_0\in\Hom_\F\bigl(E_+^{T_1}\vert_{p^{T_1}},E_-^{T_1}\vert_{p^{T_1}}\bigr)$. As $E^{T_1}$ has rank one, $a_0\ne 0$ implies that $a_0$ is an isomorphism. For $T_1<t<T_1+\ep$, we define $E^t$ to be equal to $E^{T_1}$ away from the `neck' region joining $L_+^{T_1}$ with $L_-^{T_1}$, and on the `neck' region we use the isomorphism $a_0$ to identify $E^{T_1}\vert_{L^+}$ and~$E^{T_1}\vert_{L^-}$. This choice of $E^t$ is necessary for the obstructions to $HF^*$ for $(L^t,E^t)$ from $\Si_1^t,\Si_2^t$ to cancel.

The bounding cochain $b^t$ for $T_1<t<T_1+\ep$ should be roughly equal to $b^{T_1}$ away from the `neck' region. On the `neck' region, $b^t$ should somehow encode the higher order terms in $b^{T_1}_{p^{T_1}}=a_0 +\cdots$, possibly in the form $b^t\approx\log\bigl(a_0^{-1}\ci b^{T_1}_{p^{T_1}}\bigr)\cdot [\cS^{m-1}_t]$, where $[\cS^{m-1}_t]\in C_{m-1}(L^t,\Z)$ is a fundamental cycle for the new small $(m\!-\!1)$-sphere $\cS^{m-1}_t$ spanning the `neck' in~$L^t$.

\begin{rem} We can now see an important reason why our programme requires the inclusion of the rank one $\F$-local systems $E\ra L$ in the objects $(L,E,b)$ of $D^b\sF(M)$, as mentioned in Remark \ref{bs3rem2}. We can also justify our definition of Lagrangian branes in Definition~\ref{bs2def6}. 

Firstly, note that if the initial local systems $E^t$ for $t<T_1$ above are trivial, the local systems $E^t$ for $t>T_1$ may not be trivial, as across the `neck' region $E^t$ for $t>T_1$ has holonomy $a_0\in\Hom_\F\bigl(E_+^{T_1}\vert_{p^{T_1}},E_-^{T_1}\vert_{p^{T_1}}\bigr)\cong\F$, and we need not have $a_0=1$. So this surgery can pass from trivial to nontrivial local systems $E^t$. If we omitted local systems $E$ in $D^b\sF(M)$, then the data $a_0$ in $b^{T_1}$ would be lost under the surgery, and $L^t$ for $t>T_1$ might have $HF^*$ obstructed.

Secondly, we take $\F$ to be a field (rather than say a commutative ring) so that $0\ne a_0\in\F$ implies that $a_0$ is an isomorphism.

Thirdly, observe that the argument above would not work for higher rank local systems $E\ra L$, which is why we restrict to rank one. If $E^{T_1}$ has different ranks $n_\pm$ on $L_\pm^{T_1}$, then it cannot extend across the `neck' to make $E^t$ for $t>T_1$. If $E^{T_1}$ has the same rank $n>1$ on $L_+^{T_1},L_-^{T_1}$, then $a_0\ne 0$ no longer implies that $a_0$ is an isomorphism, so we cannot use $a_0$ to extend $E^{T_1}$ across the `neck'.
\label{bs3rem5}
\end{rem}

Our discussion has shown the following rather neat:
\smallskip

\noindent{\bf Evidence for the viability of the programme of
\S\ref{bs32}.} {\it Let\/ $(M,J,g,\Om)$ be a Calabi--Yau $m$-fold
and\/ $\{(L^t,E^t):t\in[0,T)\}$ be a family of Lagrangian branes in $M$ satisfying Lagrangian MCF.

Suppose that\/ $(L^t,E^t)$ has $HF^*$ unobstructed for $0\le t<T_1<T,$ but at\/ $t=T_1$ crosses a `wall' into $HF^*$ obstructed, because at a
transverse self-intersection point\/ $p$ of\/ $L^{T_1}$ with\/
$\mu_{L_+^{T_1},L_-^{T_1}}(p)=1,$ the data\/ $b_p^t$ in the bounding cochain $b^t$ leaves\/ $\Hom_\F\bigl(E_+^t\vert_p,E_-^t\vert_p\bigr)\ot_\F\La_\nov^{\ge 0}$ in $\Hom_\F\bigl(E_+^t\vert_p,E_-^t\vert_p\bigr)\ot_\F\La_\nov$ when~$t=T_1$.

Then (at least if strict inequality holds in \eq{bs3eq9}) there is
a unique Lagrangian MCF expander in $T_pM$ asymptotic to
$T_pL_+^{T_1}\cup T_pL_-^{T_1},$ which we can (conjecturally) use to
do a surgery at\/ $t=T_1$ so that the flow can continue for $t>T_1$
with\/ $HF^*$ unobstructed, as in {\rm\S\ref{bs32}}. The analogue
does \begin{bfseries}not\end{bfseries} hold for flowing from $HF^*$
obstructed to $HF^*$ unobstructed.}
\smallskip

It also suggests a research project:

\begin{prob} Suppose $(M,J,g,\Om)$ is a Calabi--Yau $m$-fold, $L$ a compact, immersed Lagrangian in $M$ with a transverse self-intersection point at\/ $p\in M$ with local sheets $L_\pm,$ and\/ $N$ a Joyce--Lee--Tsui Lagrangian MCF expander in $T_pM$ asymptotic to $T_pL_+\cup T_pL_-$ and satisfying $H=F^\perp$. Prove that for small\/ $\ep>0,$ there is a unique family $\{L^t:t\in(0,\ep)\}$ of compact, immersed Lagrangians in $M$ satisfying Lagrangian MCF, such that\/ $\lim_{t\ra 0}L^t=L^0$ in a suitable sense, and for small\/ $t$ we have $L^t\approx \sqrt{2t}\cdot N$ near $p$ and\/ $L^t\approx L+tH_L$ away from~$p$.
\label{bs3prob2}
\end{prob}

The author has been told that Kim Moore and Tom Begley, students of Jason Lotay and Felix Schulze, have nearly completed a proof of Problem~\ref{bs3prob2}.

Ilmanen, Neves and Schulze \cite{INS} study the evolution of `networks' in the plane   (roughly, finite graphs smoothly embedded in $\R^2$) under mean curvature flow. In a similar way to Problem \ref{bs3prob2}, they show \cite[\S 7]{INS} that one can glue a certain kind of self-expander in at a singular point when $t=0$, and obtain existence of MCF for networks for $t\in[0,\ep)$. It seems likely that the methods of \cite{INS} will be helpful for Problem~\ref{bs3prob2}.

In the next example we use `opening necks' to resolve an apparent counterexample to our programme.

\begin{ex} Let $(M,J,g,\Om)$ be a Calabi--Yau $m$-fold, and $(L_1,E_1),(L_2,E_2)$ be embedded, transversely-intersecting, special Lagrangian branes in $M$ with phases $e^{i\pi\phi_1},e^{i\pi\phi_2}$ for $\phi_1<\phi_2$, with $HF^*$ unobstructed. Choose bounding cochains $b_1,b_2$ for $(L_1,E_1),(L_2,E_2)$. Let $0\!\ne\!\be\!\in\! HF^1\bigl((L_2,E_2,b_2),(L_1,E_1,b_1)\bigr)$, and $(\be_p)\in CF^1\bigl((L_2,E_2),(L_1,E_1)\bigr)$ represent $\be$, where for all $p\in L_1\cap L_2$ with $\mu_{L_2,L_1}(p)=1$ we have~$\be_p\in\Hom_\F\bigl(E_2\vert_p,E_1\vert_p\bigr)\ot_\F\La_\nov$ . 

Suppose $\be_p\in\Hom_\F\bigl(E_2\vert_p,E_1\vert_p\bigr)\ot_\F\La_\nov^{\ge 0}$ for all $p$. Set $(L,E)=(L_1,E_1)\cup (L_2,E_2)$, considered as an immersed Lagrangian brane in $M$. Then using the notation of \S\ref{bs26}, $b=b_1\op b_2\op(\be_p)$ is a bounding cochain for $(L,E)$, where $b_{\rm ch}=b_1\op b_2$ in $C_{m-1}(L,\La_\nov^+)=C_{m-1}(L_1,\La_\nov^+)\op C_{m-1}(L_2,\La_\nov^+)$, and the data $b_p$ for each $p\in M$ at which two local sheets $L_+,L_-$ of $L$ intersect transversely with $\mu_{L_+,L_-}(p)=1$ are $b_p=\be_p$ if $L_+=L_2$, $L_-=L_1$, and $b_p=0$ otherwise. We now have a distinguished triangle in the derived Fukaya category $D^b\sF(M)$ of immersed Lagrangians
\e
\smash{\xymatrix@C=27pt{ (L_1,E_1,b_1) \ar[r] & (L,E,b) \ar[r] & (L_2,E_2,b_2)
\ar[r]^(0.45){\be} & (L_1,E_1,b_1)[1]. }}
\label{bs3eq10}
\e

Let us apply the programme of \S\ref{bs32} to $(L,E,b)$. Since $L$ is a union of special Lagrangians of different phases, it is stationary under immersed Lagrangian MCF, so the obvious answer is that $(L^t,E^t,b^t)=(L,E,b)$ for all $t\in[0,\iy)$. Equation \eq{bs3eq10} gives a diagram for $(L,E,b)$ of the form \eq{bs3eq5} with $n=2$
\begin{equation*}
\xymatrix@!0@C=50pt@R=28pt{ 0=F_0 \ar[rr] && F_1=(L_1,E_1,b_1) \ar[rr] \ar[dl] && F_2=(L,E,b). \ar[dl] \\
& (L_1,E_1,b_1) \ar[ul]^(0.6){[1]}  && (L_2,E_2,b_2) \ar[ul]^(0.6){[1]}_(0.35)\be }
\end{equation*}
However, in \S\ref{bs32} we want such a diagram with $\phi_1>\phi_2$, but we assume that $\phi_1<\phi_2$. So writing $(L^t,E^t,b^t)=(L,E,b)$ for all $t\in[0,\iy)$ {\it does not satisfy\/} the programme of \S\ref{bs32}, as although we have long-time existence of Lagrangian MCF, the limiting behaviour at infinity is wrong, and this looks like a counterexample.

Here is the explanation. Although (at least initially) the $L^t,E^t$ are independent of $t$, the bounding cochains $b^t$ do evolve in time. Suppose $p\in L_1\cap L_2$ with $\mu_{L_2,L_1}(p)=1$. Then \eq{bs2eq18}--\eq{bs2eq21} with $H_{L_j}=-\th_{L_j}=-\pi\phi_j$ for $j=1,2$ shows that the data $b_p^t$ in $b^t$ should evolve according to the equation
\begin{equation*}
\frac{\d}{\d t}b_p^t=\pi(\phi_1-\phi_2)\cdot\log P\cdot b_p^t,
\end{equation*}
so as $b_p^0=\be_p$, the solution is $b_p^t=P^{\pi(\phi_1-\phi_2)t}\cdot\be_p$. Thus, we have 
\e
(L^t,E^t)=(L_1,E_1)\amalg (L_2,E_2), \quad
b^t=b_1\op b_2\op(P^{\pi(\phi_1-\phi_2)t}\cdot\be_p),
\label{bs3eq11}
\e
at least for small $t$. Write $\be_p=a_p P^{\la_p}+\cdots$ if $\be_p\ne 0$, where $0\ne a_p\in\Hom_\F\bigl(E_2\vert_p,E_1\vert_p\bigr)$ and $\la_p\ge 0$, and set $\la_p=\iy$ if $\be_p=0$. Then $b_p^t=a_pP^{\la_p+\pi(\phi_1-\phi_2)t}+\cdots$, so $b_p^t\in\Hom_\F\bigl(E_2^t\vert_p,E_1^t\vert_p\bigr)\ot_\F\La_\nov^{\ge 0}$ if $t\in[0,\la_p/\pi(\phi_2-\phi_1)]$.

Thus, at time $T=(\min_p\la_p)/\pi(\phi_2-\phi_1)$, the flow crosses a `wall' after which $b^t$ in \eq{bs3eq11} is no longer a bounding cochain, as $b_p^t$ leaves $\Hom_\F\bigl(E_2^t\vert_p,E_1^t\vert_p\bigr)\ot_\F\La_\nov^{\ge 0}$ for some $p$. We claim that the right thing to do is to `open a neck' at time $t=T$ at each $p$ with $\la_p$ minimal, gluing in a Joyce--Lee--Tsui LMCF expander. Then $L^t,E^t$ will undergo some nontrivial evolution for~$t>T$.

To see that a suitable LMCF expander exists to glue in at $p$, note that $\th_{L_j}(p)=\pi\phi_j$ for $j=1,2$, so $\th_{L_1}(p)<\th_{L_2}(p)$ by assumption, and as $\mu_{L_2,L_1}(p)=1$, the first equation of \eq{bs2eq13} gives $\th_{L_2}(p)<\th_{L_1}(p)+\pi$. These are the conditions for the existence of an LMCF expander in $T_pM$ asymptotic to~$T_pL_1\cup T_pL_2$.

\label{bs3ex1}
\end{ex}

\subsection{`Neck pinches' using Lawlor necks}
\label{bs35}

The programme of \S\ref{bs32} requires a flow $\{L^t:t\in[0,\iy)\}$
starting from a single Lagrangian $L^0=L$, but converging as
$t\ra\iy$ to a union $L_1\cup\cdots\cup L_n$ of several (possibly
intersecting) special Lagrangians of different phases, where we
regard $L_1\cup\cdots\cup L_n$ as a single immersed Lagrangian.
Thus, we need a local model for how one Lagrangian $L$ can break up
into a union $L_1\cup L_2$ of two Lagrangians under the flow, at
some singular time $t=T_i$, in the notation of~\S\ref{bs32}.

We call this local model a `neck pinch', as it involves the Lawlor
necks $L_{\bs\phi,A}$ of Example \ref{bs2ex1} as $A\ra 0$, so that
the `neck' pinches to a point. It is an example of Principles
\ref{bs3princ1}(b) and \ref{bs3princ2}, where the special Lagrangian
local models are the Lawlor necks $L_{\bs\phi,A}$. The possibility
of such pinching behaviour is clear from Thomas and Yau \cite{ThYa},
and Neves \cite[\S 4]{Neve1} proves that it occurs in an example,
where both \cite{ThYa,Neve1} work with $\SO(m)$-equivariant
Lagrangians, so that Lagrangian MCF is reduced to understanding
evolution of real curves.

\begin{conj} The following behaviour, which we call a `neck pinch',
can occur in Lagrangian MCF with surgeries in Calabi--Yau $m$-folds for $m\ge 2,$ as in\/ {\rm\S\ref{bs32}}. Furthermore, `neck pinches' are a \begin{bfseries}generic singularity\end{bfseries}. That is, if Lagrangian MCF beginning from $L^0$ develops a neck pinch, then Lagrangian MCF beginning from any sufficiently small Hamiltonian perturbation $\ti L^0$ of\/ $L^0$ also develops a neck pinch.

Let\/ $(M,J,g,\Om)$ be a Calabi--Yau $m$-fold, and extend\/
$D^b\sF(M)$ to include immersed Lagrangians, as in\/
{\rm\cite{AkJo}}. Suppose $\{(L^t,E^t):t\in(T-\ep,T+\ep)\}$ for $\ep>0$
small is a family of immersed Lagrangian branes in $M$ with\/ $HF^*$ unobstructed, and\/ $\{b^t:t\in(T-\ep,T+\ep)\}$ a corresponding family of bounding cochains, satisfying the following conditions:
\begin{itemize}
\setlength{\parsep}{0pt}
\setlength{\itemsep}{0pt}
\item[{\bf(i)}] The $(L^t,E^t,b^t)$ for $t\in(T-\ep,T+\ep)$ are all
isomorphic in $D^b\sF(M)$.
\item[{\bf(ii)}] When $t<T,$ $L^t,E^t$ depend smoothly on
$t\in(T-\ep,T),$ and\/ $\{(L^t,E^t):t\in(T-\ep,T)\}$ satisfies
Lagrangian MCF, with a finite time singularity at\/ $t=T,$ with
one singular point\/~$p\in M$.

Similarly, when $t\ge T,$ $L^t,E^t$ depend smoothly on
$t\in[T,T+\ep),$ and\/ $\{(L^t,E^t):t\in[T,T+\ep)\}$ satisfies
Lagrangian MCF. The topology of\/ $L^t$ for $t\in(T-\ep,T+\ep)$
changes discontinuously at\/ $t=T$. Nonetheless, the family
$\{L^t:t\in(T-\ep,T+\ep)\}$ is continuous at\/ $t=T$ in a
suitable sense, e.g.\ as graded Lagrangian integral currents in
Geometric Measure Theory.
\item[{\bf(iii)}] Identifying $M$ near $p$ with\/
$T_pM\cong\C^m$ near $0,$ for each\/ $t\in(T-\ep,T),$ $L^t$
approximates a `Lawlor neck' $L_{\bs\phi(t),A(t)}$ from
Example\/ {\rm\ref{bs2ex1},} after a translation and a\/ $\U(m)$
rotation in $\C^m$. Here $A(t)>0$ is small and $A(t)\ra 0$ as
$t\ra T,$ so that\/ $L_{\bs\phi(t),A(t)}$ converges to a union
$\Pi_0\cup\Pi_{\bs\phi(T)}$ of transversely intersecting special
Lagrangian planes in $\C^m$ as $t\ra T$.
\item[{\bf(iv)}] For $t\in[T,T+\ep),$ there is a
self-intersection point\/ $p^t$ of\/ $L^t$ where two local
sheets $L^t_\pm$ of\/ $L^t$ intersect transversely with\/
$\mu_{L^t_+,L^t_-}(p^t)=1$. Here $p^t,L^t_\pm$ depend smoothly
on $t\in[T,T+\ep),$ with\/~$p^T=p$.
\item[{\bf(v)}] We have $\th_{L^T_+}(p^T)=\th_{L^t_-}(p^T),$ and\/
$\th_{L^t_+}(p^t)<\th_{L^t_-}(p^t)$ for $t\in(T,T+\ep)$.
\item[{\bf(vi)}] The $\F$-local systems $E^t$ for $t\in[T,T+\ep)$ are constructed from the $\F$-local systems $E^{t'}$ for $t'\in(T-\ep,T)$ by deleting the `neck' in $L^{t'}$ and extending $E^{t'}$ over $p^t$ in $L^t_\pm$ in the unique possible way (at least for $m\ge 3$).
\item[{\bf(vii)}] When $t\in[T,T+\ep),$ the bounding cochain
$b^t$ for $L^t$ includes an element\/ $b^t_{p^t}\in\Hom_\F\bigl(E_+^t\vert_{p^t},E_-^t\vert_{p^t}\bigr)\ot_\F\La_\nov^{\ge 0}$
as in {\rm\S\ref{bs26}}. This is of the form
\begin{equation*}
b^t_{p^t}=a_0 P^{\la(t)}+\text{higher order terms,}
\end{equation*}
where $a_0\in\Hom_\F\bigl(E_+^t\vert_{p^t},E_-^t\vert_{p^t}\bigr)$ is the natural isomorphism induced from $E^{t'}$ for $t'\in(T-\ep,T)$ using {\bf(vi)\rm,} and\/ $\la(t)=\int_T^t\bigl(\th_{L^s_-}(p^s)-\th_{L^s_+}(p^s)\bigr)\d s,$ so that\/ $\la(T)=0$ and\/ $\la(t)>0$ for $t\in(T,T+\ep)$ by {\bf(v)}.
\end{itemize}
\label{bs3conj4}
\end{conj}

\begin{rem}{\bf(a)} The `neck pinching' behaviour of Conjecture
\ref{bs3conj4} is inverse to the `opening a neck' behaviour of
\S\ref{bs34}. So, for example, we can imagine a flow
$\{L^t:t\in[0,\iy)\}$ satisfying the programme of \S\ref{bs32}, with
two singular times $0<T_1<T_2$, which starts with a single $L^t$ for
$0\le t<T_1$, undergoes a `neck pinch' at $t=T_1$ and becomes a
union $L^t=L^t_1\cup L^t_2$ of Lagrangians $L^t_1,L^t_2$
intersecting at one point $p^t$ for $T_1<t<T_2$, and then at $t=T_2$
`opens the neck' at $p^t$ and turns back into a single Lagrangian
$L^t$ for~$t>T_2$.

Note that these inverse singular behaviours involve {\it
different\/} (though related) geometric local models, Lawlor necks
$L_{\bs\phi,A}$ and Joyce--Lee--Tsui expanders $L_{\bs\phi}^\al$. We
do not just na\"\i vely run the local picture for the flow in
reverse. Note too that `neck pinching' works only for $m\ge 2$, whereas `opening necks' works for $m\ge 1$, so when $m=1$, `opening necks' has no inverse behaviour.

In a similar way, the author expects that many types of finite time singularity possible in the programme of \S\ref{bs32} should have a
corresponding inverse type, so that changes in the topology of
$L^t$, and other qualitative features, are reversible. An exception to this is that when $m=1$, the flow can only decrease the number of self-intersection points, making the curve `less immersed'.
\smallskip

\noindent{\bf(b)} Theorem \ref{bs2thm3} shows that Lawlor necks
$L_{\bs\phi,A}$ are the only possible geometric local models for
such `neck pinches'.
\smallskip

\noindent{\bf(c)} The inequality $\th_{L^t_+}(p^t)<\th_{L^t_-}(p^t)$
in (v) is the opposite of \eq{bs3eq9} in \S\ref{bs34}.
Heuristically, we expect `small necks' to shrink under Lagrangian
MCF when $\th_{L^t_+}(p^t)<\th_{L^t_-}(p^t)$, and to grow when
$\th_{L^t_+}(p^t)>\th_{L^t_-}(p^t)$.
\smallskip

\noindent{\bf(d)} The case $m=2$ in Conjecture \ref{bs3conj4} is
special. For $m\ge 3$, the family $\cF$ of AC special Lagrangian
`Lawlor necks' $L$ in $\C^m$ asymptotic to $\Pi_0\cup\Pi_{\bs\phi}$
is (isomorphic to) $(0,\iy)$, and all such $L$ are exact. When
$m=2$, the family $\cF$ is $\R^2\sm\{0\}$, and the subfamily
$\cF_{\rm exact}$ of exact $L$ is $\R\sm\{0\}\subset\R^2\sm\{0\}$,
since then $\cF_{\rm exact}$ contains both the $L_{\bs\phi,A}$ for
$A>0$ and $\ti L_{\bs\phi,A}$ for $A<0$ in Example \ref{bs2ex1}. 

Also, when $m=2$ the local systems $E^{t'}$ for $t'\in(T-\ep,T)$ could have nontrivial holonomy around the `neck'. If so, the definition of $E^t$ for $t\in[T,T+\ep)$ in part (vi) no longer makes sense, since we cannot extend $E^{t'}$ over $p^t$ in~$L^t_\pm$.

One conclusion is that for $m=2$, though neck pinches should be
generic under Hamiltonian perturbations, they may be nongeneric (and
of index 1) under Lagrangian perturbations, since Lagrangian
perturbations may allow the flow to wander in $\cF=\R^2\sm\{0\}$
rather than $\cF_{\rm exact}=\R\sm\{0\}$, and will only hit the
singularity $0\in\R^2$ in real codimension 1 amongst initial
Lagrangians.

We can also ask: if Lagrangian MCF $\{L^t:t\in(T-\ep,T)\}$ develops a singularity as $t\ra T$ modelled on Lawlor necks $L_{\bs\phi,A}$ for $A\in (0,\iy)\subset\cF_{\rm exact}=\R\sm\{0\}$, rather than continuing for $t>T$ using immersed SL 2-folds as in Conjecture \ref{bs3conj4}, why not continue using Lawlor necks $\ti L_{\bs\phi,A}$ for $A\in(-\iy,0)\subset\cF_{\rm exact}=\R\sm\{0\}$, immediately opening the neck again, in a similar way to~\S\ref{bs34}? 

The author expects that this is the correct thing to do if $E^{t'}$ for $t'\in(T-\ep,T)$ has nontrivial holonomy around the `neck'. But in the trivial holonomy case, it would change the isomorphism class of $(L^t,E^t,b^t)$ in $D^b\sF(M)$, and so should be avoided according to the philosophy of~\S\ref{bs32}.
\label{bs3rem6}
\end{rem}

\subsection[Including singular Lagrangians in $D^b\sF(M)$; LMCF for
\\ Lagrangians with stable conical singularities]{Including singular
Lagrangians in $D^b\sF(M)$; LMCF for Lagrangians with stable conical
singularities}
\label{bs36}

The programme of \S\ref{bs32} involves flows $\{L^t:t\in[0,\iy)\}$
with the $L^t$ immersed Lagrangians which can be singular at the
singular times $t=T_1,T_2,\ldots,$ where we do not require
$(L^{T_i},E^{T_i},b^{T_i})$ to be objects of $D^b\sF(M)$. In this section we argue that in dimension $m\ge 3$, we must also allow the $L^t$ to
have certain kinds of `stable' singularities for $t\ne T_i$. To
complete the programme, Lagrangian MCF must work for such singular
Lagrangians, and we must include them as objects in the derived
Fukaya category~$D^b\sF(M)$.

In \cite{Joyc1,Joyc2,Joyc3,Joyc4,Joyc5} the author studied compact SL
$m$-folds $L$ with {\it isolated conical singularities\/} in a
Calabi--Yau $m$-fold $M$. That is, $L$ has singularities
$p_1,\ldots,p_k$ locally modelled on closed special Lagrangian cones
$C_1,\ldots,C_k$ in $\C^m$ which have isolated singularities at
$0\in\C^m$. As in \cite{Joyc2}, the deformation theory of $L$
involves an obstruction space $\O=\O_1\op\cdots\op\O_k$ which is the
sum of contributions $\O_i$ from each singular point $p_i$,
depending only on the cone $C_i$. We call the singularities $p_i$
and the SL cones $C_i$ {\it stable\/} \cite[Def.~3.6]{Joyc2} if the
obstruction spaces $\O_i$ are zero. By \cite[Cor.~6.11]{Joyc2}, if
$L$ has only stable isolated conical singularities, then the moduli
space $\cM_L$ of SL deformations of $L$ is a smooth manifold.

Few examples of stable SL cones are known. The SL $T^2$-cone $C$ in
$\C^3$ in equation \eq{bs2eq4} of Example \ref{bs2ex2} was shown to
be stable in \cite[\S 3.2]{Joyc1}. Ohnita \cite{Ohni} found four
more examples of stable SL cones in dimensions 5, 8, 14, and 26. In
dimension $m=2$, any irreducible, immersed SL cone in $\C^2$ is a
Lagrangian plane $\R^2$, or a finite cover of $\R^2$ branched at 0.
Nontrivial branched covers of $\R^2$ are unstable. So there are no singular stable SL cones in~$\C^2$.

\begin{princ}{\bf(a)} In the programme of\/ {\rm\S\ref{bs32},} in dimension
$m\ge 3,$ for the Lagrangians $L^t$ at nonsingular times $t\ne T_i$
we should allow Lagrangians with `stable special Lagrangian
singularities'. These should include stable isolated conical
singularities, as in\/ {\rm\cite{Joyc2},} and probably also other
classes of non-isolated or non-conical singularities.

For example, if\/ $m=k+l$ with\/ $k,l>0$ and\/ $C$ is a stable
special Lagrangian cone in $\C^k$ as above, the author expects that
Lagrangians $L$ with\/ $l$-dimensional singularities locally
modelled on $C\t\R^l$ in $\C^k\t\C^l=\C^m$ are `stable'.

In dimension $m=3,$ Lagrangians with conical singularities modelled
on the $T^2$-cone $C$ in \eq{bs2eq4} may be the only kind required.
As $m$ increases, the singularities allowed will probably become
more and more complicated.
\smallskip

\noindent{\bf(b)} For each such class of stable singularities one
should prove short time existence for Lagrangian MCF.
\smallskip

\noindent{\bf(c)} One should extend the definitions of Lagrangian
Floer cohomology, obstructions to $HF^*,$ and $D^b\sF(M)$ to include
each such class of stable singularities.

\label{bs3princ3}
\end{princ}

For (b), the author's PhD student Tapio Behrndt
proved~\cite[Th.~5.12]{Behr2}:

\begin{thm} Let\/ $(M,J,g,\Om)$ be a Calabi--Yau $m$-fold, and\/
$L$ a compact Lagrangian $m$-fold in $M$ with isolated conical
singularities modelled on stable SL cones in $\C^m$ (with any phase
$e^{i\phi}$). Then for small\/ $\ep>0$ there exists a unique smooth
family $\{L^t:t\in[0,\ep)\}$ satisfying Lagrangian MCF with\/
$L^0=L,$ where the $L^t$ are compact Lagrangians in $M$ with stable
isolated conical singularities.
\label{bs3thm1}
\end{thm}

\begin{prob} Extend the theories of Lagrangian Floer cohomology,
obstructions to $HF^*,$ and Fukaya categories $D^b\sF(M)$ to include
Lagrangians $L$ in $M$ with isolated conical singularities modelled
on stable special Lagrangian cones $C$ in\/ $\C^m,$ such as the
$T^2$-cone $C$ in $\C^3$ in \eq{bs2eq4}. The main technical issues
will involve studying moduli spaces of\/ $J$-holomorphic discs $\Si$
in $M$ whose boundaries $\pd\Si$ lie in $L$ and pass through
singular points of\/~$L$.
\label{bs3prob3}
\end{prob}

Problem \ref{bs3prob3} can be approached as an exercise in {\it Symplectic Field Theory}, as in Eliashberg et al.\ \cite{EGH}: given $L$ with conical singularities at $p_1,\ldots,p_k$ modelled on stable SL cones $C_1,\ldots,C_k\subset\C^m$, we delete $p_1,\ldots,p_k$ from $L,M$, and treat $M\sm\{p_1,\ldots,p_k\}$ as a noncompact symplectic manifold with concave cylindrical ends modelled on $\cS^{2m-1}\t(-\iy,0)$, and $L\sm\{p_1,\ldots,p_k\}$ as a noncompact Lagrangian with cylindrical ends modelled on $\Si_j\t(-\iy,0)$ for $j=1,\ldots,k$, where $\Si_j=C_j\cap\cS^{2m-1}$ is the special Legendrian link of the cone~$C_j$.

The reason we need to include Lagrangians with `stable
singularities' in the programme of \S\ref{bs32} is that (the author
expects) for $m\ge 3$ there should exist examples of flows
$\{L^t:t\in[0,T)\}$ in nonsingular Lagrangians with a finite time
singularity at $t=T$, such that one can only continue the flow for
$t>T$ by using Lagrangians with stable singularities.

Example \ref{bs2ex3} described a continuous family of exact SL
3-folds $N^t$ in $\C^3$ for $t\in(-\ep,\ep)$, such that $N^t$ is
nonsingular for $t<0$, and $N^0$ has one (non-stable) singular point
with tangent cone $\R^3\amalg_\R\R^3$, and $N^t$ for $t>0$ has two
singular points modelled on the stable SL $T^2$-cone of \eq{bs2eq4}.
By Principles \ref{bs3princ1}(a) and \ref{bs3princ2}, we should
expect there to exist similar examples of Lagrangian MCF
$L^t:t\in(-\ep,\ep)$ with surgeries, such that $L^t$ is nonsingular
for $t<0$ with a finite time singularity at $t=0$, and $L^0$ has one
singular point with tangent cone $\R^3\amalg_\R\R^3$, and $L^t$ has
two stable singularities modelled on $C$ in~\eq{bs2eq4}.

\begin{rem} We temporarily write $D^b\sF(M)_{\rm nonsing}$ for the derived Fukaya category of nonsingular immersed Lagrangians, and $D^b\sF(M)_{\rm sing}$ for the category including Lagrangians with `stable special Lagrangian singularities'. It seems likely that $D^b\sF(M)_{\rm sing}$ and $D^b\sF(M)_{\rm nonsing}$ {\it need not be equivalent categories}. If so, $D^b\sF(M)_{\rm sing}$ may be {\it preferable\/} to $D^b\sF(M)_{\rm nonsing}$, in the sense of being better behaved, more natural, or the right category to use in Mirror Symmetry. To test this, we should start in dimension $m=3$ by including Lagrangians with isolated singularities modelled on the $T^2$-cone $C$ in~\eq{bs2eq4}.

The following example was suggested to me by Ivan Smith. Harris \cite{Harr} constructs a smooth family $(M^t,\om^t):t\in[0,\ep)$
of symplectic Calabi--Yau 6-manifolds for small $\ep>0$, with the following properties:
\begin{itemize}
\setlength{\parsep}{0pt}
\setlength{\itemsep}{0pt}
\item[(i)] $M^t$ is independent of $t$, and is the result of adding a 2-handle to $T^*\cS^3$. There is an isomorphism $H^2(M^t,\R)\cong\R$ identifying $[\om^t]$ with $t$. Thus $(M^t,\om^t)$ is an exact symplectic manifold if and only if $t=0$.
\item[(ii)] For $t>0$ there is a compact, embedded Lagrangian $L^t$ in $(M^t,\om^t)$ diffeomorphic to $\cS^3$, depending smoothly on $t$, with $0\ne[L^t]\in H_3(M^t;\Z)\cong\Z$.
\item[(iii)] There are no Lagrangian $\cS^3$'s in $(M^0,\om^0)$, and in fact, no compact, exact, embedded Lagrangians in $(M^0,\om^0)$ at all.
\item[(iv)] As in \cite[Rem.~3.7]{Harr}, $L^0=\lim_{t\ra 0}L^t$ is a singular Lagrangian in $M^0$, which topologically looks like an $\cS^3$ with an $\cS^1$ collapsed to a point $p$, so that topologically $L^0$ is modelled on a $T^2$-cone near~$p$.
\end{itemize}

All this suggests that $D^b\sF(M^t)_{\rm nonsing}$ is empty for $t=0$, and nonempty for $t>0$. This counts as pathological behaviour, discontinuous in $t$, since the $D^b\sF(M^t)_{\rm nonsing}$ for small $t>0$ are not deformations of $D^b\sF(M^0)_{\rm nonsing}$ in a meaningful sense. Intuitively, one would expect objects to disappear under small deformations owing to obstructions, so that $D^b\sF(M^t)_{\rm nonsing}$ for $t>0$ should be smaller than~$D^b\sF(M^0)_{\rm nonsing}$.

It seems plausible that we can choose the $L^t$ up to Hamiltonian isotopy so that $L^0$ has one singular point $p$ locally modelled on $C$ in \eq{bs2eq4}, and $L^t$ for $t>0$ is locally modelled near $p$ on $L_1^{A(t)}$ in \eq{bs2eq5}, where $A(t)\ra 0$ as $t\ra 0$. If so, $L^0$ may give an object in $D^b\sF(M^0)_{\rm sing}$, and the derived categories $D^b\sF(M^t)_{\rm sing}$ may depend continuously on $t\in[0,\ep)$. So in this example, $D^b\sF(M)_{\rm sing}$ may be better behaved than $D^b\sF(M)_{\rm nonsing}$ under deformations of~$(M,\om)$.
\label{bs3rem7}
\end{rem}

\subsection{Collapsing zero objects in $D^b\sF(M)$}
\label{bs37}

Let $(L,E)$ be a Lagrangian brane in $\C^m$, either embedded or immersed. Since $L$ is displaceable (Hamiltonian isotopic to a disjoint Lagrangian, by translations in $\C^m$), there are two
possibilities, either:
\begin{itemize}
\setlength{\parsep}{0pt}
\setlength{\itemsep}{0pt}
\item[(A)] $(L,E)$ has $HF^*$ obstructed; or
\item[(B)] $(L,E)$ has $HF^*$ unobstructed, and for every bounding
cochain $b$ for $(L,E)$, $(L,E,b)\cong 0$ in $D^b\sF(\C^m)$. Then we call $(L,E,b)$ a {\it zero object}.

In this case $L$ must also be exact, and strictly immersed (not embedded).
\end{itemize}
For the second part of (B), note that dilation in $\C^m$ induces an infinitesimal deformation of $(L,E,b)$, corresponding to a class in $HF^1\bigl((L,E,b),(L,E,b)\bigr)$. As $(L,E,b)\cong 0$, this deformation class is zero, so dilations of $L$ are Hamiltonian isotopies, and $L$ is exact. But by an argument of Gromov there are no nonempty, compact, exact, embedded Lagrangians in $\C^m$, since then we would have~$H^*(L;\La_\nov)\cong HF^*\bigl((L,\F\t E,0),(L,\F\t E,0)\bigr)\cong 0$. 

\begin{ex} Until recently it was believed there are no compact, graded, embedded Lagrangians in $\C^m$. However, Ekholm, Eliashberg, Murphy and Smith \cite[Cor.~1.6]{EEMS} found an example of a compact, graded, embedded Lagrangian $\cS^1\t\cS^2$ in $\C^3$, and products give Lagrangian $(\cS^1\t\cS^2)^n$'s in $\C^{3n}$. These all have $HF^*$ obstructed, as they are not strictly immersed.
\label{bs3ex2}
\end{ex}

\begin{ex} Writing $\cS^m=\bigl\{(x_0,\ldots,x_m)\in\R^{m+1}:x_0^2+\cdots+x_m^2\bigr\}$, the {\it Whitney sphere\/} $L=\io(\cS^m)$ is the Lagrangian immersion $\io:\cS^m\ra\C^m$ given by
\begin{equation*}
\io:(x_0,x_1,\ldots,x_n)\longmapsto\frac{1}{1+x_0^2}\,\bigl(x_1(1+ix_0),\ldots,x_n(1+ix_0)\bigr).
\end{equation*}
It has the special property of having {\it conformal Maslov form}. It has one transverse self-intersection point at $p=(0,\ldots,0)=\io(1,0,\ldots,0)=\io(-1,0,\ldots,0)$, with $\mu_{L_-,L_+}(p)=-1$, $\mu_{L_+,L_-}(p)=m+1$. Thus if $m>2$, Lemma \ref{bs2lem2} shows that $L$ has $HF^*$ unobstructed, so as in (B), $(L,E,b)\cong 0$ in~$D^b\sF(\C^m)$.

Ekholm, Eliashberg, Murphy and Smith \cite[\S 1]{EEMS} construct Lagrangian immersions $\jmath:\cS^m\ra\C^m$ for $m$ odd, with one transverse self-intersection point $p$ with $\mu_{L_+,L_-}(p)=2$. If $m\ge 3$ it has $HF^*$ obstructed, as in~(A).
\label{bs3ex3}
\end{ex}

Next we consider graded, immersed Lagrangian MCF in an example in~$\C$.

\begin{ex} Let $L$ be a graded, immersed Lagrangian in $\C$ shaped
like an $\iy$ sign, not necessarily symmetric, bounding two
`teardrop' $J$-holomorphic curves $\Si_1,\Si_2$, as shown in Figure
\ref{bs3fig4}, and let $E\ra L$ be a rank one $\F$-local system, which is classified by its holonomy $\Hol(\nabla_E)[L]\in\F^*$ around $L$.

Then $(L,E)$ has $HF^*$ obstructed if $\area(\Si_1)\ne\area(\Si_2)$. If $\area(\Si_1)=\area(\Si_2)$, there is a unique choice of $\Hol(\nabla_E)[L]=\pm 1$ which makes the obstructions to $HF^*$ due to $\Si_1,\Si_2$ cancel, and then $(L,E)$ has $HF^*$ unobstructed.
 
\begin{figure}[htb]
\centerline{$\splinetolerance{.8pt}
\begin{xy}
0;<1mm,0mm>:
,(10,0);(-20,0)**\crv{(10,5)&(0,10)}
?(.95)="aa"
?(.88)="bb"
?(.8)="cc"
?(.72)="dd"
?(.63)="ee"
?(.53)="ff"
?(.42)="gg"
?(.27)="ii"
?(.6)="yy"
,(10,0);(-20,0)**\crv{(10,-5)&(0,-10)}
?(.95)="kk"
?(.88)="ll"
?(.8)="mm"
?(.72)="nn"
?(.63)="oo"
?(.53)="pp"
?(.42)="qq"
?(.27)="ss"
?(.6)="zz"
,"aa";"kk"**@{.}
,"bb";"ll"**@{.}
,"cc";"mm"**@{.}
,"dd";"nn"**@{.}
,"ee";"oo"**@{.}
,"ff";"pp"**@{.}
,"gg";"qq"**@{.}
,"ii";"ss"**@{.}
,(-59,0);(-20,0)**\crv{(-59,6)&(-40,10)}
?(.95)="a"
?(.85)="b"
?(.75)="c"
?(.65)="d"
?(.55)="e"
?(.45)="f"
?(.35)="g"
?(.25)="h"
?(.15)="i"
?(.05)="j"
?(.5)="y"
,(-59,0);(-20,0)**\crv{(-59,-6)&(-40,-10)}
?(.95)="k"
?(.85)="l"
?(.75)="m"
?(.65)="n"
?(.55)="o"
?(.45)="p"
?(.35)="q"
?(.25)="r"
?(.15)="s"
?(.05)="t"
?(.5)="z"
,"a";"k"**@{.}
,"b";"l"**@{.}
,"c";"m"**@{.}
,"d";"n"**@{.}
,"e";"o"**@{.}
,"f";"p"**@{.}
,"g";"q"**@{.}
,"h";"r"**@{.}
,"i";"s"**@{.}
,"j";"t"**@{.}
,(-20,0)*{\bu}
,(-40,0)*{\Si_1}
,(0,0)*{\Si_2}
,(13,0)*{L}
\end{xy}$}
\caption{`$\iy$ sign' Lagrangian $L$ in $\C$}
\label{bs3fig4}
\end{figure}

Consider the immersed Lagrangian MCF (`curve shortening flow')
$\{L^t:t\in[0,T)\}$ in $\C$ starting from $L^0=L$ with first finite
time singularity at $t=T$. The curve shortening flow is well
understood, as in Abresch and Langer \cite{AbLa}, Angenent
\cite{Ange1,Ange2}, Grayson \cite{Gray}, and others, and we can give
a good description of the flow. The difference
$\area(\Si_1^t)-\area(\Si_2^t)$ is constant during the flow, and
both $\area(\Si_1^t),\area(\Si_2^t)$ decrease until the smaller
becomes zero at $t=T$.

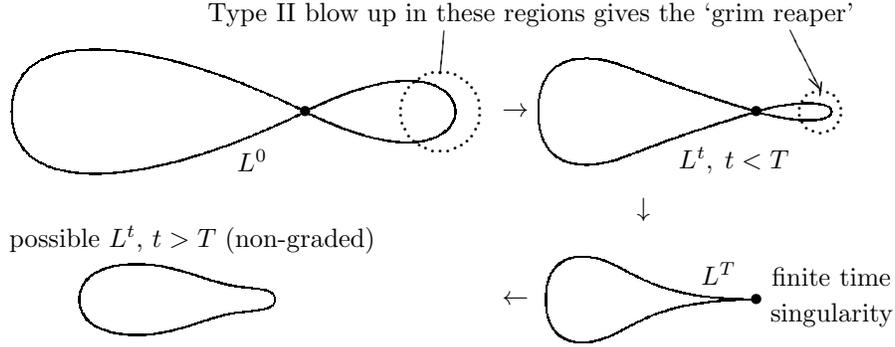
\begin{figure}[htb]
\centerline{$\splinetolerance{.8pt}
\begin{xy}
0;<1mm,0mm>:
,(0,0);(-20,0)**\crv{(0,3)&(-5,7.5)}
,(0,0);(-20,0)**\crv{(0,-3)&(-5,-7.5)}
,(-59,0);(-20,0)**\crv{(-59,12)&(-40,10)}
,(-59,0);(-20,0)**\crv{(-59,-12)&(-40,-10)}
,(-20,0)*{\bu}
,(-27,-6.5)*{L^0}
,(8,0)*{\ra}
,(50,0);(40,0)**\crv{(50,1.5)&(45,1.5)}
,(50,0);(40,0)**\crv{(50,-1.5)&(45,-1.5)}
,(11,0);(40,0)**\crv{(11,7)&(21,8)&(25,4.5)}
,(11,0);(40,0)**\crv{(11,-7)&(21,-8)&(25,-4.5)}
,(40,0)*{\bu}
,(37,-6.5)*{L^t,\; t<T}
,(25,-13)*{\downarrow}
,(-2,0)*\xycircle<15pt>{{.}}
,(48.5,0)*\xycircle<8pt>{{.}}
,(-50,-25);(-24,-25)**\crv{(-50,-21)&(-40,-19)&(-30,-24.1)&(-24,-23.3)}
,(-50,-25);(-24,-25)**\crv{(-50,-29)&(-40,-31)&(-30,-25.9)&(-24,-26.7)}
,(-35,-17)*{\text{possible $L^t$, $t>T$ (non-graded)}}
,(10,13)*{\text{Type II blow up in these regions gives the `grim reaper'}}
,(-2,5.2);(-1.5,11)**\crv{}
,(-2,5.2);(-2.7,6.4)**\crv{}
,(-2,5.2);(-1.2,6.4)**\crv{}
,(48.5,3);(44.5,11)**\crv{}
,(48.5,3);(47,4.4)**\crv{}
,(48.5,3);(48.7,4.8)**\crv{}
,(8,-25)*{\leftarrow}
,(12,-25);(40,-25)**\crv{(12,-19)&(21,-18)&(27,-25)}
,(12,-25);(40,-25)**\crv{(12,-31)&(21,-32)&(27,-25)}
,(40,-25)*{\bu}
,(35,-21.5)*{L^T}
,(50,-22)*{\text{finite time}}
,(50,-27)*{\text{singularity}}
\end{xy}$}
\caption{Lagrangian MCF when $\area(\Si_1^t)>\area(\Si_2^t)$}
\label{bs3fig5}
\end{figure}

In the case $\area(\Si_1^t)>\area(\Si_2^t)$, the flow is sketched in
Figure \ref{bs3fig5}. The loop bounding $\Si_2$ shrinks to a point
at $t=T$, and the curve develops a cusp singularity. A type II blow
up of this singularity sees only the small, highly curved regions
indicated, and yields the `grim reaper' translating soliton from
Figure \ref{bs2fig1}. Note that in this case, the type II blow up
only gives a rather incomplete picture of what is happening.

Following Angenent \cite{Ange2}, one can  continue the flow for
$t>T$ after a surgery at $t=T$ eliminating the self-intersection
point, as in the last picture of Figure \ref{bs3fig5}, but then the
$L^t$ for $t>T$ are non-graded. From the point of view of this
paper, this is the wrong thing to do, and only works as dimension
$m=1$ is so simple. A better answer is that {\it after the
singularity at\/ $t=T,$ one cannot continue the flow in graded
Lagrangian MCF for\/} $t>T$. This does not contradict the programme
of \S\ref{bs32}, as the initial Lagrangian $L$ in Figure
\ref{bs3fig4} has $HF^*$ obstructed in this case. We will discuss
this phenomenon further in~\S\ref{bs38}.

\begin{figure}[htb]
\centerline{$\splinetolerance{.8pt}
\begin{xy}
0;<.9mm,0mm>:
,(0,0);(-30,0)**\crv{(0,12)&(-20,8)}
,(0,0);(-30,0)**\crv{(0,-12)&(-20,-7)}
,(-60,0);(-30,0)**\crv{(-60,12)&(-40,8)}
,(-60,0);(-30,0)**\crv{(-60,-12)&(-40,-8)}
,(-30,0)*{\bu}
,(-30,8.5)*{L^0}
,(5.3,0)*{\ra}
,(61,0);(36,0)**\crv{(61,7)&(51,5)}
,(61,0);(36,0)**\crv{(61,-7)&(51,-5)}
,(11,0);(36,0)**\crv{(11,7)&(21,5)}
,(11,0);(36,0)**\crv{(11,-7)&(21,-5)}
,(36,0)*{\bu}
,(36,8.5)*{L^{t_1},\; 0\!<\!t_1\!<\!t_2\!<\!T}
,(60,-14)*{\downarrow}
,(-7,0)*\xycircle<25pt>{{.}}
,(-53,0)*\xycircle<25pt>{{.}}
,(57,0)*\xycircle<18pt>{{.}}
,(15,0)*\xycircle<18pt>{{.}}
,(-30,-25)*{\bu}
,(5.3,-25)*{\leftarrow}
,(51,-25);(36,-25)**\crv{(51,-23)&(46,-23.5)}
,(51,-25);(36,-25)**\crv{(51,-27)&(46,-26.5)}
,(21,-25);(36,-25)**\crv{(21,-23)&(26,-23.5)}
,(21,-25);(36,-25)**\crv{(21,-27)&(26,-26.5)}
,(36,-25)*{\bu}
,(36,-31.5)*{L^{t_2},\; 0\!<\!t_1\!<\!t_2\!<\!T}
,(-30,-29)*{L^T}
,(-45,-22)*{\text{finite time}}
,(-45,-27)*{\text{singularity}}
,(50,-25)*\xycircle<7pt>{{.}}
,(22,-25)*\xycircle<7pt>{{.}}
,(-12.2,-8.4);(-17,-12)**\crv{}
,(-12.2,-8.4);(-14.2,-8.5)**\crv{}
,(-12.2,-8.4);(-12.7,-10.2)**\crv{}
,(-47.8,-8.4);(-43,-12)**\crv{}
,(-47.8,-8.4);(-45.8,-8.5)**\crv{}
,(-47.8,-8.4);(-47.3,-10.2)**\crv{}
,(53.8,-6.4);(49,-12)**\crv{}
,(53.8,-6.4);(51.8,-7.4)**\crv{}
,(53.8,-6.4);(53.7,-8.2)**\crv{}
,(18.2,-6.4);(23,-12)**\crv{}
,(18.2,-6.4);(20.2,-7.4)**\crv{}
,(18.2,-6.4);(18.7,-8.2)**\crv{}
,(50,-22);(49,-16)**\crv{}
,(50,-22);(49,-20.6)**\crv{}
,(50,-22);(50.5,-20.4)**\crv{}
,(22,-22);(23,-16)**\crv{}
,(22,-22);(23,-20.6)**\crv{}
,(22,-22);(21.5,-20.4)**\crv{}
,(0,-14)*{\text{Type II blow up in these regions gives the `grim reaper'}}
\end{xy}$}
\caption{Lagrangian MCF when $\area(\Si_1^t)=\area(\Si_2^t)$}
\label{bs3fig6}
\end{figure}
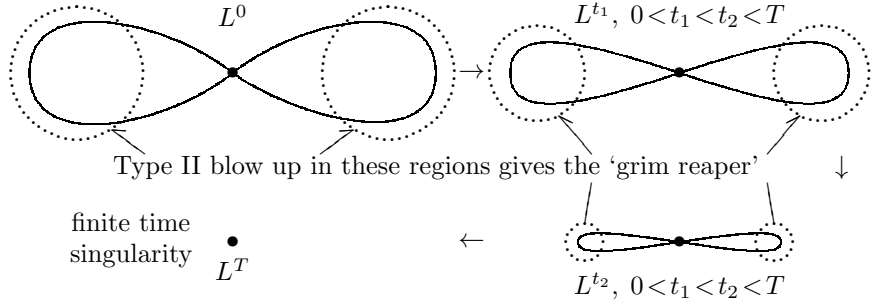

In the case $\area(\Si_1^t)=\area(\Si_2^t)$, the flow is sketched in
Figure \ref{bs3fig6}. The whole $\iy$ sign shrinks to a point at
$t=T$. It is not a type I singularity modelled on a Lagrangian MCF
shrinker, since this cannot happen in graded Lagrangian MCF as in
\S\ref{bs23}. The curve does not rescale homothetically, but as in
Figure \ref{bs3fig6} the curve shrinks faster in the vertical than
in the horizontal directions. Type II blow ups at either end of the
$\iy$ sign yield a `grim reaper' translating soliton, as in Figure
\ref{bs2fig1}, as indicated. So, {\it in this case of an immersed
curve in $\C$ with\/ $HF^*$ unobstructed, the whole curve collapses
to a point in finite time under Lagrangian MCF}.
\label{bs3ex4}
\end{ex}

More generally, for Lagrangian MCF $\{L^t:t\in[0,T)\}$ of compact,
immersed, graded Lagrangians $L^t$ in $\C^m$ with $HF^*$
unobstructed, I expect that the typical behaviour is for the whole
of $L^t$ to collapse to a point at time $t=T$ (though possibly
undergoing other surgeries along the way, as in~\S\ref{bs34}--\S\ref{bs36}).

Similarly, for immersed Lagrangian MCF $\{L^t:t\in[0,T)\}$ with
$HF^*$ unobstructed in a Calabi--Yau $m$-fold $(M,J,g,\Om)$,
connected components $L^t_1$ of $L^t=L^t_1\amalg L^t_2$ in small
open balls in $M$ may collapse to a point in finite time $t=T$. When
this happens, in the programme of \S\ref{bs32}, the correct thing to
do is to delete the collapsed component $L^t_1$, and continue
flowing the remaining components $L^t_2$ when $t>T$. As
$(L^t_1,E^t_1,b^t_1)$ is a zero object in $D^b\sF(M)$, deleting it does
not change the isomorphism class in $D^b\sF(M)$. We state this as:

\begin{princ} The following behaviour, called `collapsing a zero
object', is a possible model for finite time singularities in the
programme of\/~{\rm\S\ref{bs32}}.

Let\/ $(M,J,g,\Om)$ be a Calabi--Yau $m$-fold, and extend\/
$D^b\sF(M)$ to include immersed Lagrangians, as in\/
{\rm\cite{AkJo}}. Suppose $\{(L^t,E^t):t\in(T-\ep,T+\ep)\}$ for $\ep>0$
small is a family of Lagrangian branes in $M$ with\/ $HF^*$ unobstructed, and\/ $\{b^t:t\in(T-\ep,T+\ep)\}$ a corresponding family of bounding cochains, satisfying the following conditions:
\begin{itemize}
\setlength{\parsep}{0pt}
\setlength{\itemsep}{0pt}
\item[{\bf(i)}] The $(L^t,E^t,b^t)$ for $t\in(T-\ep,T+\ep)$ are all
isomorphic in $D^b\sF(M)$.
\item[{\bf(ii)}] When $t<T,$ $L^t,E^t$ depend smoothly on
$t\in(T-\ep,T),$ and\/ $\{(L^t,E^t):t\in(T-\ep,T)\}$ satisfies
Lagrangian MCF, with a finite time singularity at\/ $t=T,$ with
one singular point\/~$p\in M$.

Similarly, when $t\ge T,$ $L^t,E^t$ depend smoothly on
$t\in[T,T+\ep),$ and\/ $\{(L^t,E^t):t\in[T,T+\ep)\}$ satisfies
Lagrangian MCF.
\item[{\bf(iii)}] For $t\in(T-\ep,T)$ there is a decomposition
$(L^t,E^t,b^t)=(L^t_1,E^t_1,b^t_1)\amalg (L^t_2,E^t_2,b^t_2),$ with\/ $L^t_1,L^t_2$ open and closed in $L^t$. There exists a continuous $\de:(T-\ep,T)\ra(0,\iy)$ with\/ $\de(t)\ra 0$ as $t\ra T$ such
that\/ $L^t_1\subseteq B_{\de(t)}(p)$ for all\/ $t\in(T-\ep,T),$
where $B_{\de(t)}(p)$ is the open ball of radius $\de(t)$
about\/ $p$ in $M$. That is, the whole of\/ $L^t_1$ converges
uniformly to $p\in M$ as\/~$t\ra T$.
\item[{\bf(iv)}] $(L^t_1,E^t_1,b^t_1)\cong 0$ in $D^b\sF(M)$ for $t\in(T-\ep,T),$ so that\/ $(L^t,E^t,b^t)\cong(L^t_2,E^t_2,b^t_2)$ in $D^b\sF(M)$.
\item[{\bf(v)}] The family $\{(L^t_2,E^t_2,b^t_2):t\in(T-\ep,T)\}\amalg\{(L^t,E^t,b^t):t\in[T,T+\ep)\}$ is smooth in $t\in(T-\ep,T+\ep)$.
\end{itemize}
Rather than taking $L^T=L^T_2$ to be a nonsingular immersed
Lagrangian at\/ $t=T,$ we could instead write $L^T=\{p\}\amalg
L^T_2,$ where $\{p\}=\lim_{t\ra T}L^t_1$ is regarded as an extreme
example of a singular Lagrangian in~$M$.
\label{bs3princ4}
\end{princ}

Recall that a graded Lagrangian $L$ is {\it almost calibrated\/} if
it has phase variation less than $\pi$. The almost calibrated
condition is preserved by Lagrangian MCF. The next lemma implies
that `collapsing zero objects' does not happen in almost calibrated
Lagrangian MCF.

\begin{lem} Suppose $L$ is a compact, immersed, graded Lagrangian in
$\C^m,$ or in a small open ball\/ $B_\de(p)$ in a Calabi--Yau
$m$-fold\/ $(M,J,g,\Om)$. Then $L$ has phase variation greater than
$\pi$. That is, $L$ is not almost calibrated.
\label{bs3lem}
\end{lem}

To prove the lemma, assume for a contradiction that the phase
function $\th_L$ of $L$ maps $\th_L:L\ra[\phi-\frac{\pi}{2},
\phi+\frac{\pi}{2}]$, consider $\int_L(\cos\phi\Re\Om-
\sin\phi\Im\Om)$, and note that the homology class $[L]$ in
$H_m(\C^m,\Z)$ or $H_m(M,\Z)$ is zero.

\begin{prob} Find global geometric models for how such `collapsing
a zero object' finite time singularities occur in MCF for compact,
immersed, graded Lagrangians $L^t$ in $\C^m$ with\/ $HF^*$ unobstructed.

Even for\/ $m=1$ there may be something new to say.
\label{bs3prob4}
\end{prob}

\begin{ex} Let $(M,J,g,\Om)$ be a Calabi--Yau $m$-fold for $m\ge
2$, $L$ a compact Lagrangian in $M$, and $p\in L$. In \cite{Neve3},
Neves defines another Lagrangian $\ti L$ in $M$, which is
Hamiltonian isotopic to $L$ and coincides with $L$ except in a small
open neighbourhood of $p$. Here $L,\ti L$ are locally $\SO(m)$
surfaces of revolution on the curves in sketched in Figure
\ref{bs3fig7}. (Actually Neves restricts to $m=2$, but the same
ideas should work for all~$m\ge 2$.)

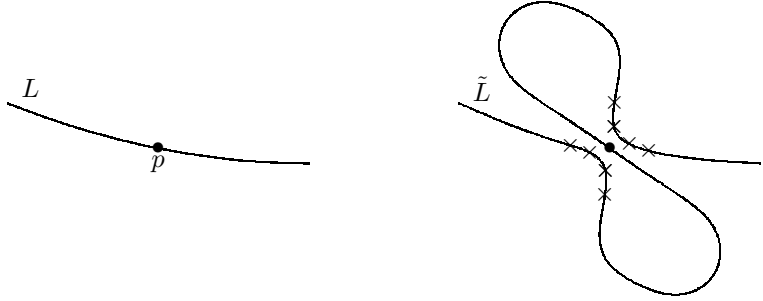
\begin{figure}[htb]
\centerline{$\splinetolerance{.8pt}
\begin{xy}
0;<1mm,0mm>:
,(-80,6);(-60,0)**\crv{(-70,2)}
,(-40,-2);(-60,0)**\crv{(-50,-2)}
,(-60,0)*{\bu}
,(-60,-2)*{p}
,(0,0)*{\bu}
,(-20,6);(0,0)**\crv{(-10,1.5)&(-2,0)&(0,-3)&(-1,-8)&(-2,-15)&(5,-19)
&(10,-20)&(15,-17)&(15,-9)&(7,-5)}
?(.09)="a"
?(.13)="b"
?(.17)="c"
?(.21)="d"
?(.29)="e"
,"a"*{\t}
,"b"*{\t}
,"d"*{\t}
,"d"*{\t}
,"e"*{\t}
,(20,-2);(0,0)**\crv{(10,-1.5)&(2,0)&(0,3)&(1,8)&(2,15)&(-5,19)
&(-10,20)&(-15,17)&(-15,9)&(-7,5)}
?(.09)="a"
?(.13)="b"
?(.17)="c"
?(.21)="d"
?(.29)="e"
,"a"*{\t}
,"b"*{\t}
,"d"*{\t}
,"d"*{\t}
,"e"*{\t}
,(-77,8)*{L}
,(-17,8)*{\ti L}
\end{xy}$}
\caption{Neves' Lagrangian with a finite time singularity under LMCF}
\label{bs3fig7}
\end{figure}

Neves' main result \cite[Th.~A]{Neve3} is that Lagrangian MCF
starting from $\ti L$ develops finite time singularities. This is
important, as it shows that finite time singularities in Lagrangian
MCF are unavoidable in many situations (although note that $\ti L$
has phase variation greater than $\pi$, so this does not show that
almost calibrated Lagrangian MCF has finite time singularities).

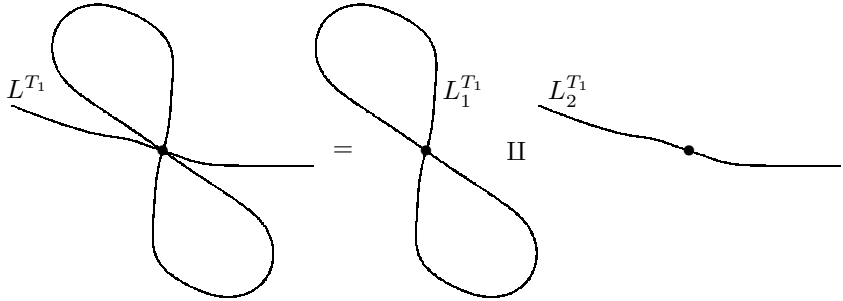
\begin{figure}[htb]
\centerline{$\splinetolerance{.8pt}
\begin{xy}
0;<1mm,0mm>:
,(-20,6);(0,0)**\crv{(-10,2)&(-5,2)&(-3,1)}
,(0,0);(0,0)**\crv{(-1,-3)&(-1,-8)&(-2,-15)&(5,-19)
&(10,-20)&(15,-17)&(15,-9)&(7,-5)}
,(20,-2);(0,0)**\crv{(10,-2)&(5,-2)&(3,-1)}
,(0,0);(0,0)**\crv{(1,3)&(1,8)&(2,15)&(-5,19)
&(-10,20)&(-15,17)&(-15,9)&(-7,5)}
,(-18,8.5)*{L^{T_1}}
,(0,0)*{\bu}
,(35,0);(35,0)**\crv{(34,-3)&(34,-8)&(33,-15)&(40,-19)
&(45,-20)&(50,-17)&(50,-9)&(42,-5)}
,(35,0);(35,0)**\crv{(36,3)&(36,8)&(37,15)&(30,19)
&(25,20)&(20,17)&(20,9)&(28,5)}
,(35,0)*{\bu}
,(40,8)*{L^{T_1}_1}
,(50,6);(70,0)**\crv{(60,2)&(65,2)&(67,1)}
,(90,-2);(70,0)**\crv{(80,-2)&(75,-2)&(73,-1)}
,(70,0)*{\bu}
,(54,8)*{L^{T_1}_2}
,(24,0)*{=}
,(47,0)*{\amalg}
\end{xy}$}
\caption{First singular time $t=T_1$ of Lagrangian MCF from $\ti L$}
\label{bs3fig8}
\end{figure}

What actually happens in Lagrangian MCF starting from $\ti L$?
Neves' proof does not tell us, as he assumes for a contradiction
that no finite time singularity occurs. The author expects a Lagrangian MCF
with surgeries $\{L^t:t\in[0,T)\}$ in $M$ with $L^0=\ti L$, with two
singular times $0<T_1<T_2<T$. For $t\in[0,T_1),$ $L^t$ looks much
like $\ti L$, but as $t\ra T_1$ in $[0,T_1)$, the region marked with
crosses `$\t$' in Figure \ref{bs3fig7} undergoes a `neck pinch'. At
$t=T_1$, as sketched in Figure \ref{bs3fig8}, $L^{T_1}$ decomposes
as $L^{T_1}_1\amalg L^{T_1}_2$, where $L^{T_1}_1$ is a small
immersed $\cS^m$ near $p$ with one transverse self-intersection
point with $\mu_{L_+,L_-}(p)=m+1$, a `Whitney sphere' as in Example \ref{bs3ex3}, and $L^{T_1}_2$ looks quite like the original~$L$.

Then as $t$ increases from $T_1$ to $T_2$, the component $L^t_1$
should shrink to a point, until at the second singular time $t=T_2$
it undergoes `collapsing a zero object' as in Principle
\ref{bs3princ4}. Meanwhile, the Lagrangian MCF of $L^t_2$ looks
quite like that of the original $L$, and continues for $t>T_2$.
Thus, at least conjecturally, Neves' examples \cite{Neve3} are not counterexamples to the programme of~\S\ref{bs32}.
\label{bs3ex5}
\end{ex}

\subsection{What goes wrong in LMCF of obstructed Lagrangians}
\label{bs38}

The programme of \S\ref{bs32} claims that if $(M,J,g,\Om)$ is a Calabi--Yau $m$-fold and $L$ a compact, immersed, graded Lagrangian in $M$ with $HF^*$ unobstructed, then graded Lagrangian MCF with surgeries $\{L^t:t\in[0,\iy)\}$ with $L^0=L$ should exist for all time. But if $L$ has $HF^*$ obstructed, the author expects that Lagrangian MCF $\{L^t:t\in[0,T)\}$ can develop finite time singularities at $t=T$ such that {\it one cannot continue the flow for\/} $t>T$, even after a surgery.

In dimension $m=1$, we met an example of this in Example \ref{bs3ex4}: if $L$ is the `$\iy$ sign' Lagrangian in $\C$ from Figure \ref{bs3fig4} with $\area(\Si_1)\ne\area(\Si_2)$, then Lagrangian MCF starting from $L$ has a finite time singularity after which one cannot continue in graded Lagrangian MCF (though in this case one can continue in non-graded Lagrangian MCF after a surgery).

We now discuss the nature of these terminal singularities of $HF^*$-obstructed Lagrangian MCF. I expect they should be impossible in $HF^*$-unobstructed flow, and so the obstructions should be present locally as the singularity forms. As in \S\ref{bs25}--\S\ref{bs26}, obstructions to $HF^*$ for a Lagrangian $L$ or brane $(L,E)$ are caused by `bad' $J$-holomorphic discs $\Si$ in $M$ with boundary in $L$, of two kinds:
\begin{itemize}
\setlength{\parsep}{0pt}
\setlength{\itemsep}{0pt}
\item[(i)] For $L$ embedded, moduli spaces $\oM_1^A$ of $J$-holomorphic discs $\Si$ with area $A>0$ and one boundary marked point, whose virtual classes $\bigl[\bigl[\oM_1^A\bigr]{}_{\rm virt}\bigr]$ are nonzero in $H_{m-2}(L,\Q)$. (This is oversimplified.)
\item[(ii)] For $L$ immersed, $\Si$ of type (i), and also `teardrop-shaped'
$J$-holomorphic discs $\Si$ of the form shown in Figure
\ref{bs2fig3}, with one corner at $q\in M$, and with
$\mu_{L_+,L_-}(q)=2$, where $L_\pm$ are the local sheets of $L$
intersecting at $q$.
\end{itemize}
Thus an obvious guess is that the singularities we are interested in occur when such a `bad' $\Si$ shrinks to a point, and $\area(\Si)\ra 0$. As $L$ is graded, $\Si$ of type (i) have constant area under Lagrangian MCF, so they are not relevant. For $\Si$ of type (ii), as for \eq{bs3eq8} under Lagrangian MCF we have
\begin{equation*}
\frac{\d}{\d t}\area(\Si)=-\int_{\pd\Si}\d\th_L=
\th_{L_-}(q)-\th_{L_+}(q),
\end{equation*}
so $\area(\Si)$ will decrease under Lagrangian MCF if $\th_{L_+}(q)>\th_{L_-}(q)$.

Therefore we propose:

\begin{princ} In contrast to\/ {\rm\S\ref{bs32},} Lagrangian MCF\/ $\{L^t:t\in[0,T)\}$ of compact, immersed, graded Lagrangians $L$ or branes $(L,E)$ with\/ $HF^*$ obstructed in a Calabi--Yau $m$-fold may develop finite time singularities at\/ $t=T,$ such that one cannot continue the flow for $t>T$ in graded LMCF, even after a surgery.

A typical way in which this occurs is that for $t\in(T-\ep,T),$ there exists a `teardrop' $J$-holomorphic curve $\Si^t$ with boundary in $L^t$ of the form shown in Figure\/ {\rm\ref{bs2fig3},} and\/ $\area(\Si^t)\ra 0$ as $t\ra T,$ where $\Si^t$ causes $L^t$ to have $HF^*$ obstructed if\/ $\area(\Si^t)$ is small enough.

In dimension $m\ge 2,$ this should be possible for $L^0$ with arbitrarily small phase variation.
\label{bs3princ5}
\end{princ}

Note that this is exactly what happens in Example \ref{bs3ex4} in
dimension~$m=1$.

\begin{rem} We are restricting to graded Lagrangians, so as above,
discs $\Si$ of type (i) have constant area under the flow, and
cannot cause singularities.

We could generalize the programme of \S\ref{bs32} to oriented
Lagrangians rather than graded Lagrangians, so that
$HF^*\bigl((L,E,b),(L',E',b')\bigr)$ is $\Z_2$-graded rather than
$\Z$-graded. In this case, curves of type (i) can cause
singularities. For non-graded $L$, the area of curves $\Si$ of type
(i) change under Lagrangian MCF by
\e
\frac{\d}{\d t}\area(\Si^t)=-\mu_L\cdot[\pd\Si^t],
\label{bs3eq12}
\e
where $\mu_L\in H^1(L,\R)$ is the Maslov class from \S\ref{bs21},
and $[\pd\Si^t]\in H_1(L,\R)$. As the r.h.s.\ of \eq{bs3eq12} is
independent of $t$, if $\mu_L\cdot[\pd\Si^0]>0$ then unless other
singularities happen first, the area of $\Si^t$ shrinks to zero at
time $T=\area(\Si^0)/(\mu_L\cdot[\pd\Si^0])$. So in the non-graded
analogue of Principle \ref{bs3princ5}, we should also include
shrinking of type (i) discs $\Si$. Groh, Schwarz, Smoczyk and
Zehmisch \cite{GSSZ} used this idea to study singularities of
Lagrangian MCF for monotone Lagrangians in~$\C^m$.
\label{bs3rem8}
\end{rem}

\begin{ex} Wolfson \cite{Wolf} constructed an example of a
Calabi--Yau 2-fold $(M,J,g,\Om)$ (a $K3$ surface) with the following
properties:
\begin{itemize}
\setlength{\parsep}{0pt}
\setlength{\itemsep}{0pt}
\item[(i)] There exists $\al\in H_2(M,\Z)$ with
$\al\cdot\al=-4$, such that every compact, immersed Lagrangian
$L$ in $M$ has $[L]\in\Z\cdot\al\subset H_2(M,\Z)$.
\item[(ii)] There exists an immersed Lagrangian two-sphere $L$
in $M$ with $[L]=\al$.
\item[(iii)] There does not exist a compact, immersed SL 2-fold
$L'$ in $M$ with homology class $\al$ (even if one allows branch
point singularities in $L'$).
\end{itemize}
Here (iii) is proved as follows: $L'$ must be connected, as we
cannot split $\al=\be+\ga$ for $\be\ne 0\ne\ga$ homology classes
represented by SL 2-folds. Suppose $L'$ has genus $g$, and for
simplicity has $k$ transverse self-intersection points. An easy
calculation shows that $[L']\cdot[L']=2g+2k-2\ge -2$. But $[L']=\al$
and~$\al\cdot\al=-4$.

So we can ask: what happens to Lagrangian MCF $\{L^t:t\in[0,T)\}$ in
$M$ with $L^0=L$? I expect that $L$ has $HF^*$ obstructed, and that
a finite time singularity develops at $t=T$ after which one cannot
continue the (graded) flow, as in Principle \ref{bs3princ5}. As
evidence for this, note that if Lagrangian MCF with surgeries
$\{L^t:t\in[0,\iy)\}$ existed for all time, one would expect
$L'=\lim_{t\ra\iy}L^t$ to be an SL 2-fold in homology class $\al$,
which is excluded by~(iii).

Wolfson uses his example to prove something different. Schoen and
Wolfson \cite{ScWo} show that by minimizing volume amongst (not
necessarily graded) compact, immersed, oriented Lagrangians $L$ in a
Calabi--Yau 2-fold in a fixed homology class $\al$ and taking a
limit, one can construct a singular Lagrangian $L'$ with minimal
volume in homology class $\al$, such that $L'$ is Hamiltonian
stationary and has finitely many singular points of two kinds:
\begin{itemize}
\setlength{\parsep}{0pt}
\setlength{\itemsep}{0pt}
\item[(a)] Branch points, like those of Riemann surfaces, and
\item[(b)] Singularities modelled on certain Lagrangian cones
$C_{p,p+1}$ in $\C^2$ for $p\ge 1.$ These $C_{p,p+1}$ are
Hamiltonian stationary, but not Maslov zero, or graded.
\end{itemize}
If there are only singular points of type (a), then $L'$ is special
Lagrangian. Wolfson deduces \cite[Th.~3.3]{Wolf} that in his
example, the minimizer $L'$ must have singular points of type (b).
But then $L'$ is not graded, so it is not a possible limit
$\lim_{t\ra\iy}L^t$ for graded Lagrangian MCF.
\label{bs3ex6}
\end{ex}

The next example gives a heuristic description of how the author
expects the finite time singularities in Principle \ref{bs3princ5}
may form geometrically.

\begin{ex} Example \ref{bs2ex6} described a family of Lagrangian
MCF translators $L$ in $\C^m$ given in equation \eq{bs2eq10},
asymptotic to the union of two Lagrangian planes $\Pi_0,\Pi_{\bs\phi}\cong\R^m$ intersecting in $\R$. We have sketched $L$ in Figure \ref{bs3fig9}
(not easy to draw in only two dimensions).
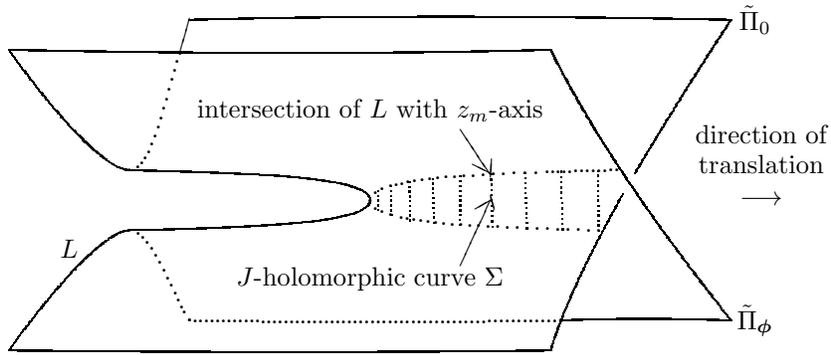
\begin{figure}[htb]
\centerline{$\splinetolerance{.8pt}
\begin{xy}
0;<.8mm,0mm>:
,(-30,30);(60,30)**\crv{(10,31)}
,(-60,25);(30,25)**\crv{(-10,24)}
,(-60,25);(-40,5)**\crv{(-53,15)&(-44,5)}
,(-40,5);(-31.6,24.5)**\crv{~*=<3pt>{.} (-36,5)&(-34,17.5)}
,(-30,30);(-31.6,24.5)**\crv{}
,(-40,5);(0,0)**\crv{(-20,4.5)&(0,4)}
,(-40,-5);(0,0)**\crv{(-20,-4.5)&(0,-4)}
,(-60,-25);(-40,-5)**\crv{(-53,-15)&(-44,-5)}%
,(-40,-5);(-30,-20)**\crv{~*=<3pt>{.} (-36,-5)&(-32,-17.5)}
,(-60,-25);(30,-25)**\crv{(-10,-25.5)}
,(30,-25);(42,1)**\crv{(33,-15)}
,(60,30);(44,4.5)**\crv{(52,18)}
,(30,25);(60,-20)**\crv{(40,5)}
,(-30,-20);(31.7,-20)**\crv{~*=<3pt>{.} (0,-20.5)}
,(31.7,-20);(60,-20)**\crv{(45,-19.5)}
,(38,-5);(0,0)**\crv{~*=<3pt>{.} (20,-4.5)&(0,-4)}
?(.845)="b"
?(.745)="c"
?(.645)="d"
?(.54)="e"
?(.44)="f"
?(.33)="g"
?(.22)="h"
?(.11)="i"
?(.0)="j"
,(41,5);(0,0)**\crv{~*=<3pt>{.} (20,4.5)&(0,4)}
?(.85)="l"
?(.75)="m"
?(.65)="n"
?(.55)="o"
?(.45)="p"
?(.35)="q"
?(.25)="r"
?(.15)="s"
?(.05)="t"
,"b";"l"**@{.}
,"c";"m"**@{.}
,"d";"n"**@{.}
,"e";"o"**@{.}
,"f";"p"**@{.}
,"g";"q"**@{.}
,"h";"r"**@{.}
,"i";"s"**@{.}
,"j";"t"**@{.}
,(65,0)*{\longra}
,(65,11)*{\text{direction of}}
,(65,6)*{\text{translation}}
,(20,4.5);(15,12)**\crv{}
,(20,4.5);(17,6)**\crv{}
,(20,4.5);(20.2,7.5)**\crv{}
,(0,15)*{\text{intersection of $L$ with $z_m$-axis}}
,(20,0);(15,-11)**\crv{}
,(20,0);(17.5,-1.5)**\crv{}
,(20,0);(20.7,-2.5)**\crv{}
,(0,-13)*{\text{$J$-holomorphic curve $\Si$}}
,(64,30)*{\ti\Pi_0}
,(64,-20)*{\ti\Pi_{\bs\phi}}
,(-50,-8)*{L}
\end{xy}$}
\caption{Joyce--Lee--Tsui Lagrangian MCF translator from Example \ref{bs2ex6}}
\label{bs3fig9}
\end{figure}
We indicate the intersection of $L$ with the $z_m$-axis, the curve
\begin{align*}
L\,\cap\, &\bigl\{(0,\ldots,0,z_m):z_m\in\C\bigr\}=\bigl\{\bigl(0,\ldots,0,\\
&\ha y^2-\ts\frac{i}{\al}\sum_{j=1}^{m-1}\psi_j(y)-\ts\frac{i}{\al}
\arg(y+iP(y)^{-1/2})\bigr):y\in\R\bigr\},
\end{align*}
which bounds a noncompact $J$-holomorphic curve $\Si$ in the
$z_m$-axis as shown.

We will try and describe a type II singularity of Lagrangian MCF
$\{L^t:t\in[0,T)\}$ with a singularity at $x\in M$ modelled on these
LMCF translators $L$, using Principle \ref{bs3princ1}(b).
Identifying $M$ with $T_xM\cong\C^m$ near $x\in M$, each $L^t$
should to `first order' approximate an LMCF translator $L$ from
Example \ref{bs2ex6}, and as $t\ra T$ these LMCF translators should
slowly shrink homothetically, as well as translate. What interests
us is the `second order' changes to $L$ which cause this shrinking.

Far to the right in Figure \ref{bs3fig9}, the LMCF translator $L$
approximates two non-intersecting affine Lagrangian planes $\ti\Pi_0,\ti\Pi_{\bs\phi}$ in $\C^m$ from \eq{bs2eq11}, just as far to the right in Figure \ref{bs2fig1}, the `grim reaper' approximates two non-intersecting parallel lines in $\C$. I suggest that to `second order' in $L^t$, the two planes $\ti\Pi_0,\ti\Pi_{\bs\phi}$ should be bent towards each other by a small angle, introducing a new immersed self-intersection point, and so that the noncompact $J$-holomorphic curve $\Si$ becomes a compact `teardrop'
as in Figure \ref{bs2fig3}, which makes $HF^*$ obstructed. This modification $\ti L$ of $L$ is sketched in Figure~\ref{bs3fig10}.

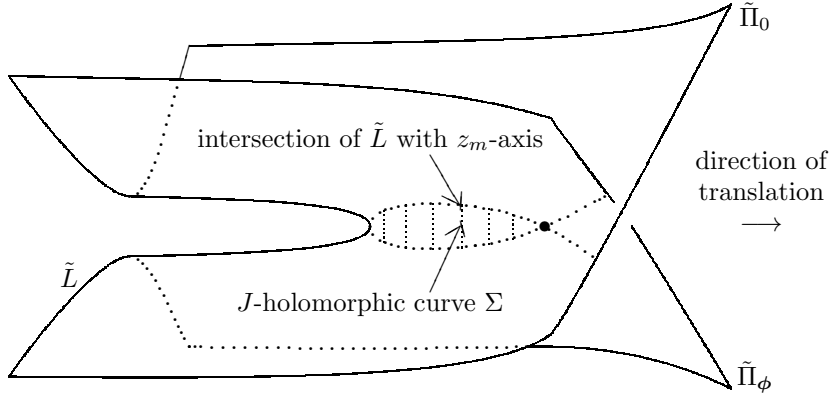
\begin{figure}[htb]
\centerline{$\splinetolerance{.8pt}
\begin{xy}
0;<.8mm,0mm>:
,(-30,30);(60,37)**\crv{(10,31)&(50,30)}
,(-60,25);(30,18)**\crv{(-10,24)&(10,24)}
,(-60,25);(-40,5)**\crv{(-53,15)&(-44,5)}
,(-40,5);(-31.6,24.4)**\crv{~*=<3pt>{.} (-36,5)&(-34,17.5)}
,(-30,30);(-31.6,24.4)**\crv{}
,(-40,5);(0,0)**\crv{(-20,4.5)&(0,4)}
,(-40,-5);(0,0)**\crv{(-20,-4.5)&(0,-4)}
,(-60,-25);(-40,-5)**\crv{(-53,-15)&(-44,-5)}%
,(-40,-5);(-30,-20)**\crv{~*=<3pt>{.} (-36,-5)&(-32,-17.5)}
,(-60,-25);(30,-18)**\crv{(-10,-25.5)&(20,-25)}
,(30,-18);(60,37)**\crv{(33,-15)}
,(30,18);(40.5,4)**\crv{(35,11)}
,(43.5,0);(60,-27)**\crv{(52,-13.5)}
,(-30,-20);(26,-20)**\crv{~*=<3pt>{.} (0,-20.5)}
,(26,-20);(60,-27)**\crv{(45,-19.5)}
,(39,5);(0,0)**\crv{~*=<3pt>{.} (30,0)&(20,-4.5)&(0,-4)}
?(.85)="b"
?(.75)="c"
?(.65)="d"
?(.55)="e"
?(.45)="f"
?(.35)="g"
,(37,-5);(0,0)**\crv{~*=<3pt>{.} (30,0)&(20,4.5)&(0,4)}
?(.85)="l"
?(.75)="m"
?(.65)="n"
?(.55)="o"
?(.45)="p"
?(.35)="q"
,"b";"l"**@{.}
,"c";"m"**@{.}
,"d";"n"**@{.}
,"e";"o"**@{.}
,"f";"p"**@{.}
,"g";"q"**@{.}
,(29,-.2)*{\bu}
,(65,0)*{\longra}
,(65,11)*{\text{direction of}}
,(65,6)*{\text{translation}}
,(15,3.5);(10,12)**\crv{}
,(15,3.5);(12,5)**\crv{}
,(15,3.5);(15.2,6.5)**\crv{}
,(0,15)*{\text{intersection of $\ti L$ with $z_m$-axis}}
,(15,0);(10,-11)**\crv{}
,(15,0);(12.5,-1.5)**\crv{}
,(15,0);(15.7,-2.5)**\crv{}
,(0,-13)*{\text{$J$-holomorphic curve $\Si$}}
,(64,35)*{\ti\Pi_0}
,(64,-25)*{\ti\Pi_{\bs\phi}}
,(-50,-8)*{\ti L}
\end{xy}$}
\caption{Modification $\ti L$ of Joyce--Lee--Tsui LMCF translator}
\label{bs3fig10}
\end{figure}

I expect that this `bending' of $\ti\Pi_0,\ti\Pi_{\bs\phi}$ towards one another is both the `outside influence' in Principle \ref{bs3princ1}(b) which makes $L$ shrink and causes the finite time singularity, and also
the cause of the self-intersection point, the `teardrop' curve
$\Si$, and the obstructions to~$HF^*$.
\label{bs3ex7}
\end{ex}

\begin{conj} In dimension $m\ge 2,$ Example\/ {\rm\ref{bs3ex7}}
describes a possible finite time singularity of graded, immersed
Lagrangian MCF with $HF^*$ obstructed, after which one cannot
continue the flow in graded Lagrangian MCF.

Such finite time singularities admit type II blow ups, as in
Theorem\/ {\rm\ref{bs2thm5},} which are Lagrangian MCF translators
from Example\/~{\rm\ref{bs2ex6}}.

This is a \begin{bfseries}generic singularity\end{bfseries} of
Lagrangian MCF, that is, if Lagrangian MCF starting from $L^0$
develops such a singularity, then so does Lagrangian MCF starting
from any sufficiently small Hamiltonian perturbation $\ti L^0$
of\/~$L^0$.

All this is possible for Lagrangians with arbitrarily small phase
variation.
\label{bs3conj5}
\end{conj}

\subsection{A Thomas--Yau type conjecture}
\label{bs39}

Finally we state our second main conjecture, about the programme of \S\ref{bs32}, which summarizes the discussion of \S\ref{bs32}--\S\ref{bs37}. We call it a `Thomas--Yau type conjecture', as it aims to update the conjectures of Thomas and Yau~\cite{Thom,ThYa}. 

Our focus here is mostly on the unique long-time existence of immersed Lagrangian MCF with surgeries, although proving the conjecture would go some way to proving Conjecture \ref{bs3conj1} on Bridgeland stability conditions. To simplify the possible finite time singularities, we take $L$ generic in its Hamiltonian isotopy class. To minimize the singular Lagrangians to be included in $D^b\sF(M)$, we do not require $(L^{T_i},E^{T_i},b^{T_i})$ or $\lim_{T\ra\iy}(L^t,E^t,b^t)$ to be objects in~$D^b\sF(M)$.

\begin{conj} Let\/ $(M,J,g,\Om)$ be a Calabi--Yau $m$-fold, either compact or suitably convex at infinity, and\/  $D^b\sF(M)$ an enlarged version of the derived Fukaya category of Lagrangian branes in\/ $M$ from {\rm\cite{FOOO},} including classes of immersed or singular Lagrangians, depending on the dimension\/~{\rm$m$:}
\begin{itemize}
\setlength{\parsep}{0pt}
\setlength{\itemsep}{0pt}
\item[{\bf(i)}] When $m=1,$ $D^b\sF(M)$ can be the usual derived Fukaya category of nonsingular, embedded Lagrangian branes.
\item[{\bf(ii)}] When $m\ge 2,$ $D^b\sF(M)$ must include \begin{bfseries}immersed\end{bfseries} Lagrangians, as in Akaho and Joyce\/ {\rm\cite{AkJo}} and\/ {\rm\S\ref{bs26}}. For $m=2,$ these are all of\/~$D^b\sF(M)$.
\item[{\bf(iii)}] When $m\ge 3,$ $D^b\sF(M)$ must also include singular Lagrangians with \begin{bfseries}stable special Lagrangian singularities\end{bfseries}, as in\/ {\rm\S\ref{bs36}}. When $m=3,$ these include Lagrangians with \begin{bfseries}isolated conical singularities\end{bfseries} in the sense of\/ {\rm\cite{Joyc1,Joyc2,Joyc3,Joyc4,Joyc5}} modelled on the special Lagrangian $T^2$-cone from {\rm\eq{bs2eq4},} and this may be the only kind of stable singularity when $m=3$. When $m\ge 4,$ stable singularities may be more complicated, and need not be isolated.
\end{itemize}

Let\/ $(L,E)$ be a Lagrangian brane in $M$ with\/ $HF^*$ unobstructed, and suppose $L$ is generic in its Hamiltonian isotopy class. Let\/ $b$ be a bounding cochain for $(L,E)$. Then there is a unique family $\bigl\{(L^t,E^t,b^t):t\in[0,\iy)\bigr\}$ satisfying:
\begin{itemize}
\setlength{\parsep}{0pt}
\setlength{\itemsep}{0pt}
\item[{\bf(a)}] $(L^0,E^0,b^0)=(L,E,b)$.
\item[{\bf(b)}] There is a finite series of \begin{bfseries}singular times\end{bfseries} $0<T_1<T_2<\cdots<T_N$
such that if\/ $t\in [0,\iy)\sm\{T_1,\ldots,T_N\}$ then
$(L^t,E^t,b^t)$ is an object in $D^b\sF(M)$ isomorphic to $(L,E,b),$
with\/ $L^t$ a (possibly immersed or singular) compact, graded
Lagrangian in $(M,\om),$ with\/ $HF^*$ unobstructed.
\item[{\bf(c)}] The family $\bigl\{L^t:t\in[0,\iy)\sm\{T_1,\ldots,T_N\}\bigr\}$ satisfies Lagrangian mean curvature flow, and\/ $\bigl\{E^t:t\in[0,\iy)\sm\{T_1,\ldots,T_N\}\bigr\}$ is locally constant in $t$. The bounding cochains $b^t$ also change by a kind of `parallel transport' for $t\in[0,\iy)\sm\{T_1,T_2,\ldots\}$ as in\/ {\rm\S\ref{bs25}--\S\ref{bs26},} to ensure that the isomorphism class of\/ $(L^t,E^t,b^t)$ in $D^b\sF(M)$ remains constant.
\item[{\bf(d)}] At each singular time $T_1,\ldots,T_N,$ the flow undergoes a surgery, which may involve a finite time singularity of Lagrangian MCF, and a change in the topology of\/ $L^t$. The kinds of surgery allowed include `opening a neck' as in\/ {\rm\S\ref{bs34}} when $m\ge 1,$  `neck pinches' as in\/ {\rm\S\ref{bs35}} when $m\ge 2,$ transitions to and from Lagrangians with `stable special Lagrangian singularities' as in\/ {\rm\S\ref{bs36}} when $m\ge 3,$ and `collapsing zero objects' as in\/ {\rm\S\ref{bs37}} for $m\ge 1$ (the latter is excluded for almost calibrated Lagrangians). 

We do not require $(L^{T_i},E^{T_i},b^{T_i})$ to be an object in
$D^b\sF(M),$ as the singularities of\/ $L^{T_i}$ may be too bad,
and if so, $b^{T_i}$ is meaningless.

\item[{\bf(e)}] The family $\{L^t:t\in[0,\iy)\}$ is continuous as graded Lagrangian integral currents in $M$ in Geometric Measure Theory.

In graded Lagrangian integral currents, we have $\lim_{t\ra\iy}L^t=L_1+\cdots+L_n$ for some $n\ge 0,$ where $0\ne L_j$ for $j=1,\ldots,n$ is a nonzero, compactly-supported, graded, special Lagrangian integral current with phase $e^{i\pi\phi_j}$ and grading $\th_{L_j}=\pi\phi_j,$ with\/ $\phi_1>\cdots>\phi_n$.

For the Bridgeland stability condition $(Z,\cP)$ on $D^b\sF(M)$ discussed in Conjecture\/ {\rm\ref{bs3conj1},} if\/ $n=1$ then\/ $(L,E,b)\in\cP(\phi_1),$ and otherwise $(L,E,b)\notin\cP(\phi)$ for any $\phi\in\R$.
\end{itemize}
\label{bs3conj6}
\end{conj}

\begin{rem}{\bf(i)} The most feasible case of the conjecture is that of Lagrangian MCF in dimension $m=2$, starting from an almost calibrated Lagrangian $L$ generic in its Hamiltonian isotopy class.
\smallskip

\noindent{\bf(ii)} Assuming the initial object $(L,E,b)$ is semistable or stable in the sense of Conjectures \ref{bs3conj1} and \ref{bs3conj2} means in part (e) that the limit $\lim_{t\ra\iy}L^t=L_1$ is only one (singular) special Lagrangian, rather than a finite union $L_1\cup\cdots\cup L_n$ of special Lagrangians with different phases, but otherwise it does not simplify things: we still expect nontrivial finite time singularities, and surgeries.
\smallskip

\noindent{\bf(iii)} It is an interesting question whether there are useful extra assumptions on $L$ which limit the kinds of singularities occurring at the singular times $T_1,T_2,\ldots.$ For example, if $L$ is generic in its Hamiltonian isotopy class then only singularities of `index zero' appear, as in Remark \ref{bs3rem3}(iv), and if $L$ is almost calibrated, then as in \S\ref{bs37} `collapsing zero objects' cannot happen. Thomas and Yau give conditions \cite[(7.1) or (7.2)]{ThYa} preventing `neck pinches' in \S\ref{bs35} dividing $L$ into two pieces $L_1\amalg L_2$ from happening, although I expect other singularities can.

There are some very special situations in which Lagrangian MCF is known to exist for all time without singularities, such as the Lagrangian $T^m$-graphs in $T^{2m}$ studied by Smoczyk and Wang \cite{SmWa}, or Lagrangian MCF starting from a small perturbation of a smooth special Lagrangian. But apart from these, I do not know of any useful, nontrivial conditions on $L$ under which I expect the flow to exist for all time without singularities, as hoped for in~\cite[Conj.~7.3]{ThYa}.

\label{bs3rem9}
\end{rem}

\section{Generalizations}
\label{bs4}

Finally we discuss two variations on the conjectural picture of \S\ref{bs3}. The first in \S\ref{bs41} is a simplification: for (immersed) exact Lagrangians in exact Calabi--Yau $m$-folds, we can use a (not yet written down) version of Seidel's Lagrangian Floer theory \cite{Seid2} in place of that of Fukaya et al.\ \cite{FOOO}, and in particular, we can dispense with Novikov rings $\La_\nov$. The second, in \S\ref{bs42}, considers whether we can generalize our picture to Lagrangian MCF in K\"ahler--Einstein manifolds rather than Calabi--Yau $m$-folds.

\subsection{Exact Lagrangians in exact Calabi--Yau $m$-folds}
\label{bs41}

Consider the situation of \S\ref{bs3} in the special case in which the 
Calabi--Yau $m$-fold $(M,J,g,\Om)$ is noncompact, and $(M,\om)$ is a Liouville manifold in the sense of \cite{Seid2}, so that the K\"ahler form $\om=\d\la$ for some Liouville form $\la$. Then we can restrict our attention to {\it exact\/} graded Lagrangians $L$ in $M$, and for each such $L$ we can choose a {\it potential\/} $f_L:L\ra\R$ with $\d f_L=\la\vert_L$. 

Then Seidel \cite{Seid2} defines Lagrangian Floer cohomology $HF^*(L_1,L_2)$ and Fukaya categories $D^b(\sF(M))_{\rm ex},D^\pi(\sF(M))_{\rm ex}$ for embedded, exact, graded Lagrangians $L_1,L_2$ in $M$, which are simpler and more complete than Fukaya, Oh, Ohta and Ono \cite{Fuka1,Fuka2,FOOO}. Here are two differences between the two theories:
\begin{itemize}
\setlength{\parsep}{0pt}
\setlength{\itemsep}{0pt}
\item[(i)] Seidel does not work over the Novikov ring $\La_\nov$. This is because $\La_\nov$ in \cite{FOOO} keeps track of infinite sums $\sum_{i=0}^\iy n_iP^{A_i}$, where $n_i$ is the number of $J$-holomorphic discs of area $A_i$ and $A_i\ra\iy$ as $i\ra\iy$ in a problem, but in the exact case, the area $A$ for $J$-holomorphic discs is uniquely determined.
\item[(ii)] Seidel does not include bounding cochains in his theory. This is because bounding cochains $b$ for $L$ are there in \cite{FOOO} to compensate for $J$-holomorphic discs $\Si$ with boundary in $L$, but for $L$ exact, there are no such non-constant $\Si$. In effect, Seidel sets $b=0$ throughout. 
\end{itemize}

I expect that for exact Lagrangian branes $(L,E)$ in Calabi--Yau Liouville manifolds, one should be able to carry out a simplified version of the programme of \S\ref{bs3}, using an extended version of Seidel's theory \cite{Seid2} in place of an extended version of Fukaya et al.\ \cite{Fuka1,Fuka2,FOOO}. To do this, one must extend Seidel's derived Fukaya category $D^\pi(\sF(M))_{\rm ex}$ of exact Lagrangians to include immersed Lagrangians, as in \S\ref{bs26}, and some classes of singular Lagrangians, as in~\S\ref{bs36}.

We need not work over $\La_\nov$ as in (i), instead one should take just coefficients in the field $\F$ of Definition \ref{bs2def5}. However, including immersed Lagrangians will mean that one has to consider obstructions to $HF^*$ and bounding cochains, as in (ii). The next definition explains how to do this.

\begin{dfn} Let $(M,J,g,\Om)$ be a Calabi--Yau Liouville manifold, and $(L,E)$ a compact, graded, immersed Lagrangian brane in $M$ with only transverse self-intersection points, with potential $f_L:L\ra\R$. Here $E\ra L$ is a rank one $\F$-local system as in~\S\ref{bs25}. 

Define a {\it bounding cochain\/} $b=(b_p)$ for $L$ to assign to each point $p\in M$ at which two local sheets $L^p_+,L^p_-$ of $L$ intersect transversely with $\mu_{L^p_+,L^p_-}(p)=1$ and $f_{L^p_+}(p)\ge f_{L^p_-}(p)$, an element $b_p\in\Hom_\F\bigl(E_+\vert_p,E_-\vert_p\bigr)$, as in \eq{bs2eq17} but without the `$\ot_\F\La_\nov^{\ge 0}$'. This data $(b_p)$ must satisfy the following condition: for each point $q\in M$ at which two local sheets $L^q_+,L^q_-$ of $L$ intersect transversely with $\mu_{L^q_+,L^q_-}(q)=2$, we should have
\e
\sum_{\begin{subarray}{l} \text{$k\ge 0$, $p_1,\ldots,p_k\in M$: local sheets $L^{p_i}_+,L^{p_i}_-$ of $L$ intersect}\\
\text{transversely at $p_i$, $\mu_{L^{p_i}_+,L^{p_i}_-}(p_i)=1$, $f_{L^{p_i}_+}(p_i)\ge f_{L^{p_i}_-}(p_i),$ $i=1,\ldots,k$}\end{subarray} 
\!\!\!\!\!\!\!\!\!\!\!\!\!\!\!\!\!\!\!\!\!\!\!\!\!\!\!
\!\!\!\!\!\!\!\!\!\!\!\!\!\!\!\!\!\!\!\!\!\!\!\!\!\!\!
\!\!\!\!\!\!\!\!\!\!\!\!\!\!\!\!\!\!\!\!\!\!\!\!\!\!\!
\!\!\!\!\!\!\!\!\!\!\!\!\!\!\!\!\!\!\!\!\!\!\!\!\!\!\!
\!\!\!\!\!\!\!\!\!\!\!\!\!\!\!\!\!\!\!\!\!\!\!\!\!\!\!
\!\!\!\!\!\!\!\!\!\!\!\!\!\!\!} N_{q,p_1,\ldots,p_k}\cdot (b_{p_1}\ot\cdots\ot b_{p_k})= 0 \quad\text{in $\Hom_\F\bigl(E_+\vert_q,E_-\vert_q\bigr)$.}
\label{bs4eq1}
\e

\begin{figure}[htb]
\centerline{$\splinetolerance{.8pt}
\begin{xy}
0;<1mm,0mm>:
,(20,0);(0,12)**\crv{(10,4)}
?(.1)="a"
?(.3)="b"
?(.5)="c"
?(.7)="d"
?(.9)="e"
?(.6)="xx"
,(-20,0);(0,12)**\crv{(-10,4)}
?(.9)="f"
?(.7)="g"
?(.5)="h"
?(.3)="i"
?(.1)="j"
?(.6)="zz"
,(0,12);(4,17)**\crv{(2.5,14)}
,(0,12);(-4,17)**\crv{(-2.5,14)}
,(-20,0);(-30,-6)**\crv{(-30,-5)}
,(20,0);(30,-6)**\crv{(30,-5)}
,(20,0);(0,-12)**\crv{(10,-4)}
?(.1)="k"
?(.3)="l"
?(.5)="m"
?(.7)="n"
?(.9)="o"
?(.6)="x"
,(-20,0);(0,-12)**\crv{(-10,-4)}
?(.9)="p"
?(.7)="q"
?(.5)="r"
?(.3)="s"
?(.1)="t"
?(.6)="z"
,(0,-12);(4,-17)**\crv{(2.5,-14)}
,(0,-12);(-4,-17)**\crv{(-2.5,-14)}
,(20,0);(30,3)**\crv{(25,2)}
,(-20,0);(-30,3)**\crv{(-25,2)}
,"a";"k"**@{.}
,"b";"l"**@{.}
,"c";"m"**@{.}
,"d";"n"**@{.}
,"e";"o"**@{.}
,"f";"p"**@{.}
,"g";"q"**@{.}
,"h";"r"**@{.}
,"i";"s"**@{.}
,"j";"t"**@{.}
,(0,-12)*{\bu}
,(0,12)*{\bu},
(-20,0)*{\bu}
,(-20,-3)*{p_2}
,(20,0)*{\bu}
,(20,-3)*{q}
,(9,-12.3)*{p_3=p_k}
,(4,12.3)*{p_1}
,(0,0)*{\Si}
,(-32,3.5)*{L^{p_2}_+}
,(-32,-5.5)*{L^{p_2}_-}
,(33.4,3.5)*{L^q_+}
,(32.5,-5.5)*{L^q_-}
,(6.5,-16)*{L^{p_3}_-}
,(-6.2,-16)*{L^{p_3}_+}
,(6.5,16)*{L^{p_1}_+}
,(-6.2,16)*{L^{p_1}_-}
\end{xy}$}
\caption{Holomorphic $(k+1)$-gon $\Si$ with boundary in $L$, case $k=3$}
\label{bs4fig1}
\end{figure}
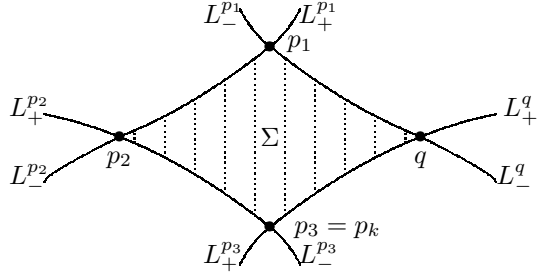

Here the term
\e
\begin{split}
N_{q,p_1,\ldots,p_k}\in \Hom(E_+\vert_q,E_+\vert_{p_k})\ot_\F \Hom(E_-\vert_{p_k},E_+\vert_{p_{k-1}})\ot_\F&\\
\cdots\ot_\F \Hom(E_-\vert_{p_2},E_+\vert_{p_1})\ot_\F \Hom(E_-\vert_{p_1},E_-\vert_q)&
\end{split}
\label{bs4eq2}
\e
in \eq{bs4eq1} `counts' $J$-holomorphic discs $\Si$ of the form shown in Figure \ref{bs4fig1}, in the case $k=3$, weighted by the parallel transport in $E$ around $\pd\Si$ in the clockwise direction in Figure \ref{bs4fig1}, where the parallel transport from $q$ to $p_k$ contributes the term $\Hom(E_+\vert_q,E_+\vert_{p_k})$ in \eq{bs4eq2}, and so on. Contracting $N_{q,p_1,\ldots,p_k}$ in \eq{bs4eq2} with $b_{p_1}\in \Hom_\F\bigl(E_+\vert_{p_1},E_-\vert_{p_1}\bigr),\ldots,b_{p_k}\in \Hom_\F\bigl(E_+\vert_{p_k},E_-\vert_{p_k}\bigr)$ yields an element of $\Hom_\F\bigl(E_+\vert_q,E_-\vert_q\bigr)$, as desired.

The area of $\Si$ in Figure \ref{bs4fig1} is
\begin{equation*}
\area(\Si)=f_{L_+^q}(q)-f_{L_-^q}(q)-\ts\sum_{i=1}^k\bigl(f_{L_+^{p_i}}(p_i)-f_{L_-^{p_i}}(p_i)\bigr)\le f_{L_+^q}(q)-f_{L_-^q}(q).
\end{equation*}
Since these areas are uniformly bounded, standard results on $J$-holomorphic curves tell us that the family of such $\Si$ for all $k,p_1,\ldots,p_k$ is compact, and therefore the sum \eq{bs4eq1} has only finitely many nonzero terms, and is well defined.

If $L$ is embedded, then there are no possibilities for $p$ or $q$ above, so $b=\es$ is trivially the unique bounding cochain.

We say that $(L,E)$ has $HF^*$ {\it unobstructed\/} if a bounding cochain $b$ exists, and has $HF^*$ {\it obstructed\/} otherwise. To work with Lagrangians $L$ rather than Lagrangian branes $(L,E)$, we take $E$ to be the trivial $\F$-local system $\F\t L\ra L$.
\label{bs4def1}
\end{dfn}

The next conjecture, similar to Akaho and Joyce \cite{AkJo}, should not be difficult. 

\begin{conj} Seidel's theory of Lagrangian Floer cohomology and Fukaya categories\/ {\rm\cite{Seid2}} may be extended to include immersed Lagrangians\/ $L$ and Lagrangian branes\/ $(L,E)$ with transverse self-intersections with choices of bounding cochains, in the sense of Definition\/~{\rm\ref{bs4def1}}. 
\label{bs4conj1}
\end{conj}

The author expects that the programme of \S\ref{bs3} can be carried out for exact Lagrangians in this extended Seidel-style derived Fukaya category. Sections \ref{bs34} and \ref{bs35} should be modified as follows. Suppose as in \S\ref{bs34} that $\{(L^t,E^t):t\in [0,T)\}$ is a family of {\it exact\/} Lagrangian branes satisfying Lagrangian MCF, with transverse self-intersections, potential functions $f_{L^t}$, and bounding cochains $b^t$ for $t\in[0,T)$ as in Definition \ref{bs4def1}.

Write $p^t$ for an intersection of local sheets $L^t_+,L^t_-$ of $L^t$ for $t\in[0,T)$, depending smoothly on $t$, with $\mu_{L^t_+,L^t_-}(p^t)=1$. Then $b^t$ includes an element $b^t_{p^t}\in\Hom_\F\bigl(E^t_+\vert_{p^t},E^t_-\vert_{p^t}\bigr)$ if $f_{L^t_+}(p^t)\ge f_{L^t_-}(p^t)$. Then `opening a neck' at $p^t$ as in \S\ref{bs34} should happen at time $t=T_1$ if $b^t_{p^t}\ne 0$ and we have $f_{L^t_+}(p^t)>f_{L^t_-}(p^t)$ for $T_1-\ep<t<T_1$ and $f_{L^{T_1}_+}(p^{T_1})=f_{L^{T_1}_-}(p^{T_1})$.

Similarly, for the `neck pinch' in \S\ref{bs35}, the analogue of Conjecture \ref{bs3conj4} should say that the immersed Lagrangian $L^T$ at the singular time $t=T$ should have potential $f_{L^T}$ with $f_{L^T_+}(p)=f_{L^T_-}(p)$, and $f_{L^t_+}(p^t)>f_{L^t_+}(p^t)$ for $T<t<T+\ep$, and $b^t(p^t)=a_0$ for $T\le t<T+\ep$, where $a_0$ is as in Conjecture~\ref{bs3conj4}(vii).

\subsection{`Balanced' Lagrangians in K\"ahler--Einstein manifolds}
\label{bs42}

Let $(M,J,g)$ be a K\"ahler manifold. Then mean curvature flow of compact submanifolds preserves Lagrangian submanifolds not only if $g$ is Ricci-flat (the Calabi--Yau case), but also if $g$ is Einstein, that is, if $(M,J,g)$ is a {\it K\"ahler--Einstein manifold}. Thus, it is natural to ask whether, and to what extent, the programme of \S\ref{bs3} can be generalized to K\"ahler--Einstein manifolds.

To do this, we should answer the following questions:

\begin{quest}{\bf(a)} Compact, graded Lagrangians in Calabi--Yau $m$-folds have the property that Lagrangian MCF stays in the same Hamiltonian isotopy class. Is there an interesting class of Lagrangians in K\"ahler--Einstein manifolds with the same property?
\smallskip

\noindent{\bf(b)} For a Calabi--Yau $m$-fold\/ $M,$ one can define derived Fukaya categories $D^b\sF(M),D^\pi\sF(M)$ of graded Lagrangians, which are $\Z$-graded triangulated categories. What is the appropriate analogue in the K\"ahler--Einstein case?
\smallskip

\noindent{\bf(c)} As in {\rm\S\ref{bs31},} for Calabi--Yaus there conjecturally exist Bridgeland stability conditions on $D^b\sF(M)\simeq D^\pi\sF(M),$ such that special Lagrangians are semistable. In the K\"ahler--Einstein case, is there some useful notion of `stability condition' on $D^b\sF(M)$ or $D^\pi\sF(M),$ such that minimal Lagrangians are semistable?
\smallskip

\noindent{\bf(d)} In the Calabi--Yau case, as {\rm\cite[Th.~4.3]{ThYa}} special Lagrangians are unique in their isomorphism classes in $D^b\sF(M)$. In the K\"ahler--Einstein case, are minimal Lagrangians unique in their isomorphism classes in $D^b\sF(M)$?
\smallskip

\noindent{\bf(e)} As in {\rm\S\ref{bs39},} for Calabi--Yaus we conjecture unique long-time existence of Lagrangian MCF with surgeries starting from a Lagrangian brane $(L,E)$ with $HF^*$ unobstructed. Should we expect an analogue in the K\"ahler--Einstein case?

\label{bs4quest1}
\end{quest}

Note that the Einstein condition on $(M,J,g)$ is that the Ricci curvature $R_{ab}$ of $g$ satisfies $R_{ab}=\la g_{ab}$ for some $\la\in\R$, where $\la=0$ is the Calabi--Yau case. By rescaling $g$ we can take $\la=1,0$ or $-1$. The cases $\la>0$ (so $g$ has positive scalar curvature) and $\la<0$ (so $g$ has negative scalar curvature) are likely to behave differently, and should be considered separately. When $\la>0$, $[\om]$ is a positive multiple of $c_1(M)$ in $H^2(M,\R)$, and so the symplectic manifold $(M,\om)$ is called {\it monotone}. Monotone symplectic manifolds have been extensively studied.

Here are some partial answers to Question \ref{bs4quest1}:
\begin{itemize}
\setlength{\parsep}{0pt}
\setlength{\itemsep}{0pt}
\item[(a)] The appropriate class of Lagrangians are known as `balanced' or `Bohr--Sommerfeld monotone' Lagrangians, as in Seidel \cite[\S 6]{Seid4}, for instance. For a balanced Lagrangian $L$, the mean curvature $H$ corresponds to an exact 1-form on $L$, so LMCF stays in the same Hamiltonian isotopy class.

When $\la>0$, balanced Lagrangians are sometimes called {\it monotone Lagrangians}, although definitions of these vary (balanced implies monotone, but the definition of monotone Lagrangians involving $\pi_2(M,L)$ rather than $H_2(M,L;\R)$ does not imply balanced).
\item[(b)] We would like to define a Fukaya category $D^b\sF(M)$ of balanced Lagrangians in a K\"ahler--Einstein manifold $(M,J,g)$. When $\la<0$ (negative scalar curvature), a model for this is provided by Seidel \cite{Seid4}, who defines $D^b\sF(M)$ when $M$ is the genus two curve.

When $\la>0$ (positive scalar curvature, monotone), a Fukaya category $D^b\sF(M)$ is defined by Sheridan \cite[\S 3]{Sher2}, following Oh \cite{Oh} and Seidel \cite{Seid3}. In this case one can include an extra parameter $w\in\F$, where $\F$ is the base field as in \S\ref{bs25}, to get categories $D^b\sF(M)^w$ for $w\in\F$ generated by Lagrangian branes $(L,E)$, such that for a generic $p\in L$ we have
\begin{equation*}
\sum_{\begin{subarray}{l}\text{Maslov 2 $J$-holomorphic discs $\Si$ in $M$ with $p\in\pd\Si\subset L$}\end{subarray}\!\!\!\!\!\!\!\!\!\!\!\!\!\!\!\!\!\!\!\!\!} \mathop{\rm sign}(\Si)\Hol_{\pd\Si}(E)=w,\end{equation*}
where $\Hol_{\pd\Si}(E)\in\F$ is the holonomy of the $\F$-local system $E\ra L$ around $\pd\Si\subset L$. One expects $D^b\sF(M)^w=0$ for all but finitely many $w\in\F$. Lagrangians in $D^b\sF(M)^w$ for $w\ne 0$ are called {\it weakly unobstructed}.

In both cases $D^b\sF(M)$ will not be a $\Z$-graded triangulated category: for oriented Lagrangians it will be $\Z_2$-graded, though as in Seidel \cite{Seid1} we can improve this to $\Z_{2k}$-graded if $k$ divides~$c_1(M)$.
\item[(c)] As $D^b\sF(M)$ is not $\Z$-graded, Bridgeland stability conditions on $D^b\sF(M)$ do not make sense. Also the central charge map $Z$ in \eq{bs3eq2}, and phases of Lagrangians, make no sense in the K\"ahler--Einstein case.

Nonetheless, it seems possible that some features of Bridgeland stability may survive to the K\"ahler--Einstein case, in particular, being given a set of `semistable objects' (represented by minimal Lagrangians) satisfying some axioms, such that every object in $D^b\sF(M)$ has a unique decomposition into semistable objects of the form \eq{bs3eq1}. This may be worth further study.

\item[(d)] When $\la>0$ (positive scalar curvature), minimal Lagrangians need not be unique in their isomorphism classes in $D^b\sF(M)$. For example, a minimal Lagrangian $L$ in $\CP^m$ (such as the Clifford torus) is isomorphic to all of its images under the Lie group $\mathop{\rm Aut}(\CP^m)=\mathop{\rm PU}(m+1)$.

When $\la<0$ (negative scalar curvature), the author expects minimal Lagrangians to be unique in their isomorphism classes in $D^b\sF(M)$.

This is parallel to uniqueness for K\"ahler--Einstein metrics: if $(M,J)$ is a compact complex manifold, then K\"ahler--Einstein metrics on $M$ are unique in their K\"ahler classes if $\la=0$ or $\la<0$, but need not be unique if~$\la>0$.
\item[(e)] The author's guess is that the answer is yes: or at least, in any dimension $m$ in which Conjecture \ref{bs3conj6} holds for Calabi--Yaus, some analogous conjecture on long-time existence of LMCF with surgeries should also hold in dimension $m$ for the K\"ahler--Einstein case.

Some justification is that if the conjecture is false, it will be because of singular behaviour of LMCF developing locally near a point $x\in M$. But locally, LMCF in Calabi--Yau $m$-folds and in K\"ahler--Einstein manifolds looks essentially the same, the differences are global.
\end{itemize}

\begin{ex}{\bf(i)} Let $M$ be $\CP^1$, a K\"ahler--Einstein manifold with $\la>0$, and $L$ an embedded $\cS^1$ in $\CP^1$, which is automatically Lagrangian. Then $\CP^1\sm L$ is two open discs $D_1\amalg D_2$, and $L$ is balanced if and only if $\area(D_1)=\area(D_2)$.
 
One can show using known results on the curve-shortening flow \cite{AbLa,Ange1,Ange2,Gray} that Lagrangian MCF in $\CP^1$ starting from $L^0=L$ exists for all time and converges as $t\ra\iy$ to a great circle (a minimal Lagrangian in $\CP^1$) if $L$ is balanced, and collapses to a point in finite time if $L$ is unbalanced.
\smallskip

\noindent{\bf(ii)} Let $M$ be a Riemann surface of genus $g>1$ with a hyperbolic metric, a K\"ahler--Einstein manifold with $\la>0$, and $L$ an embedded $\cS^1$ in $M$ with $[L]\ne 0$ in $H_1(M,\Z)$. Then Lagrangian MCF in $M$ starting from $L^0=L$ exists for all time and converges as $t\ra\iy$ to the unique closed geodesic $\ga$ in $M$ with homology class $[L]$. Here $\ga$ is balanced, and the flow stays in a fixed Hamiltonian isotopy class (the Hamiltonian isotopy class of $\ga$) if and only if $L$ is balanced.
\smallskip

This suggests that in (i),(ii), all (indecomposable?) objects in $D^b\sF(M)$ are `semistable', so the notion of stability in (c) above is not interesting for~$m=1$.
\label{bs4ex1}
\end{ex}

\begin{ex} Let $M$ be $\CP^2$, a K\"ahler--Einstein manifold with $\la>0$, and take $\F=\C$. Mirror Symmetry \cite{AKO} predicts that the mirror of $M$ is the Landau--Ginzburg model $f:U=(\C^*)^2\ra\C$, $f(x,y)=x+y+1/xy$. Thus for $w\in\C$ we expect an equivalence of $\Z_2$-periodic triangulated categories $D^\pi\sF(\CP^2)^w\simeq\mathop{\rm MF}(f-w:U\ra\C)$, where $D^\pi\sF(\CP^2)^w$ is the idempotent-completed Fukaya category of oriented, balanced Lagrangian branes $(L,E)$ in $\CP^2$ as in (b) above, and $\mathop{\rm MF}(\cdots)$ is the matrix factorization category. Since $f$ is Morse, this indicates that $D^\pi\sF(\CP^2)^w\simeq D^b_{\Z_2}(\mathop{\rm Vect}_\C)$ if $w$ is one of the 3 critical values $3,3e^{2\pi i/3},3e^{-2\pi i/3}$ of $f$, and $D^\pi\sF(\CP^2)^w\simeq 0$ otherwise.

Few examples of oriented, balanced Lagrangians in $\CP^2$ are known, up to Hamiltonian isotopy. They include the minimal Lagrangian {\it Clifford torus}
\begin{equation*}
T^2_{\rm Cl}=\bigl\{[z_0,z_1,z_2]\in\CP^2:\md{z_0}=\md{z_1}=\md{z_2}\bigr\},
\end{equation*}
and two exotic examples, the Chekanov torus $T^2_{\rm Ch}$ \cite{ChSc}, and the Vianna torus $T^2_{\rm Vi}$ \cite{Vian}. All three are pairwise non Hamiltonian isotopic~\cite{ChSc,Vian}.

I expect the following picture. Each of $T^2_{\rm Cl},T^2_{\rm Ch},T^2_{\rm Vi}$ should have $\C$-local systems $E^w_{\rm Cl},E^w_{\rm Ch},E^w_{\rm Vi}$ for $w=3,3e^{2\pi i/3},3e^{-2\pi i/3}$, such that $(T^2_{\rm Cl},E^w_{\rm Cl}),(T^2_{\rm Ch},E^w_{\rm Ch}),\ab(T^2_{\rm Vi},E^w_{\rm Vi})$ are nontrivial objects in $D^\pi\sF(\CP^2)^w$ for each $w$. Furthermore, we should have $(T^2_{\rm Cl},E^w_{\rm Cl})\cong(T^2_{\rm Ch},E^w_{\rm Ch})\cong(T^2_{\rm Vi},E^w_{\rm Vi})$, because there is basically only one interesting object in~$D^\pi\sF(\CP^2)^w\simeq D^b_{\Z_2}(\mathop{\rm Vect}_\C)$.

As an aside, note that $\mathop{\rm End}^*(T^2_{\rm Cl},E^w_{\rm Cl})\cong \mathop{\rm End}^*(F\op F[1])$ as a $\Z_2$-graded ring, where $F$ is a simple object with $\mathop{\rm End}(F)=\C$. Thus $D^b\sF(\CP^2)^w$ is not idempotent complete, since $(T^2_{\rm Cl},E^w_{\rm Cl})$ is indecomposable in $D^b\sF(\CP^2)^w$, but decomposes as $F\op F[1]$ in $D^\pi\sF(\CP^2)^w$. This suggests that the analogue of Conjecture \ref{bs3conj3} may be false for~$\CP^2$.

Lagrangian MCF starting from $(T^2_{\rm Cl},E^w_{\rm Cl})$ is stationary, as $T^2_{\rm Cl}$ is minimal. I conjecture that Lagrangian MCF starting from $(T^2_{\rm Ch},E^w_{\rm Ch})$ or $(T^2_{\rm Vi},E^w_{\rm Vi})$ exists, with surgeries, for all time, as in (e) above, and converges as $t\ra\iy$ to $(T^2_{\rm Cl},E^w_{\rm Cl})$, or one of its images under $\mathop{\rm PU}(3)$. If this is so, then the exotic tori $T^2_{\rm Ch},T^2_{\rm Vi}$ and their local systems $E^w_{\rm Ch},E^w_{\rm Vi}$ can be obtained from $T^2_{\rm Cl}$ and $E^w_{\rm Cl}$ by applying one or more Lagrangian surgeries which change the Hamiltonian isotopy class of $T^2_{\rm Cl}$, but not the isomorphism class of $(T^2_{\rm Cl},E^w_{\rm Cl})$ in~$D^b\sF(\CP^2)^w$. 

My guess is that the right kind of surgery is that described in Remark \ref{bs3rem6}(d) in the special case $m=2$: at the singular time $t=T$ there should be a `neck pinch' as in \S\ref{bs35}, where $E^w_{\rm Ch},E^w_{\rm Vi}$ have nontrivial holonomy around the $\cS^1$ of the `neck', and for $t>T$ one immediately `opens the neck' as in \S\ref{bs34}, but in the opposite direction. 

If this is correct, then there should exist an immersed Lagrangian $L^T$ in $\CP^2$ diffeomorphic to $\cS^2$ (a Whitney sphere?), with one transverse self-intersection point $p$, such that resolving $L^T$ by connected sum at $p$ in one way gives $T^2_{\rm Cl}$ (up to Hamiltonian isotopy), and in the other way gives $T^2_{\rm Ch}$ or $T^2_{\rm Vi}$. Immersed Lagrangians with few double points are discussed by Ekholm et al.~\cite{EEMS}.
\label{bs4ex2}
\end{ex}

\medskip

\noindent{\small\sc The Mathematical Institute, Radcliffe
Observatory Quarter, Woodstock Road, Oxford, OX2 6GG, U.K.

\noindent E-mail: {\tt joyce@maths.ox.ac.uk.}}

\end{document}